\documentclass[leqno,11pt]{article}
\usepackage{latexsym}
\usepackage{amssymb}
\usepackage{amsfonts}
\usepackage{amsmath}
\usepackage{epsfig}
\usepackage{color}

\makeatletter
\long\def\unmarkedfootnote#1{{\long\def\@makefntext##1{##1}\footnotetext{#1}}}
\makeatother

\newtheorem{definition}{Definition}[section]

\newtheorem{lemma}[definition]{Lemma}

\newtheorem{theorem}[definition]{Theorem}

\newtheorem{proposition}[definition]{Proposition}

\newtheorem{corollary}[definition]{Corollary}

\newtheorem{remark}[definition]{Remark}

\def\a{\alpha}
\def\b{\beta}
\def\o{\Omega}

\def\u{{\bf u}}
\def\v{{\bf v}}
\def\f{\bf f}

\def\mo{|\Omega |}
\def\m2{|\Omega | /2}
\def\M2{\frac{|\Omega |}{2}}

\def\-p{\overline{p}}

\def\w0{{W_0^{1,p}(\Omega)}}

\def\R{\mathbb R}
\def\N{\mathbb N}

\def\rn{{{\R}^n}}
\def\rN{{{\R}^N}}

\def\MR{\mathcal R}

\newcommand{\hh}{{\cal H}^{n-1}}
\newcommand{\medint}{-\kern  -,395cm\int}
\newcommand{\medintinrigo}{-\kern  -,315cm\int}
\newcommand{\medelle}{-\kern  -,235cm L}
\newcommand{\medellenrigo}{-\kern  -,180cm L}
\newcommand{\qed}{\thinspace\null\nobreak\hfill
\hbox{\vbox{\kern-.2pt\hrule height.2pt
depth.2pt\kern-.2pt\kern-.2pt \hbox to1.8mm {\kern-.2pt\vrule
width.4pt \kern-.2pt\raise1.8mm\vbox to.2pt{} \lower0pt\vtop
to.2pt{}\hfil\kern-.2pt \vrule
width.4pt\kern-.2pt}\kern-.2pt\kern-.2pt \hrule height.2pt
depth.2pt \kern-.2pt}}\par\medbreak}

\title{Global boundedness of the gradient   for a class of nonlinear
  elliptic systems
} 
\numberwithin{equation}{section}

\author{
  Andrea Cianchi\\
 {\it Dipartimento di Matematica U.Dini, Universit\`a di Firenze}\\ {\it Piazza Ghiberti
27, 50122 Firenze, Italy} 
\bigskip
\\
  Vladimir G. Maz'ya \\
{\it   Department of Mathematics, Link\"oping University, SE-581
83 Link\"oping, Sweden}
}

\date{}

\begin{document}
\maketitle
\begin{abstract}

Gradient boundedness up to the boundary for solutions to Dirichlet
and Neumann problems for  elliptic systems with Uhlenbeck type
structure is established. Nonlinearities of possibly non-polynomial
type are allowed, and minimal regularity on the data and on the
boundary of the domain is assumed. The case of arbitrary bounded
convex domains is also included.
%
\end{abstract}

\unmarkedfootnote {
\par\noindent {\it Mathematics Subject
Classifications:} 35B45, 35J25.
\par\noindent {\it Keywords:} Nonlinear elliptic systems,
Dirichlet problems, Neumann problems, everywhere regularity,
boundedness of the gradient, Lipschitz continuity of solutions,
isoperimetric inequalities, convex domains, Orlicz-Sobolev spaces,
Lorentz spaces.}

\section{Introduction}\label{intro}
%


We are concerned with second-order nonlinear elliptic systems
of the form
\begin{equation}\label{system}
%
- {\rm {\bf div}} (a(|\nabla {\bf u}|)\nabla {\bf u} ) ={\bf f}(x)\quad
{\rm in}\,\,\, \o\,,
\end{equation}
coupled with either the Dirichlet condition $\u =0$, or the Neumann
condition $\frac{\partial {\bf u}}{\partial \nu} =0$ on $\partial
\Omega$.
Here, $\Omega $ is a domain, namely an open bounded connected  set
 in
$\rn$, $n \geq 2$,  ${\bf u}: \Omega \to \rN$,
$N \geq 1$, is a vector-valued unknown function, $\nabla {\bf u}:
\Omega \to \mathbb R^{Nn} $ denotes its  gradient, ${\bf f}: \Omega
\to \rN$ is a datum,
 ${\rm {\bf div}}$
stands for the $\rN$-valued divergence operator, and $\nu$  for the
outward unit normal to $\partial \Omega$.
\par
We prove the   boundedness of the gradient, or, equivalently, the
Lipschitz continuity, of the solutions to the relevant
boundary-value problems
 in the whole of $\Omega$. Quite general
nonlinearities of the differential operator, non-necessarily of
power type, are allowed, and  essentially weakest possible
integrability conditions on ${\bf f}$, and minimal regularity
assumptions on $\partial \Omega$ are imposed.
%
%
%
%
%
%
%
%
%
%
In the case when $\Omega$ is convex, no regularity on $\partial
\Omega$ is assumed at all.
%
%
The boundary value problems to be considered are the Euler equation
of variational problems for strictly convex integral functionals
depending on the  gradient only through its modulus, and hence the
solutions to the former agree with the minimizers of the latter. In
particular, our results on convex domains provide a version in the
vectorial case ($N>1$) of the so called semi-classical Hilbert-Haar
theory of minimization of strictly convex scalar integral
functionals of the modulus of gradient on convex domains in classes
of Lipschitz functions (see e.g. \cite[Chapter 1]{Giusti}).
\par
Nonlinear elliptic systems
 involving differential operators as in \eqref{system},  whose
 coefficient only depends on the modulus of the gradient,
 are sometimes referred to as systems with Uhlenbeck structure in the literature.
 Indeed, the  regularity theory of their solutions
%
can be  traced
back to the celebrated paper \cite{Uhl}.  Results from \cite{Uhl}
imply that,   if, for instance, $a(t)=t^{p-2}$ for some $p \geq 2$,
a choice which turns \eqref{system} into the $p$-Laplacian system,
%
then any local weak solution $\u$ to \eqref{system} satisfies
$\nabla \u \in L^\infty _{\rm loc}(\Omega, \mathbb R^{Nn})$, and, in
addition, $\nabla \u \in C^\alpha _{\rm loc}(\Omega, \mathbb
R^{Nn})$ for some $\alpha \in (0, 1)$. The regularity of solutions
to the  $p$-Laplacian  equation in the scalar case ($N=1$) had
instead  earlier been derived in \cite{Ur} from theorems on the
$C^\alpha$-regularity of solutions to non-uniformly elliptic
quasilinear systems.
 The contribution \cite{Uhl} was subsequently
extended to the situation when $1<p<2$ in \cite{AF, CDiB}. Further
generalization to elliptic systems  with non polynomial growth are
the subject of \cite{BreitSV, DieningSV, MingSiepe, Mar}.
 Precise inner  pointwise gradient estimates, via nonlinear and
linear potentials, for local solutions to nonlinear elliptic
equations, and to systems with Uhlenbeck structure, are the subject
of the recent papers \cite{DM, DMnew, KuM}.
\par 
Recall that, in striking contrast with the scalar case \cite{Di, Ev,
Le, To},
 solutions to nonlinear elliptic systems with a more general
structure than that of \eqref{system} can be irregular. Examples in
this connection are produced in \cite{SY}, where the existence of
nonlinear elliptic systems, with smooth coefficients depending only
on the gradient, but endowed with  solutions which are not even
bounde, is established. Earlier contributions on irregular solutions
to elliptic systems are
rooted in the paper
%
%
\cite{DeG}, and include \cite{GM} and \cite{Ne}.
%
%
Related examples, in the scalar case, of  linear and nonlinear
higher-order elliptic equations with irregular solutions were
independently exhibited in \cite{Macounter}.
\par
Let us incidentally mention that, however, solutions to elliptic
systems with general structure  are   well known to enjoy partial
regularity properties, in the sense that they are locally regular in
some open subset of $\Omega$ whose complement has zero Lebesgue
measure. This is the subject of a rich literature, starting with the
contributions \cite{HKW, GiaMo, Iv} -- see the monographs \cite{BF,
Gia, Giusti} for a comprehensive treatment of this topic. Recent
improvements of these results in terms of Hausdorff measures can be
found in \cite{Mi1, Mi2, KrM}.
\par
The study of global regularity, that is to say up to the boundary,
in boundary value problems for nonlinear elliptic systems has a more
recent history. Global gradient boundedness, and  H\"older
regularity, for the $p$-Laplacian elliptic system, under homogeneous
Dirichlet boundary conditions, were obtained in \cite{CDiB}, as a
consequence of analogous results for the associated parabolic
problem. Right-hand sides which are bounded in $x$, and domains
whose boundary is of class $C^{1,\alpha}$ are considered in that
paper. Further contributions on gradient regularity up to the
boundary
 for  systems and variational problems with
 Uhlenbeck structure, or perturbations of it,
 are \cite{BeiraoCrispo1,  Fo, FPV}. Partial boundary
regularity, i.e. regularity at the boundary outside subsets of zero
$(n-1)$-dimensional Hausdorff measure, for nonlinear elliptic
systems with general structure, is proved in  \cite{JM, DGK, KrM1}.
\par
The techniques employed in the   literature mentioned above,
  for both inner and boundary regularity, have a local
character. In particular,
the proofs of global results entail  distinguishing  between points
inside the domain, and boundary points, and   reducing the treatment
of the latter to the former. Such   reduction requires the use of
suitable local coordinates in which the boundary of the domain is
flat, and the structure of the differential operator is still close
enough to the Uhlenbeck type to ensure  everywhere inner regularity.
Techniques for inner regularity then apply after a reflection
argument, which allows to extend the solution beyond the flattened
boundary.
%
%
\par
By contrast, the approach of the present paper
 is   global in nature. Loosely speaking, an underlying
idea in our  proof of the global boundedness of the gradient
consists in integrating the system \eqref{system},  after
multiplication by $\Delta \u$, over the level sets of $|\nabla \u|$.
In particular, no localization via cut-off functions is employed.
Besides allowing for weak regularity assumptions on ${\bf f}$ and
$\Omega$,
%
%
%
 on approach of this kind enables us to deal not only with Dirichlet, but also
with Neumann boundary conditions, for which results seem  to be
missing in the literature. The  proof of the global boundedness of
the gradient in arbitrary convex domains
also relies upon the fact
%
that no local change of coordinates
near the
%
%
boundary is required. \par\noindent  Let us finally notice that
integration on the level sets of partial derivatives was used in
\cite{Magradient, Maconvex} to show gradient boundedness for linear
scalar problems, and in \cite{CMlipschitz} for  nonlinear scalar
problems.
 The approximation scheme exploited in those papers does not apply
  to the vectorial case. Here, we  follow an alternative
outline, which provides  a more self-contained proof also in the
scalar case.
%
%
%

%
%
%
%
%
%


 \section{Main results}\label{sec1}

Our assumptions on the system  \eqref{system} amount to what
follows. The function $a : (0, \infty ) \to (0, \infty )$ is
required to be
 monotone (either non-decreasing or non-increasing), of
class $C^1(0, \infty )$, and to fulfil
\begin{equation}\label{infsup}
-1 < i_a \leq s_a < \infty,
\end{equation}
where
 \begin{equation}\label{ia}
i_a= \inf _{t >0} \frac{t a'(t)}{a(t)} \qquad \hbox{and} \qquad
s_a= \sup _{t >0} \frac{t a'(t)}{a(t)}.
\end{equation}
In particular, the standard  $p$-Laplace operator for vector-valued
functions, corresponding to the choice $a(t)=t^{p-2}$, with $p >1$,
falls within this framework, since $i_a=s_a=p-2$ in this case.
\par\noindent
Thanks to the first inequality in \eqref{ia}, the function $b : [0,
\infty ) \to [0, \infty )$,  defined as
\begin{equation}\label{b} b(t)=a(t)t  \qquad  \hbox{if $t > 0$,\,\, and \, \, $b(0)=0$,}
\end{equation}
 turns out to be strictly increasing, and hence the
function $B: [0, \infty ) \to [0, \infty )$,  given by
\begin{equation}\label{B}
B(t) = \int _0^t b (\tau ) \, d\tau \qquad \hbox{for $t \geq 0$},
\end{equation}
is strictly convex. The  Orlicz-Sobolev space $W^{1, B}(\Omega ,
\rN)$ associated with the function $B$, or its subspace $W^{1,
B}_0(\Omega , \rN)$ of those functions vanishing in the suitable
sense on $\partial \o$, are  appropriate functional settings where
to define weak solutions to the boundary value problems associated
with the system \eqref{system}.
Precise definitions of function spaces and weak solutions are given
in Sections \ref{sec2} and \ref{sec4}, respectively;  existence and
 uniqueness of such solutions is also discussed in Section
 \ref{sec4}.
%
\par The right-hand side ${\bf f}$ in \ref{system} is assumed  to belong to the
Lorentz space $L^{n,1}(\o , \rN)$. This  space is borderline, in a
sense, for the family of Lebesgue spaces $L^q(\o , \rN)$ with $q>n$,
since $L^q(\o , \rN) \subsetneqq L^{n,1}(\o , \rN) \subsetneqq
L^{n}(\o , \rN)$ for every $q>n$. Let us mention that membership of
the right-hand side to the same Lorentz space is already known to
yield global gradient boundedness for solutions to scalar boundary
value problems \cite{CMlipschitz}. It has also been shown, via a
different approach relying upon potential theory, to ensure the
inner local boundedness of the gradient of local solutions to
equations, and to systems with Uhlenbeck structure  \cite{DMnew},
and also its continuity \cite{DMcont}.
\par The regularity of $\partial \Omega$ is  prescribed in terms
of a Lorentz space as well. We  impose that $\partial \Omega \in
W^{2}L^{n-1, 1}$. This means that $\Omega$
 is locally the subgraph of a function of $n-1$ variables whose
second-order distributional derivatives belong to the Lorentz space
$L^{n-1,1}$. This is the weakest possible integrability assumption
on the second-order derivatives of such a function for its
first-order derivatives to be continuous, and hence for $\partial
\Omega \in C^{1,0}$. Note that, by contrast, the available
regularity results at the boundary require  $\partial \Omega \in
C^{1,\alpha}$ for some $\alpha \in (0, 1]$.
\par
Let us emphasize that both the assumption on ${\bf f}$ and that on
$\Omega$ cannot be essentially relaxed for our conclusions to hold
-- see Remarks \ref{f} and \ref{omega} at the end of this section.
\par
Our  result for the Dirichlet problem
\begin{equation}\label{eqdirichlet}
\begin{cases}
- {\rm {\bf div}} (a(|\nabla {\bf u}|)\nabla {\bf u} ) ={\bf f}(x) &
{\rm in}\,\,\, \o\,,
\\
 \displaystyle {\bf u} =0 &
{\rm on}\,\,\,
\partial \o \,
\end{cases}
\end{equation}
 reads as follows.

\begin{theorem}\label{dirichletc2}
Let $\Omega $ be a domain in $\rn$, $n \geq 3$, such that $\partial
\Omega \in W^{2}L^{n-1, 1}$. Assume that ${\bf f} \in L^{n,1}(\Omega
, \rN)$. Let $\u$ be the (unique) weak solution to the Dirichlet
problem
%
%
%
%
%
\eqref{eqdirichlet}. Then there exists a constant $C$
$C=C(i_a, s_a,
\Omega )$
such that
\begin{equation}\label{gradienteqdirichlet}
\|\nabla \u\|_{L^\infty (\Omega , \R^{Nn})} \leq C b^{-1}\big(\|{\bf
f}\|_{L^{n,1}(\Omega , \rN)}\big).
\end{equation}
\par\noindent In particular, $\u$ is Lipschitz
continuous in $\Omega$.
\end{theorem}

\par
An interesting variant of Theorem \ref{dirichletc2} asserts that the
%
%
the regularity assumption on $\partial \Omega$ can be replaced by the convexity of $\Omega$.
This
is stated in the next result.

\begin{theorem}\label{dirichletconvex}
The same conclusion as in Theorem \ref{dirichletc2} holds if %
$\Omega$ is any  convex domain in $\rn$, $n \geq 3$.
\end{theorem}

Problem \eqref{eqdirichlet} is the Euler equation  of the
minimization problem for the strictly convex functional
\begin{equation}\label{functional}
J(u) = \int _\Omega \big( B(|\nabla \u|) - {\bf f}\cdot \u\big)\,
dx\,
\end{equation}
among trial functions $\u$ in $W^{1, B}_0(\Omega , \rN)$. Note that
$J$ is well defined in this function space under our assumption on
$\f$ -- the beginning of Section \ref{sec4}. The interpretation of
Theorem \ref{dirichletconvex} as an existence result for minimizers
of the functional $J$
in the space ${\rm Lip} _0(\Omega , \rN)$ of $\rN$-valued Lipschitz
continuous functions in $\Omega$ vanishing on $\partial \Omega$, to
which we alluded in Section \ref{intro}, is the content of the
following corollary.

%

\begin{corollary}\label{hilberhaar}
Let $\Omega$ be any convex domain in $\rn$, $n \geq 3$.
 Assume that ${\bf f} \in L^{n,1}(\Omega
, \rN)$. Then the functional $J$ admits a (unique) minimizer in the
space
${\rm Lip} _0(\Omega , \rN)$.
\end{corollary}

Results parallel to Theorems
\ref{dirichletc2}--\ref{dirichletconvex}  and Corollary
\ref{hilberhaar} hold for the solutions to the Neumann problem
\begin{equation}\label{eqneumann}
\begin{cases}
- {\rm {\bf div}} (a(|\nabla {\bf u}|)\nabla {\bf u} ) ={\bf f}(x) &
{\rm in}\,\,\, \o\,,
\\
\displaystyle\frac{\partial {\bf u}}{\partial \nu} =0 & {\rm
on}\,\,\,
\partial \o \,.
\end{cases}
\end{equation}
Clearly, here ${\bf f}$ has to fulfil the compatibility condition
\begin{equation}\label{intf0}
\int _\o {\bf f}(x)\, dx = 0.
\end{equation}

\begin{theorem}\label{neumannc2}
Let  $\Omega$ and ${\bf f}$ be as in Theorem \ref{dirichletc2}.
Assume, in addition, that \eqref{intf0} holds.
%
%
Let $\u$ be the (unique up to additive constant vectors) weak
solution to problem \eqref{eqneumann}. Then there exists a constant
$C=C(i_a, s_a, \Omega )$ such that
\begin{equation}\label{gradienteqneumann}
\|\nabla \u\|_{L^\infty (\Omega , \R^{Nn})} \leq C b^{-1}\big(\|{\bf
f}\|_{L^{n,1}(\Omega , \rN)}\big).
\end{equation}
In particular, $\u$ is Lipschitz continuous in $\Omega$.
\end{theorem}

A counterpart of Theorem \ref{neumannc2} for convex domains is
contained in the next result.

\begin{theorem}\label{neumannconvex}
The same conclusion as in Theorem \ref{neumannc2} holds if
 %
 $\Omega$ is any  convex domain in $\rn$, $n \geq 3$.
\end{theorem}


The minimization problem for the functional $J$  in the whole of
$W^{1, B}(\Omega , \rN)$ leads to the Euler equation
\eqref{eqneumann}. Hence, we have the following corollary of Theorem
\ref{neumannconvex}.

%
%
%

\begin{corollary}\label{hilberhaarneumann}
Let $\Omega$ be any  convex domain in $\rn$, $n \geq 3$.
 Assume that ${\bf f} \in
L^{n,1}(\Omega , \rN)$ and fulfils \eqref{intf0}. Then the
functional $J$ admits a (unique up to additive constant vectors)
minimizer in the class ${\rm Lip} (\Omega , \rN)$.
\end{corollary}

%

\begin{remark}\label{n=2}
{\rm Versions of the above results can be established via our
approach also in the case when $n=2$, under the slightly stronger
assumption that  ${\bf f} \in L^{q}(\Omega , \rN)$ for some $q>n$.
This becomes clear from a close inspection of the proofs.
}
\end{remark}

\begin{remark}\label{f}
{\rm The sharpness of assumption ${\bf f} \in L^{n,1}(\Omega , \rN)$
for the boundedness of the gradient of the solution to the Dirichlet
problem follows, in the linear case corresponding to the choice
$a=1$, from a result of \cite{Cmaxim} dealing with the scalar
Laplace equation in a ball.}
\end{remark}

\begin{remark}\label{omega}
{\rm The assumption  $\partial \Omega \in W^2 L^{n-1 , 1}$  is
 optimal for the boundedness of the gradient, as long as the regularity of $\Omega$ is prescribed in terms of integrability properties of its curvature. This
 can be demonstrated, again even just for scalar problems, by ad hoc
 examples of
 Dirichlet and Neumann problems for the $p$-Laplace equation in domains whose boundaries have conical
 singularities -- see e.g. \cite{CMJEMS}. Examples of the same nature  also show that the conclusion of Theorems \ref{dirichletconvex}
 and \ref{neumannconvex} may fail  under slight local non-smooth perturbations of
 convex domains  \cite{CMJEMS}.}
 \end{remark}

\section{Function spaces}\label{sec2}

\subsection{Spaces of measurable functions and rearrangements}\label{lorentzspaces}

\par
Let $(\MR , m )$ be a  positive, finite, non-atomic measure space.
The decreasing rearrangement $v^* : [0, \infty ) \to [0, \infty ]$
of    a real-valued $m$-measurable function $v$ on $\MR$ is the
unique right-continuous non-increasing function in $[0, \infty ) $
equidistributed with $v$. Namely, on defining the distribution
function $\mu _v: [0, \infty ) \to [0, \infty )$  of $v$ as
\begin{equation}\label{mu}
\mu _v (t) =m (\{x \in \MR :|v(x)|
>t\}) \qquad \quad \hbox{for $t\geq 0$,}
\end{equation}
we have that
\begin{equation}\label{rearr}
v^* (s) = \sup\{ t\geq 0 : \mu _v (t)> s\} \qquad {\rm {for}}\, \, s \in [0, \infty ).
\end{equation}
Clearly, $v^*(s)=0$ if $s \geq m (\MR )$.
\par\noindent
The function $v^{**} : (0, \infty ) \to [0, \infty )$, defined by
$$v^{**}(s) = \frac 1s\int _0^s v^*(r)\, dr \qquad \hbox{for
$s>0$},$$ is also nondecreasing, and such that $v^*(s) \leq
v^{**}(s)$ for $s
>0$.
\par\noindent The Hardy-Littlewood inequality is a basic property of rearrangements, and asserts that
\begin{equation}\label{HL}
\int_{\MR}|v(x) w(x)| dm (x) \leq \int_{0}^{\infty
}v^{\ast}(s)w^{\ast}(s)ds
\end{equation}
for all measurable functions $v$ and $w$ on $\MR$.
\par
Roughly speaking, a rearrangement-invariant space is a Banach
function space whose norm only depends on the rearrangement of
functions -- see e.g. \cite[Chapter 2]{BS} for a more precise
definition. Besides the Lebesgue spaces, their generalizations
provided by the Lorentz and the Orlicz spaces are classical
instances of rearrangement-invariant spaces which will play a role
in our discussion.
\par\noindent
Given $q \in (1, \infty )$ and $\sigma \in [1, \infty ]$, the
Lorentz space $L^{q, \sigma}(\MR )$ is the set of all real-valued
measurable functions $v$ on $\MR$ for which the quantity
\begin{equation}\label{lorentz}
\|v\|_{L^{q, \sigma}(\MR )} = \|s^{\frac 1q - \frac 1\sigma } v^*(s)
\|_{L^\sigma (0, m (\MR))}
\end{equation}
is finite. One has that  $L^{q, \sigma}(\MR )$ is a Banach space,
equipped with the norm, equivalent to $\|\cdot \|_{L^{q, \sigma}(\MR
)}$, obtained on replacing $v^*$ with $v^{**}$ on the right-hand
side of \eqref{lorentz}. Furthermore,
\begin{equation}\label{lorentz10}
L^{q,q}(\MR ) = L^q(\MR )\qquad \hbox{if $q \in (1, \infty)$,}
\end{equation}
\begin{equation}\label{lorentz11}
  L^{q, \sigma _1}(\MR )
\to L^{q, \sigma _2}(\MR ) \qquad \hbox{if $q \in (1, \infty )$ and
$\sigma _1 < \sigma _2$},
\end{equation}
and \begin{equation}\label{lorentz12}
L^{q_1, \sigma _1}(\MR ) \to
L^{q_2, \sigma _2}(\MR ) \qquad \hbox{if $q_1>q_2$ and $\sigma _1,
\sigma _2 \in [1, \infty]$}.
\end{equation}
 Here, and in what
follows, the arrow $\lq\lq \to "$ stands for continuous embedding.
Let us notice that the norm of the embedding \eqref{lorentz11}
depends on $q$, $\sigma _1$, $\sigma _2$, and the norm of the
embedding \eqref{lorentz12} depends  on $q_1$, $q_2$, $\sigma _1$,
$\sigma _2$ and $m(\MR)$.
\par\noindent
Denote by $q'$ and $\sigma '$
 the usual H\"older's conjugate
exponents of $q$ and $\sigma$. A H\"older type inequality in Lorentz
spaces tells us that
%
%
there exists a constant $C=C(q,  \sigma )$ such that
\begin{equation}\label{holderlorentz}
\int _{\MR} |v(x) w(x)|\, dm (x) \leq C \|v\|_{L^{q, \sigma }(\MR )}
\|w\|_{L^{q', \sigma ' }(\MR )}
\end{equation}
for every $v \in L^{q, \sigma }(\MR )$ and $w \in L^{q', \sigma '
}(\MR )$.
\par\noindent Given $N>1$, the Lorentz space $L^{q, \sigma}(\MR ,
\rN)$ of $\rN$-valued measurable functions on $\MR$ is defined as
$L^{q, \sigma}(\MR , \rN)  = (L^{q, \sigma}(\MR ))^N$, and is
endowed with the norm defined as $\|\v \|_{L^{q, \sigma}(\MR ,
\rN)}=\|\, |\v |\, \|_{L^{q, \sigma}(\MR)}$ for $\v \in L^{q,
\sigma}(\MR , \rN)$.
%
\par The
Orlicz spaces extend the Lebesgue spaces in the sense that the role
of powers in the definition of the norms is instead played by Young
functions. A Young  function $B : [0, \infty ) \to [0, \infty ]$ is
a convex function such that $B (0)=0$. If, in addition, $0< B(t) <
\infty$ for $t>0$ and
$$\lim _{t \to 0}\frac {B (t)}{t} =0 \qquad \hbox{and} \qquad \lim _{t \to \infty }\frac {B (t)}{t}
=\infty,$$ then $B$ is called an $N$-function. The Young conjugate
of a Young function $B$ is the Young function $\widetilde B$ defined
as
\begin{equation*}
\widetilde B (t) = \sup \{st - B (s): s \geq 0\} \qquad \hbox{for $t
\geq 0$.}
\end{equation*}
In particular, if  $B$ is an $N$-function, then $\widetilde B$ is an
$N$-function as well. Moreover, if $B$ is given by \eqref{B}, then
\begin{equation}\label{Btilde}
\widetilde B (t) = \int _0^t b^{-1}(s)\, ds \qquad \hbox{for $t \geq
0$.}
\end{equation}
Notice that
\begin{equation}\label{youngprop}
s \leq B^{-1}(s) \widetilde B^{-1}(s) \leq 2s \qquad \hbox{for $s
\geq 0$.}
\end{equation}
A Young function (and, more generally, an increasing function) $B$
is said to belong to the class $\Delta _2$ if there exists a
constant $C>1$ such that
\begin{equation*}
B (2t ) \leq C B (t) \qquad \hbox{for $t>0$.}
\end{equation*}
The Orlicz space $L^{B}(\MR )$ is the Banach function space of those
real-valued  measurable functions $v$ on $\MR$   whose Luxemburg
norm
\begin{equation*}
\|v\|_{L^B (\MR )} = \inf \bigg\{\lambda >0: \int _\MR B \bigg(\frac
{|v(x)|}\lambda \bigg) \,d m (x) \leq 1\bigg\}
\end{equation*}
is finite. The H\"older type inequality
\begin{equation}\label{holder}
\int _\MR |v(x) w(x)|\, dm(x) \leq 2 \|v\|_{L^B (\MR )}
\|w\|_{L^{\widetilde B} (\MR )}
%
\end{equation}
holds for every  $v \in L^B (\MR )$ and $w \in L^{\widetilde B} (\MR
)$. Let $B_1$ and $B_2$ be Young functions. Then
\begin{equation}\label{inclusion}
L^{B_1}(\MR ) \to L^{B_2}(\MR ) \,\,\, \hbox{if and only if} \,\,\,
\hbox{there exist $c, t_0 >0$ such that $B_2(t) \leq B_1(ct)$ for
$t>t_0$.}
\end{equation}
\par\noindent The Orlicz space $L^{B}(\MR ,
\rN)$, with $N>1$, of $\rN$-valued measurable functions on $\MR$ is
defined as $L^{B}(\MR , \rN)  = (L^{B}(\MR ))^N$, and is equipped
with the norm given by $\|\v \|_{L^B(\MR , \rN)}=\|\,|\v |\,
\|_{L^B(\MR)}$ for $\v \in L^B(\MR , \rN)$.

\subsection{Spaces of Sobolev type}\label{os}

Let $\Omega $ be a domain in $\rn$, with $n \geq 2$, and let $m \in
\N$. Sobolev type spaces of $m$-th order weakly differentiable
functions in $\Omega$, built upon Lorentz and Orlicz spaces, are
defined as follows.
\par\noindent
Given $q \in (1, \infty )$ and $\sigma \in [1, \infty ]$, the
Lorentz-Sobolev space
\begin{multline*}
W^mL^{q, \sigma }(\Omega ) = \{u \in L^{q, \sigma }(\Omega ):
\hbox{is $m$-times  weakly differentiable in $\Omega$} \\
\hbox{and $|\nabla ^k u| \in L^{q, \sigma }(\Omega )$ for $1 \leq k
\leq m$}\}\, \end{multline*}
  is a Banach space equipped with the norm
$\|u\|_{W^mL^{q, \sigma } (\Omega )} = 
\sum _{k=0}^m\|\, |\nabla ^k u|\,\|_{L^{q, \sigma }(\Omega )}.$
Here, $\nabla ^k u$ denotes the vector of all weak derivatives of
$u$ of order $k$. By $\nabla ^0 u$ we mean $u$. Moreover,
when $k=1$ we simply write $\nabla u$ instead of $\nabla ^1 u$.
\par\noindent
If $\sigma < \infty$, the space  $C^\infty (\Omega )\cap W^mL^{q,
\sigma }(\Omega )$ is dense in $W^mL^{q, \sigma }(\Omega )$. This
fact follows via an easy variant of a standard argument for
classical Sobolev spaces, and makes use the density of $C^\infty
_0(\Omega )$  in $L^{q, \sigma }(\Omega )$, and of a  version of
Young convolution inequality in Lorentz spaces due to O'Neil
\cite[Theorem 2.10.1]{Z}.
%
%
%
\par\noindent
A limiting case of the Sobolev embedding theorem asserts that if
$\Omega $ has a Lipschitz boundary, then $W^1L^{n,1}(\Omega ) \to
C^0(\Omega )$;
moreover, $L^{n,1}(\Omega )$ is optimal, in the sense that it is the
largest
  rearrangement invariant
 space enjoying this property  \cite{CP}. Hence, in particular,
\begin{equation}\label{lorentzemb}
W^2L^{n,1}(\Omega ) \to C^{1,0}(\Omega ), \end{equation}
 and
$L^{n,1}(\Omega )$ is optimal in the same sense as above.
\par\noindent
The Lorentz-Sobolev space $W^mL^{q, \sigma }(\Omega , \rN)$, $N>1$,
of $\rN$-valued functions in $\Omega$ is   defined as $W^mL^{q,
\sigma }(\Omega , \rN) = \big(W^mL^{q, \sigma }(\Omega )\big)^N$,
and endowed with the norm given by
$\|\u\|_{W^mL^{q, \sigma } (\Omega )} = 
\sum _{k=0}^m\|\, |\nabla ^k \u|\,\|_{L^{q, \sigma }(\Omega )}$ for
$\u \in W^mL^{q, \sigma }(\Omega , \rN)$.
\par
 The Orlicz-Sobolev space $W^{m, B
}(\Omega )$ is the Banach space
\begin{multline*}
W^{m, B}(\Omega ) =
 \{u \in L^B (\Omega ):
\hbox{is $m$-times  weakly differentiable in $\Omega$} \\
\hbox{and $|\nabla ^k u| \in L^B (\Omega )$ for $1 \leq k \leq
m$}\}\,, \end{multline*}
 and is equipped with the norm
$\|u\|_{W^{m, B} (\Omega )} = 
\sum _{k=0}^m\|\,|\nabla ^k u|\,\|_{L^B(\Omega )}.$
%
%
%
%
%
%
%
%
In what follows, we shall  only make use of first-order
Orlicz-Sobolev spaces $W^{1, B }(\Omega )$. By $W^{1, B }_0(\Omega
)$ and $W^{1, B }_\bot(\Omega )$ we denote the subspaces of $W^{1, B
}(\Omega )$ given by
\begin{equation*}
W^{1, B }_0(\Omega ) = \{u \in W^{1,B}(\Omega ): \hbox{the
continuation of $u$ by $0$ outside $\Omega$ is  weakly
differentiable  in $\rn$}\}\,,
\end{equation*}
and
\begin{equation*}
W^{1, B }_\bot (\Omega ) = \bigg\{u\in W^{1, B }(\Omega ):  \int
_\Omega u(x)\,dx =0\bigg\}.
\end{equation*}
A theorem of \cite{DT} ensures that, if  $B \in \Delta _2$, then
the space $C^\infty _0( \Omega )$ is dense in $W^{1,B}_0(\Omega)$,
and that, if $\Omega$ is a Lipschitz domain, then $C^\infty
(\overline \Omega )$ is dense in $W^{1,B}(\Omega)$.
\par\noindent
Let  $B$ be  a Young function such that
\begin{equation}\label{conv}
\int _0 \bigg(\frac {t}{B (t)}\bigg)^{\frac 1{n-1}} \,dt < \infty.
\end{equation}
The Sobolev conjugate of $B$, introduced in \cite{Cbound} (and, in
an equivalent form, in \cite{Csharp}), is the Young function $B_n$
defined as
\begin{equation}\label{sobolevconj}
B_n (t) = B\big(H_n^{-1}(t)\big) \qquad \hbox{for $t \geq 0$,}
\end{equation}
where
 \begin{equation}\label{5.106}
H_n (s)= \bigg( \int _0 ^s \bigg(\frac {t}{B (t)}\bigg)^{\frac
1{n-1}} \,dt \bigg)^{1/n'} \qquad \hbox{for $s \geq 0$,}
\end{equation}
and $H_n^{-1}$ denotes the (generalized) left-continuous inverse of
$H_n$.
\par\noindent
An embedding theorem for Orlicz-Sobolev spaces \cite{Csharp, Cbound}
tells us that, if
%
%
$B$  fulfils \eqref{conv}, then there exists a constant $C=C(n,
|\Omega |)$ such that
\begin{equation}\label{embedding0}
\|u\|_{L^{B_n}(\Omega )} \leq C \|\, |\nabla u|\,\|_{L^B(\Omega )}
\end{equation}
for every $u \in  W^{1, B}_0(\Omega )$. Moreover, if has a Lipschitz
boundary, then inequality \eqref{embedding0} holds for every $u \in
W^{1, B}_\bot(\Omega )$. The space $L^{B_n}(\Omega )$ is optimal in
\eqref{embedding0} among all Orlicz spaces. Note that assumption
\eqref{conv} is, in fact, immaterial in \eqref{embedding0}, since,
owing to \eqref{inclusion}, the Young function $B$ can be replaced,
if necessary, with another Young function fulfilling \eqref{conv} in
such a way that $ W^{1, B}(\Omega )$ remains unchanged, up to
equivalent norms.
\par \noindent
If $B$ is a   Young function, then
\begin{equation}\label{incl}
L^\infty (\Omega ) \to L^{B_n}(\Omega ) \to L^{n'}(\Omega ).
\end{equation}
The first embedding in \eqref{incl} is trivial. As for the second
one, since $B$ is a  Young function, there exist constants $c_0$ and
$t_0
>0$  such that $t \leq B(c_0t)$ if $t >t_0$. As a consequence,  there  exist constants $c_1$ and $t_1$ such that $t^{n'} \leq
B_n(c_1t)$ for some $t>t_1$. Hence, the second embedding in
\eqref{incl} follows via  \eqref{inclusion}.
 \par \noindent Let us also observe that,
if $B$ grows so fast near infinity that
\begin{equation}\label{5.104'}
\int ^\infty \bigg(\frac {t}{B (t)}\bigg)^{\frac 1{n-1}} \,dt
<\infty,
\end{equation}
then  equality holds in the first embedding in \eqref{incl}. Indeed,
under \eqref{5.104'}, $H_n^{-1}(t)= \infty$ for large $t$, and hence
$B_n(t)= \infty$ for large $t$ as well.
Hence, $L^{B_n}(\Omega )= L^\infty (\Omega )$, up to equivalent
norms.
\par\noindent
The Orlicz-Sobolev space $W^{m, B} (\Omega , \rN)$ is defined, for
$N>1$, as $W^{m, B}(\Omega , \rN) = \big(W^{m, B}(\Omega )\big)^N$,
and equipped with the norm $\|\u\|_{W^{m, B} (\Omega )} =  \sum
_{k=0}^m\|\,|\nabla ^k \u|\,\|_{L^B(\Omega )}.$ The spaces $W^{1,
B}_0(\Omega , \rN)$ and $W^{1, B}_\bot (\Omega , \rN)$ are defined
accordingly.

\section{The function $a$}\label{seca}

This section is devoted to the proof of some properties of the
function $a$ appearing in \eqref{system}. Hereafter, $b$ and $B$
 denote the
functions associated with $a$ as in \eqref{b} and  \eqref{B}.
%
%
 Furthermore, we define the function
  $H: [0, \infty ) \to [0, \infty )$ as
\begin{equation}\label{Hdef}
H(t) = \int _0^t a(\tau )b(\tau ) \, d\tau \qquad \quad \hbox{for $t
\geq 0$,}
\end{equation}
and the function $F: [0, \infty ) \to [0, \infty )$ as
\begin{equation}\label{F}
F(t) = \int _0^t b(\tau )^2 \, d\tau \qquad \quad \hbox{for $t \geq
0$.}
\end{equation}

\begin{proposition}\label{abis}
Assume that the function $a : (0, \infty ) \to (0, \infty )$ is of
class  $C^1$ and fulfils \eqref{infsup}. Let $b$ and $B$ be the
functions given  by \eqref{b} and \eqref{B}, respectively, and let
$H$ and $F$ be defined as above. Then:
\par\noindent
(i)
\begin{equation}\label{compatti}
a(1)  \min\{t^{i_a},  t^{s_a}\} \leq a(t) \leq a(1) \max\{t^{i_a},
t^{s_a}\} \quad \hbox{for $t>0$.}
\end{equation}
\par\noindent
(ii) $b$ is increasing,
\begin{equation}\label{limiteb}
\lim _{t \to 0} b(t) = 0, \quad \hbox{and} \quad \lim _{t \to
\infty} b(t) = \infty.
\end{equation}
\par\noindent
(iii) $B$  is a strictly convex $N$-function,
\begin{equation}\label{delta2}
B\in \Delta _2 \qquad \hbox{and} \qquad \widetilde B\in \Delta _2.
\end{equation}
\par\noindent
(iv)  There exists a constant $C=C(i_a , s_a)$ such that
\begin{equation}\label{conj}
\widetilde{B}(b(t)) \leq C B(t) \qquad \quad \hbox{for $t \geq 0$.}
\end{equation}
(v) For every $C>0$, there exists  a positive constant $C_1=C_1(s_a,
C)>0$ such that
\begin{equation}\label{B7}
C b^{-1}(s) \leq b^{-1}(C_1 s) \qquad \qquad \hbox{for $s>0$,}
\end{equation}
and a positive constant $C_2=C_2(i_a, C)>0$ such that
\begin{equation}\label{B7bis}
 b^{-1}(Cs) \leq C_2b^{-1}(s) \qquad \qquad \hbox{for $s>0$.}
\end{equation}
(vi)  There exists a positive constant $C=C(s_a)$ such that
\begin{equation}\label{B5}
B(t) \leq tb(t) \leq C B(t) \qquad \hbox{for $t\geq 0$.}
\end{equation}
\par\noindent
(vii) There exists a positive constant $C=C(i_a, s_a)$ such that
\begin{equation}\label{B1}
F(t) \leq t b(t)^2 \leq C F(t) \qquad \hbox{for $t\geq 0$.}
\end{equation}
\par\noindent
(viii) There exist  positive constants $C_1=C_1(i_a)$  and
$C_2=C_2(i_a, s_a)$such that
\begin{equation}\label{H}
 C_1 H(t) \leq   b(t)^2 \leq C_2 H(t) \qquad \hbox{for $t\geq 0$.}
\end{equation}
\end{proposition}
{\bf Proof}. Assertions (ii)--(vii) are proved in \cite[Propositions
2.9 and 2.15]{CMlipschitz}. Property (i) can be shown on
distinguishing the case when $t \in (0,1)$ and $t \in [1, \infty)$,
and integrating the inequality
$$\frac{i_a}{\tau} \leq \frac{a'(\tau)}{a(\tau)} \leq \frac{s_a}{\tau} \quad \hbox{for $\tau >0$,}$$
on $(t, 1)$ and on $(1, t)$, respectively.
\par\noindent As far as (viii) is concerned, since $b'(t) = a(t) + t a'(t)$ for $t >0$, one has that
\begin{equation}\label{ab'}
\frac{b'(t)}{1+s_a} \leq a(t) \leq \frac{b'(t)}{1+\min \{i_a , 0\}}
\quad \hbox{for $t >0$,}
\end{equation}
and hence
\begin{equation}\label{ab'bis}
\frac{b (t)}{1+s_a} \leq \int _0^ta(\tau)\, d\tau \leq \frac{b
(t)}{1+\min \{i_a , 0\}} \quad \hbox{for $t >0$.}
\end{equation}
Thus, inasmuch as $b$ is increasing,
$$H(t) \leq b(t) \int _0^t a(\tau )\, d\tau \leq \frac{b
(t)^2}{1+\min \{i_a , 0\}} \quad \hbox{for $t >0$,}
$$
whence the first inequality in \eqref{H} follows.
 On the other hand, integration by parts and inequalities
 \eqref{ab'} and
 \eqref{ab'bis}
yield
\begin{align*}
H(t) & = \int _0^t a(\tau) b(\tau)\, d\tau   =   b(t)\int _0^t
a(\tau)\, d\tau - \int _0^t b'(\tau)\int _0^\tau a(r)\, dr\, d\tau \\
\nonumber & \geq \frac{b (t)^2}{1+s_a} - \int _0^t (1+s_a)
a(\tau)\frac{b (\tau)}{1+\min
\{i_a , 0\}}d\tau 
= \frac{b (t)^2}{1+s_a} - \frac{1+s_a}{1+\min \{i_a , 0\}} H(t)
\quad \hbox{for $t >0$.}
\end{align*}
This implies the second inequality in \eqref{H}. \qed

\begin{lemma}\label{elliptic}
Let $a$ be as in Lemma \ref{abis}.
%
Then
\begin{equation}\label{elliptic1}
(1 + \min\{i_a , 0\}) a(|\xi |) |\eta|^2 \leq \sum _{\alpha , \beta
=1}^N \sum _{i,j =1}^n \frac{\partial (a(|\xi |)\xi
_{i}^\alpha)}{\partial \xi _{j}^\beta}\eta _{i}^\alpha \eta
_{j}^\beta \leq (1 + \max\{s_a , 0\}) a(|\xi |) |\eta|^2
\end{equation}
for $\xi, \eta \in \mathbb R^{Nn}$. Moreover,
\begin{equation}\label{elliptic2}
[a(|\xi |)\xi   - a(|\eta |)\eta] \cdot (\xi - \eta ) \geq (1 +
\min\{i_a , 0\})  |\xi - \eta|^2 \int _0^1 a(|\eta + s(\xi - \eta
)|)ds
\end{equation}
for $\xi, \eta \in \mathbb R^{Nn}$.
\end{lemma}
\par\noindent {\bf Proof}. Given $i, j, \alpha , \beta$, one has that
 \begin{equation}\label{elliptic3} \frac{\partial (a(|\xi
|)\xi _{i}^\alpha)}{\partial \xi _{j}^\beta} = \frac{a'(|\xi|)}{|\xi
|} \xi _i^\alpha \xi _j^\beta + a(|\xi |) \delta _{ij} \delta
_{\alpha \beta} \quad \hbox{for $\xi \in \mathbb R^{Nn}.$}
\end{equation}
Thus
\begin{multline}\label{elliptic4}
\sum _{\alpha , \beta =1}^N \sum _{i,j =1}^n \frac{\partial (a(|\xi
|)\xi _{i}^\alpha)}{\partial \xi _{j}^\beta}\eta _{i}^\alpha \eta
_{j}^\beta
 = \sum _{\alpha , \beta
=1}^N \sum _{i,j =1}^n \frac{a'(|\xi|)}{|\xi |} \xi _i^\alpha \xi
_j^\beta \eta _i^\alpha \eta_j^\beta  + a(|\xi |) \sum _{\alpha=1}^N
\sum _{i=1}^n (\eta_i^\alpha )^2 \\ = \frac{a'(|\xi|)}{|\xi |}  (\xi
\cdot \eta )^2 + a(|\xi |) |\eta|^2 \quad \hbox{for $\xi , \eta \in
\mathbb R^{Nn}.$}
\end{multline}
Since $|\xi \cdot \eta | \leq |\xi | |\eta |$, equation
\eqref{elliptic1} follows via \eqref{infsup}.
\par\noindent
Next, set $A_{ij}^{ \alpha \beta}(\zeta ) = \frac{\partial (a(|\zeta
|)\zeta _{i}^\alpha)}{\partial \zeta _{j}^\beta}$ for $\zeta \in
\mathbb R ^{Nn}$. Given $i$ and $\alpha$, we have that
$$[a(|\xi |)\xi _i^\alpha  - a(|\eta |)\eta _i^\alpha ] = \sum _{\beta =1}^N \sum _{j=1}^n  (\xi _j^\beta- \eta _j^\beta
)\int _0^1 A_{ij}^{ \alpha \beta}(\eta + s(\xi - \eta ))\, ds \quad
\hbox{for $\xi , \eta \in \mathbb R^{Nn}.$}.$$
Therefore,
  $$[a(|\xi |)\xi   - a(|\eta |)\eta] \cdot (\xi - \eta )
= \sum _{\alpha , \beta =1}^N \sum _{i,j =1}^n (\xi _i^\alpha- \eta
_i^\alpha ) (\xi _j^\beta- \eta _j^\beta )\int _0^1 A_{ij}^{ \alpha
\beta}(\eta + s(\xi - \eta ))\, ds \quad  \hbox{for $\xi , \eta \in
\mathbb R^{Nn},$}
$$
and hence inequality \eqref{elliptic2} follows via
\eqref{elliptic1}. \qed

\begin{lemma}\label{ineqa}
Let $a$ be as in Lemma \ref{abis}. Assume in addition that $a$ is
non-decreasing.
%
 Then
\begin{align}\label{ineqa2}
[a(| \xi|)\xi -a(|\eta |)\eta ]\cdot (\xi-\eta) \geq \frac 13 [a(|
\xi|) +a(|\eta |)]\, |\xi - \eta |^2 \quad \hbox{for $\xi, \eta \in
\mathbb R^{Nn}$.}
\end{align}
%
%
%
%
\end{lemma}
\noindent {\bf Proof}.
Since inequality \eqref{ineqa2} is invariant under replacements of
$\xi$ and $\eta$ by each other, we may assume, without loss of
generality, that $|\xi | \geq |\eta|$, and hence $a(|\xi |) \geq
a(|\eta|)$. Consider first the case when $a(|\xi |) \leq 2 a(|\eta
|)$. Then, given any  $\xi \neq \eta $,
\begin{align}\label{ineqa3}
\frac{(a(|\xi |)\xi - a(|\eta|)\eta ) \cdot (\xi - \eta )}{(a(|\xi
|) + a(|\eta |)) |\xi - \eta |^2}
  & =
\frac{(\frac {a(|\xi |)}{a(|\eta|)}\xi - \eta) \cdot (\xi - \eta
)}{(\frac{a(|\xi |)}{a(|\eta |)} + 1) |\xi - \eta |^2}
  \geq
\frac{(\frac {a(|\xi |)}{a(|\eta|)}\xi - \eta) \cdot (\xi - \eta
)}{3 |\xi - \eta |^2}
\\ \nonumber & =
\frac{|\xi - \eta |^2 + (\frac {a(|\xi |)}{a(|\eta|)} - 1)\xi \cdot
(\xi - \eta )}{3 |\xi - \eta |^2}
 =
\frac 13 + \frac{(\frac {a(|\xi |)}{a(|\eta|)} - 1) (|\xi |^2 - \xi
\cdot \eta )}{3 |\xi - \eta |^2}
  \\ \nonumber & \geq
\frac 13 + \frac{(\frac {a(|\xi |)}{a(|\eta|)} - 1) (|\xi |^2 -
|\xi| |\eta | )}{3 |\xi - \eta |^2}
 \geq
\frac 13\,.
\end{align}
Assume  now that $a(|\xi |) \geq 2 a(|\eta |)$.  Then, given any
$\xi \neq \eta $,
\begin{align}\label{ineqa4}
\frac{(a(|\xi |)\xi - a(|\eta|)\eta ) \cdot (\xi - \eta )}{(a(|\xi
|) + a(|\eta |)) |\xi - \eta |^2}
  & =
\frac{(\frac {a(|\xi |)}{a(|\eta|)}\xi - \eta) \cdot (\xi - \eta
)}{(\frac{a(|\xi |)}{a(|\eta |)} + 1) |\xi - \eta |^2}
 \geq
\inf _{s \geq 2} \frac{(s \xi - \eta )\cdot (\xi - \eta )}{(s+1)
|\xi - \eta |^2}
 =
\frac{(2 \xi - \eta )\cdot (\xi - \eta )}{3 |\xi - \eta |^2}
\\ \nonumber & =
\frac{2 |\xi |^2  - 3 \xi \cdot \eta + |\eta |^2 }{3 |\xi - \eta
|^2}
 =
\frac{|\xi - \eta |^2  + |\xi |^2 -  \xi \cdot \eta }{3 |\xi - \eta
|^2}
\\ \nonumber & \geq
\frac 13 + \frac{|\xi |^2 - |\xi| |\eta | }{3 |\xi - \eta |^2}
 \geq
\frac 13\,.
\end{align}
Inequality \eqref{ineqa2} is fully proved. \qed

\begin{lemma}\label{approxbis}
Let $a$  be as in Lemma \ref{abis}. Assume, in addition, that $a$ is
monotone (either non-decreasing or non-increasing). Then, for every
$t, \tau >0$, there exists a positive constant $\vartheta =
\vartheta (a, t, \tau )$ such that
\begin{equation}\label{approxbis1}
\inf\{[a(| \xi|)\xi -a(|\eta |)\eta ]\cdot (\xi-\eta): \,\xi, \eta
\in \mathbb R^{Nn},
 |\xi - \eta |\geq t, |\xi |\le \tau , |\eta |\le \tau \} >\vartheta\,.
 \end{equation}
\end{lemma}

\noindent {\bf Proof}.
Assume first that $a$  is non-decreasing. In particular, $s_a \geq
i_a \geq  0$. Then, by Lemma \ref{ineqa},
\begin{align}\label{approxbis3}
[a (| \xi|)\xi -a (|\eta |)\eta ]\cdot (\xi-\eta) \geq \tfrac 13 [a
(| \xi|) +a (|\eta |)]\, |\xi - \eta |^2 \quad \hbox{for $\xi, \eta
\in \mathbb R^{Nn}$}.
\end{align}
Since, by \eqref{delta2} and \eqref{B5}, $a\in \Delta _2$,
\begin{equation}\label{approxbis4}
a(|\xi - \eta |) \leq a(|\xi | + |\eta |) \leq a(2|\xi |) + a(2|\eta
|) \leq C[a(|\xi |) + a(|\eta |)] \quad \hbox{for $\xi, \eta \in
\mathbb R^{Nn}$.}
\end{equation}
By \eqref{approxbis3}, \eqref{approxbis4} and  the first inequality
in \eqref{compatti},
\begin{align}\label{approxbis5}
[a(| \xi|)\xi -a(|\eta |)\eta ]\cdot (\xi-\eta)  \geq C \min\{ |\xi
- \eta |^{i_{a}+2}, |\xi - \eta |^{s_{a}+2}\}
\quad \hbox{for $\xi, \eta \in \mathbb R^{Nn}$,}
\end{align}
for some positive constant $C=C(a)$. Hence \eqref{approxbis1}
follows, since $\{\xi, \eta \in \mathbb R^{Nn}:\,  |\xi - \eta |\geq
t, |\xi |\le \tau , |\eta |\le \tau \}$ is a compact set.
\par\noindent
Assume next that
 $a$  is non-increasing.
 Thus, $0 \geq
s_{a} \geq i_a$. By \eqref{elliptic2},
\begin{align}\label{approxbis9}
[a (| \xi|)\xi -a (|\eta |)\eta ]\cdot (\xi-\eta) \geq (1 + \min
\{i_{a} , 0\})|\xi - \eta |^2 \int _0^1a (|\eta + t(\xi - \eta )|)\,
dt \quad \hbox{for $\xi, \eta \in \mathbb R^{Nn}$.}
\end{align}
Owing to \eqref{compatti},
\begin{align}\label{approxbis10}
a (|\eta + t(\xi - \eta )|)
 & \geq C \min\{ |\eta +
t(\xi - \eta ) |^{i_{a}}, |\eta + t(\xi - \eta ) |^{s_{a}}\}
 \\ \nonumber & \geq C \min\{ (|\eta| +
|\xi - \eta |) ^{i_{a}}, (|\eta| + |\xi - \eta |)^{s_{a}}\}
\\ \nonumber & \geq C' \min\{ (|\eta| + |\xi |) ^{i_{a}}, (|\eta| + |\xi |) ^{s_{a}}\},
\end{align}
for some positive constants $C=C(a)$ and $C'=C'(a)$. Coupling
\eqref{approxbis9} with \eqref{approxbis10} yields
$$[a(| \xi|)\xi -a(|\eta |)\eta ]\cdot
(\xi-\eta) \geq C \min\{ (|\eta| + |\xi |) ^{i_{a}}, (|\eta| + |\xi
|) ^{s_{a}}\} |\xi - \eta |^2 \quad \hbox{for $\xi, \eta \in \mathbb
R^{Nn}$,}$$ for some positive constant $C=C(a)$. Hence,
\eqref{approxbis1} follows also in this case. \qed

\medskip
\par
In the following lemma, any function $a$ as in the statement of
Theorem \ref{dirichletc2} is approximated by a
 family $\{a_\varepsilon
\}$ of  functions enjoying the  additional property of being bounded
from above and from below by positive constants.
%
%

\begin{lemma}\label{approx}
Let $a$  be as in Lemma \ref{abis}. Assume, in addition, that $a$ is
monotone (either non-decreasing or non-increasing). Given
$\varepsilon \in (0,1)$, define $a_\varepsilon : [0, \infty ) \to
(0, \infty )$ as
\begin{equation}\label{aeps}
a_\varepsilon (t) = \frac{a(\sqrt{\varepsilon + t^2}) +
\varepsilon}{1 + \varepsilon a(\sqrt{\varepsilon + t^2})} \quad
\quad \hbox{for $t \geq 0$.}
\end{equation}
 Then $a_\varepsilon $ has the same monotonicity property as $a$,
%

\begin{equation}\label{cinf}
a_\varepsilon \in C^1 ([0 , \infty )),
\end{equation}


\begin{equation}\label{abound}
\varepsilon \leq    a_\varepsilon (t)   \leq \varepsilon ^{-1} \quad
\hbox{for $t \geq 0$,}
\end{equation}

\begin{equation}\label{indici}
\min \{i_a , 0\} \leq i_{a_\varepsilon} \leq s_{a_\varepsilon} \leq
\max \{s_a , 0\},
\end{equation}

\begin{equation}\label{conva}
 \lim _{\varepsilon \to 0} a_\varepsilon (|\xi |) \xi = a (|\xi |) \xi \qquad
\hbox{uniformly in $\{\xi \in \R ^{Nn} : |\xi | \leq M\}$  for every
$M>0$.}
\end{equation}
\par\noindent
Moreover, if $b_\varepsilon $ and $B_\varepsilon$ are defined as in
\eqref{b} and \eqref{B}, respectively, with $a$ replaced with
$a_\varepsilon$, then

\begin{equation}\label{convb}
 \lim _{\varepsilon \to 0} b_\varepsilon = b \qquad
\hbox{uniformly in $[0, M]$  for every $M>0$,}
\end{equation}
and hence
\begin{equation}\label{convB}
 \lim _{\varepsilon \to 0} B_\varepsilon = B \qquad \hbox{uniformly in $[0, M]$   for every $M>0$.}
\end{equation}


\end{lemma}

\par\noindent
{\bf Proof}. Property \eqref{cinf}  trivially follows from the fact
that $a \in C^1(0, \infty)$. Since
\begin{equation}\label{aeps'}
a_\varepsilon '(t) = \frac{(1 - \varepsilon
^2)\,a'(\sqrt{\varepsilon + t^2})\, t}{\big(1 + \varepsilon
a(\sqrt{\varepsilon + t^2})\big)^2\,\sqrt{\varepsilon + t^2}} \quad
\hbox{for $t \geq 0$,}
\end{equation}
$a'$ and $a_\varepsilon'$ have like signs, and hence $a$ and
$a_\varepsilon$ share the same   monotonicity property. Equation
\eqref{indici} is an easy consequence of \eqref{aeps'} and of the
very definitions of $i_{a_\varepsilon}$ and $s_{a_\varepsilon}$.
Equation \eqref{abound} follows from the definition of
$a_\varepsilon$, and from the fact that the function $[0, \infty )
\ni s \mapsto \frac{s+\varepsilon}{1+ \varepsilon s}$ is increasing
for every $\varepsilon \in (0, 1)$.
\par\noindent
Next, note that
 $$\lim _{\varepsilon \to 0} a_\varepsilon = a \qquad
\hbox{uniformly in $[L, M]$  for every $M>L>0$.}$$
Hence,
 \begin{equation}\label{unif2}
 \lim _{\varepsilon \to 0}  b _\varepsilon = b \quad \hbox{uniformly in $[L, M]$ for every $M>L>0$.}
 \end{equation}
On the other hand, by \eqref{compatti} with $a$ replaced with
$a_\varepsilon$ and by \eqref{indici},
\begin{equation}\label{unif3}
0 \leq b_\varepsilon (t) = t a_\varepsilon (t) \leq a_\varepsilon
(1) t^{1+ \min\{i_a , 0\}} \leq (\max\{a(\sqrt{2}), a(1)\}+1)t^{1+
\min\{i_a , 0\}}
%
%
\,\, \quad \hbox{if  $0< t <1$,}
\end{equation}
whence
\begin{equation}\label{unif4}
\lim _{t \to 0}  b_\varepsilon (t) =0 \qquad \quad \hbox{uniformly
for $\varepsilon \in (0, 1)$.}
\end{equation}
Combining \eqref{unif2}, \eqref{unif4} and \eqref{limiteb} yields
\eqref{convb}.
 \par\noindent The proof of \eqref{conva} is
analogous. \qed

\section{Fundamental geometric and differential inequalities}\label{differential}

Here, we enucleate some inequalities of geometric and functional
nature which are needed in the proofs of our main results.
%
\par\noindent We begin with a relative isoperimetric inequality, which tells us that
if $\Omega$ is an open subset of $\rn$, $n \geq 2$, with a Lipschitz
boundary, then there exists a constant $C$ such that
\begin{equation}\label{relisop}
|E|^{1/n'} \leq C \hh (\partial ^ME  \cap \Omega)
\end{equation}
for every measurable  set $E \subset \Omega$
such that $|E| \leq \m2$ \cite[Corollary 5.2.1/3]{Mazlibro}. Here,
$|E|$ denotes  the Lebesgue measure of $E$,  $\partial ^ME$ its
essential boundary, and $\hh$ stands for
   $(n-1)$-dimensional Hausdorff measure.
\par\noindent Inequality \eqref{relisop} can be derived via another
geometric inequality, which holds in any Lipschitz domain $\Omega$,
and asserts that
\begin{equation}\label{isopboundary}
\hh (\partial ^M E \cap \partial \Omega) \leq C \hh (\partial ^ME
\cap \Omega)
\end{equation}
 for some constant $C=C(\Omega)$,  and for every measurable set $E
\subset \Omega$
such that $|E| \leq \m2$  \cite[Chapter 6]{Mazlibro}. Indeed,
\begin{equation}\label{boundary}
\partial ^M E = (\partial ^ME \cap \partial \o )\cup  (\partial ^ME \cap  \o
),
\end{equation}
for every measurable  set $E \subset \Omega$, the union being  disjoint, and hence
\begin{equation}\label{boundarymeasure}
\hh (\partial ^M E) = \hh(\partial ^ME \cap \partial \o ) + \hh
(\partial ^ME \cap  \o ),
\end{equation}
inasmuch as $\hh$ is a measure when restricted to Borel sets.
%
%
Thus, inequality \eqref{relisop} follows
from \eqref{isopboundary}, \eqref{boundarymeasure}, and the
 classical isoperimetric inequality in $\rn$, which takes the form
 $$|E|^{1/n'} \leq C \hh ( \partial ^ME )$$
 for some constant $C=C(n)$, and for every measurable set $E$ in $\rn$
 with  finite measure. Notice that, in particular, the constant in
 \eqref{relisop} depends on $n$ and on the constant in
 \eqref{isopboundary}.

A trace inequality for functions from the Sobolev space $
W^{1,2}(\Omega )$, whose support has measure not exceeding $\m2$, is
the content of the following lemma. In the statement,  ${\rm Tr}\,
v$ denotes the trace  on $\partial \Omega$ of a function $v$, and
${\rm supp}\, v$ its support.

\begin{lemma}\label{trace}
Let $\Omega $ be a domain with a Lipschitz boundary  in $\rn$, $n
\geq 2$. Assume that either $1 \leq q \leq \frac{2(n-1)}{n-2}$, or
$1 \leq q < \infty$, according to whether $n\geq 3$ or $n=2$.
 Then there exists a
constant $C$, depending on $n$, $q$ and on the constant   in
\eqref{isopboundary},
%
%
such that
\begin{equation}\label{11}
\bigg(\int _{\partial \Omega} |{\rm Tr} \,v|^{q}\, d\hh
(x)\bigg)^{\frac 1q}\leq C |{\rm supp}\, v |^{\frac {n-1}{qn} -
\frac{n-2}{2n} } \bigg(\int _\Omega |\nabla v|^2 dx\bigg)^{\frac 12}
\end{equation}
for every $v \in W^{1,2}(\Omega )$ satisfying $|{\rm supp}\, v |
\leq \m2$.
\end{lemma}
\par\noindent
{\bf Proof of Lemma \ref{trace}}
There exists a constant $C$, depending on $n$, $q$
and on the constant in \eqref{isopboundary},
%
such that
\begin{equation}\label{12}
\bigg(\int _{\partial \Omega} |{\rm Tr}\, v|^{q} d\hh \bigg)^{1/q}
\leq C \bigg(\int _\Omega |\nabla v|^{\frac{nq}{q+n-1}}
dx\bigg)^{\frac{q+n-1}{nq}}
\end{equation}
%
for every $v \in W^{1,2}(\Omega )$ fulfilling $|{\rm supp}\, v |
\leq \m2$. Inequality \eqref{12} can be derived from a subsequent
application of \cite[Theorem 6.11.4/1]{Mazlibro} to $|v|^q$ (with
$p=1$),   of H\"older's inequality, and of \cite[Lemma
6.2]{Mazlibro}. With inequality \eqref{12} in place, inequality
\eqref{11} follows  via H\"older's inequality. \qed
%

\medskip
\par
If $\u \in W^{2,1}(\Omega , \rN)$, then $|\nabla \u| \in
W^{1,1}(\Omega)$, by the chain rule for vector-valued functions
\cite[Theorem 2.1]{MarcusMizel}. 
An application of the coarea formula for Sobolev functions in the
form of \cite{BZ}  then tells us that, for every Borel function $g :
\Omega \to [0, \infty )$,
\begin{equation}\label{coareanabla}
 \int _{\{|\nabla \u | >t\}} g(x) |\nabla |\nabla \u || dx = \int _{t}^{\infty}
\int _{\{|\nabla \u | =\tau\}} \,g(x)\, d\hh (x)\,d\tau \quad
\hbox{for $t \geq 0$,}
\end{equation}
provided that a suitable precise representative of the function
$|\nabla \u |$ is employed.  Hence, if the left-hand side is finite
for $t>0$, then it is a (locally) absolutely continuous function of
$t$, and
\begin{equation}\label{38}
  - \frac{d}{dt} \int _{\{|\nabla \u |
>t\}}g(x)|\nabla
|\nabla \u || dx = \int _{\{|\nabla \u | =t\}} g(x)\, d\hh (x)
\qquad \hbox{for a.e. $t>0$}.
\end{equation}
The use of the coarea formula again tells us that $\hh (\{|\nabla \u
| =t\} \cap \{|\nabla |\nabla \u || =0\})=0$ for a.e. $t >0$, and
that if $g$ is as above, then
\begin{equation}\label{coarea}
 \int _{\{|\nabla \u | >t\}} g(x)\, dx = \int _{\{|\nabla \u | >t\}\cap \{|\nabla |\nabla \u || =0\}} g(x)\, dx
+
 \int _{t}^{\infty} \int _{\{|\nabla \u | =\tau\}} \frac{g(x)}{|\nabla |\nabla \u ||}\, d\hh
(x)\,d\tau
\end{equation}
for $t \geq 0$.In particular, equation \eqref{coarea} entails that,
if $g \in L^1(\Omega )$, then
\begin{equation}\label{coareadiff}
  - \frac{d}{dt} \int _{\{|\nabla \u |
>t\}}g(x)  dx \geq  \int _{\{|\nabla \u | =t\}} \frac{g(x)}{|\nabla |\nabla \u ||}\, d\hh
(x)
\qquad \hbox{for a.e. $t>0$}.
\end{equation}

%

The following differential inequality involving integrals over the
level sets of a Sobolev function relies upon the coarea formula and
the relative isoperimetric inequality \eqref{relisop}, and is
established in \cite{Ma2}.

\begin{lemma}\label{talenti}
Let $\Omega $ be a domain with a Lipschitz boundary in $\rn$, $n
\geq 2$.
Let $v$ be a nonnegative function from $W^{1,2}(\Omega )$, and let
$\mu _v$ and $v^*$ denote the distribution function and the
decreasing rearrangement of $v$ defined as in \eqref{mu} and
\eqref{rearr}, respectively.
Then there exists a constant $C$, depending on the constant in
\eqref{isopboundary}, such that
\begin{equation}\label{tal2}
1 \leq C (-\mu _v'(t))^{1/2}\mu _v(t)^{-1/n'} \bigg(- \frac{d}{dt}
\int _{\{v>t \}}|\nabla v|^2 dx\bigg)^{1/2} \quad \hbox{for a.e. $t
\geq v^*(\m2 )$.}
\end{equation}
\end{lemma}

%
%

The next lemma provides us with a lower estimate for the scalar
product between the differential operator appearing on the left-hand
side of the equation in \eqref{eqdirichlet}, evaluated at some
 $\rN$-valued smooth function, and its Laplacian, via terms in
divergence form and a nonnegative term.

\begin{lemma}\label{differentialinequality}
Assume that  $a : (0, \infty ) \to (0, \infty )$ is of class  $C^1$,
and satisfies the
 first inequality in
\eqref{infsup}. Let $\Omega $ be an open set in $\rn$, $n \geq 2$,
and  let ${\bf  v} \in C^3(\Omega , \mathbb R ^N)$, with ${\bf v} =
(v^1, \dots , v^N)$.
Then
\begin{align}\label{1}
  \sum _{\alpha =1}^N \Delta v^\alpha   \,{\rm div} (a(|\nabla {\bf
v}|)\nabla v^\alpha ) & \geq \sum _{\alpha =1}^N{\rm div} (\Delta
v^\alpha \, a(|\nabla {\bf v}|)\nabla v^\alpha )
\\ \nonumber &-  \sum _{\alpha =1}^N \sum _{i,j}^n \big(v_{x_i x_j}^\alpha a(|\nabla
{\bf v}|)v_{x_i}^\alpha \big)_{x_j} + (1 + \min \{i_a , 0\})
a(|\nabla {\bf v}|) \sum _{\alpha =1}^N |\nabla ^2 v^\alpha|^2
\end{align}
in  $\{\nabla {\bf v} \neq 0\}$.
\end{lemma}
{\bf Proof}. In   $\{\nabla {\bf v} \neq 0\}$, we have that
\begin{align}\label{2}
\sum _{\alpha =1}^N \Delta v^\a  \,{\rm div} (a(|\nabla {\bf
v}|)\nabla v ^\a) & = \sum _{\alpha =1}^N {\rm div} (\Delta v^\a \,
a(|\nabla {\bf v}|)\nabla v^\a ) - \sum _{\alpha =1}^N \sum
_{i,j=1}^nv_{x_i x_j x_j}^\a a(|\nabla {\bf
v}|)v_{x_i}^\alpha \\
\nonumber & = \sum _{\alpha =1}^N {\rm div} (\Delta v^\a \,
a(|\nabla {\bf v}|)\nabla v^\a ) - \sum _{\alpha =1}^N \sum _{i,j}^n
\big(v_{x_i x_j}^\a a(|\nabla {\bf v}|)v_{x_i}^\a \big)_{x_j}
\\ \nonumber & \quad + \sum _{\alpha =1}^N \sum _{i,j=1}^n\big( v_{x_i x_j}^\a\big)^2 a(|\nabla {\bf
v}|) + \sum _{\alpha =1}^N \sum _{i,j=1}^nv_{x_i x_j}^\a a(|\nabla
{\bf v}|)_{x_j} v_{x_i}^\a\,.
\end{align}
Now,
\begin{align}\label{3}
\sum _{\alpha =1}^N \sum _{i,j=1}^n& \big( v_{x_i x_j}^\a\big)^2
a(|\nabla {\bf v}|) + \sum _{\alpha =1}^N \sum _{i,j=1}^n v_{x_i
x_j}^\a a(|\nabla {\bf v}|)_{x_j} v_{x_i}^\a \\ \nonumber & = \sum
_{i,j=1}^n \big( v_{x_i x_j}^\a \big)^2 a(|\nabla {\bf v}|) + \sum
_{\alpha , \b =1}^N\sum _{i,j,k=1}^nv_{x_i x_j}^\a a'(|\nabla {\bf
v}|) \frac{v_{x_k}^\b}{|\nabla \v|} v_{x_k x_j}^\b v_{x_i}^\a \\
&= a(|\nabla {\bf v}|) \bigg( \sum _{\alpha =1}^N \sum
_{i,j=1}^n(v_{x_i x_j}^\a)^2 + \sum _{\alpha , \b =1}^N\sum
_{i,j,k=1}^n\frac{a'(|\nabla {\bf v}|)|\nabla {\bf v}|}{a(|\nabla {\bf
v}|)} \frac{v_{x_i}^\a}{|\nabla {\bf v}|} v_{x_i x_j}^\a
\frac{v_{x_k}^\b}{|\nabla {\bf v}|} v_{x_k x_j}^\b\bigg).\nonumber
\end{align}
On setting $V^{j} = (v_{ x_1 x_j}^1, \dots , v_{ x_n x_j}^1, \dots ,
v_{ x_1 x_j}^N, \cdots , v_{ x_n x_j}^N)$ and $\omega = \frac{\nabla
{\bf v}}{|\nabla {\bf v}|}$, and making use of the first inequality
in \eqref{infsup}, one obtains  that
\begin{align}\label{4}
a(|\nabla {\bf v}|) & \bigg( \sum _{\alpha =1}^N\sum
_{i,j=1}^n(v_{x_i x_j}^\a)^2 + \sum _{\alpha , \b =1}^N\sum
_{i,j,k=1}^n\frac{a'(|\nabla {\bf v}|)|\nabla {\bf v}|}{a(|\nabla
{\bf v}|)} \frac{v_{x_i}^\a}{|\nabla {\bf v}|} v_{x_i x_j}^\a
\frac{v_{x_k}^\b}{|\nabla {\bf v}|} v_{x_k x_j}^\b \bigg)
\\ \nonumber & = a(|\nabla {\bf v}|)  \sum _{j=1}^n\bigg(|V^{j }|^2 +
\frac
{a'(|\nabla {\bf v}|) |\nabla {\bf v}|}{a(|\nabla {\bf v}|)} (V^{j} \cdot \omega )^2 \bigg)\\
\nonumber & \geq a(|\nabla {\bf v}|) \sum _{j=1}^n\big(|V^j|^2 + i_a
(V^j \cdot \omega )^2 \big)  \geq a(|\nabla {\bf v}|)  (1 + \min
\{i_a , 0\})\sum _{j}|V^j|^2.
\end{align}
Inequality \eqref{1} follows from \eqref{2}-\eqref{4}. \qed

The last two results of this section provide us with key
inequalities involving integrals on level sets and on level surfaces
of $\rN$-valued smooth functions in a smooth domain $\Omega$,
satisfying either Dirichlet, or Neumann homogenous boundary
conditions.

\begin{lemma}\label{keystep}
Let $\Omega$ be a domain  with
$\partial \Omega \in C^2$  in $\rn$, $n \geq 2$, and let $a$ be as in Theorem
\ref{dirichletc2}.  Assume that $\v\in C^\infty (\Omega
, \rN)\cap C^2(\overline \Omega , \rN)$, and $\v = 0$ on $\partial
\Omega$. Let $\mathcal B$ denote the second fundamental form on
$\partial \Omega$, and let ${\rm tr} \mathcal B$ be its trace.
 Then
\begin{multline}\label{37old}
\frac {(1 + \min \{i_a , 0\})^2}2
  b(t) \int _{\{|\nabla \v | =t\}}  |\nabla
|\nabla \v || \, d\hh (x)
 \leq t  \int _{\{|\nabla \v | =t\}}|{\rm {\bf div}}
(a(|\nabla {\bf v}|)\nabla \v)| d\hh (x) \\ +
 \frac {\|\nabla {\bf v} \|_{L^\infty (\Omega , \mathbb R^{Nn})}}{ b(t)}  \int _{\{|\nabla \v |>t\}}
|{\rm {\bf div}} (a(|\nabla {\bf v}|)\nabla \v)|^2  dx \\ +
a\big(\|\nabla {\bf v} \|_{L^\infty (\Omega , \mathbb
R^{Nn})}\big)\|\nabla {\bf v} \|_{L^\infty (\Omega , \mathbb
R^{Nn})}^2 \int_{\partial \Omega \cap
\partial\{|\nabla{\bf v} |>t\}}|{\rm tr} \mathcal B (x)| d\hh (x)
\end{multline}
for a.e. $t > 0$. Moreover, if $r >n-1$, then
\begin{multline}\label{37}
\frac {(1 + \min \{i_a , 0\})^2}2
  b(t) \int _{\{|\nabla \v | =t\}}  |\nabla
|\nabla \v || \, d\hh (x)
 \leq t  \int _{\{|\nabla \v | =t\}}|{\rm {\bf div}}
(a(|\nabla {\bf v}|)\nabla \v)| d\hh (x) \\ +
  \int _{\{|\nabla \v |>t\}}
\frac {1}{ a(|\nabla {\bf v}|)}|{\rm {\bf div}} (a(|\nabla {\bf
v}|)\nabla \v)|^2  dx +  a(t)t^2 \int_{\partial \Omega \cap
\partial\{|\nabla{\bf v} |>t\}}|{\rm tr} \mathcal B (x)| d\hh (x)
\end{multline}
for a.e. $t \geq t_\v$, where $t_\v = |\nabla \v|^*(\alpha
|\Omega|)$,  and $\alpha \in (0, \tfrac 12]$ is a constant depending
on $i_a$, $s_a$, $n$, $r$, $\|{\rm tr} \mathcal B \|_{L^r(\partial
\Omega )}$, $|\Omega |$, and on the constant  in  inequality
\eqref{isopboundary}.
%
%
\par\noindent
If $\Omega$ is convex, then the integral involving ${\rm tr}
\mathcal B$ can be dropped on the right-hand sides of
 inequalities \eqref{37old} and \eqref{37}, and the constant $\alpha$ neither depends on $r$, nor on  $\|{\rm tr} \mathcal B
\|_{L^r(\partial \Omega )}$.
  \end{lemma}
  \par \noindent {\bf Proof}.
The level set $\{|\nabla {\bf v} | >t\}$ is  open  for $t>0$.
Moreover,
  for a.e. $t>0$, the level surface $\partial\{|\nabla \v | >t\}$ is an
$(n-1)$-dimensional manifold of class $C^1$  outside a set of $\hh$
measure zero, and
$$\partial\{|\nabla {\bf v} | >t\} = \{|\nabla {\bf v} |
=t\} \cup \big( \partial \Omega \cap \partial   \{|\nabla {\bf v} |
>t\}\big).$$
By inequality \eqref{1} and the divergence theorem we have that
\begin{align}\label{21}
\sum _{\alpha = 1}^N &\int _{\{|\nabla {\bf v} | >t\}} \Delta
v^\alpha  \,{\rm div} (a(|\nabla {\bf v}|)\nabla v^\alpha )dx \geq
\sum _{\alpha = 1}^N\int _{\{|\nabla {\bf v} |
>t\}} {\rm div}
(\Delta v^\a \, a(|\nabla {\bf v}|)\nabla v^\a )dx
\\ \nonumber &-   \sum _{\alpha = 1}^N \int _{\{|\nabla {\bf v} |
>t\}}\sum _{i,j=1}^n\big( v_{x_i x_j}^\a a(|\nabla
{\bf v}|)v_{x_i}^\a\big)_{x_j}dx +  (1 + \min \{i_a , 0\}) \sum
_{\alpha = 1}^N \int _{\{|\nabla {\bf v} |
>t\}}a(|\nabla {\bf v}|) |\nabla ^2
v^\a|^2 dx
\\ \nonumber & = \sum _{\alpha = 1}^N
\int _{\partial \{|\nabla {\bf v} |
>t\}}
\Delta v^\a \, a(|\nabla {\bf v}|) \frac{\partial v^\a}{\partial \nu
} d\hh (x) - \sum _{\alpha = 1}^N \int _{\partial \{|\nabla {\bf v}
|
>t\}}\sum _{i,j=1}^n  v_{x_i x_j}^\a a(|\nabla
{\bf v}|)v_{x_i}^\a \nu _j d\hh (x) \\ \nonumber & \quad + (1 + \min
\{i_a , 0\}) \sum _{\alpha = 1}^N \int _{\{|\nabla {\bf v} |
>t\}}a(|\nabla {\bf v}|) |\nabla ^2
v^\a |^2dx \qquad \qquad \hbox{for a.e. $t>0$}.
\end{align}
Here, $\nu _j$ denotes the $j-th$ component of the outer normal
vector $\nu$ to $\partial \{|\nabla {\bf v} |
>t\}$. Now, observe that, for a.e. $t
>0$,  $$\nu = -\frac{\nabla |\nabla {\bf v}|}{|\nabla |\nabla
{\bf v}||} \quad \hbox{on } \,\,
 \{|\nabla {\bf v} |
=t\}.$$  Moreover,
$$\sum _{\alpha = 1}^N\sum _{i=1}^n  v_{x_i x_j}^\a v_{x_i}^\a = |\nabla {\bf v}|_{x_j}|\nabla {\bf v}|.$$
Thus,
\begin{align}\label{22}
\sum _{\alpha = 1}^N & \int _{\partial \{|\nabla {\bf v} |
>t\}}
\Delta v^\a \, a(|\nabla {\bf v}|) \frac{\partial v^\a}{\partial \nu
} d\hh (x) - \sum _{\alpha = 1}^N \int _{\partial \{|\nabla {\bf v}
|
>t\}}\sum _{i,j=1}^n  v_{x_i x_j}^\a a(|\nabla
{\bf v}|)v_{x_i}^\a \nu _j d\hh (x)
 \\ \nonumber & =
a(t)  \sum _{\alpha = 1}^N \int _{\{|\nabla {\bf v} | =t\}} \Delta
v^\a \, \frac{\partial v^\a}{\partial \nu } d\hh (x)+ a(t)t \int
_{\{|\nabla {\bf v} | =t\}} |\nabla |\nabla {\bf v} || \, d\hh (x) \\
\nonumber &
 \quad + \sum _{\alpha = 1}^N  \int_{\partial \Omega \cap \partial \{|\nabla
{\bf v} |
>t\}} a(|\nabla {\bf v}|) \Big(\Delta v^\a \frac{\partial v^\a}{\partial \nu
}- \sum _{i,j=1}^n  v_{x_i x_j}^\a v_{x_i}^\a \nu _j\Big) d\hh (x)
\quad \hbox{for a.e. $t>0$}.
\end{align}
Let us focus on the integrals on the right-hand side of \eqref{22}.
Since, for each $\alpha = 1, \dots , N$,
\begin{equation}\label{25}
{\rm div} (a(|\nabla {\bf v}|)\nabla v^\a ) = a(|\nabla {\bf
v}|)\Delta v^\a + a'(|\nabla {\bf v} |) \nabla v^\a \cdot \nabla
|\nabla {\bf v}|\,,
\end{equation}
one has that
\begin{align}\label{26}
& a(t)  \sum _{\alpha = 1}^N  \int _{\{|\nabla {\bf v} | =t\}}
\Delta v^\a \, \frac{\partial v^\a}{\partial \nu } d\hh
 \\
\nonumber & =  \sum _{\alpha = 1}^N \int _{ \{|\nabla {\bf v} |
=t\}} {\rm div} (a(|\nabla {\bf v}|)\nabla v^\a ) \, \frac{\partial
v^\a}{\partial \nu } d\hh (x) - a'(t) \sum _{\alpha = 1}^N \int _{
\{|\nabla {\bf v} | =t\}} \nabla v^\a \cdot \nabla |\nabla {\bf v}|
\frac{\partial v^\a}{\partial \nu } d\hh (x)
\\ \nonumber &=   \sum _{\alpha = 1}^N \int _{ \{|\nabla {\bf v} |
=t\}} {\rm div} (a(|\nabla {\bf v}|)\nabla v^\a ) \, \frac{\partial
v^\a}{\partial \nu } d\hh (x) + a'(t) \sum _{\alpha = 1}^N \int _{
\{|\nabla {\bf v} | =t\}} |\nabla |\nabla {\bf v}||
\bigg(\frac{\partial v^\a}{\partial \nu}\bigg)^2 d\hh (x)\,
\\ \nonumber &
\leq t\, \int _{ \{|\nabla {\bf v} | =t\}} {\rm {\bf div}}
(a(|\nabla {\bf v}|)\nabla {\bf v}) \,  d\hh (x) + a'(t) \int _{
\{|\nabla {\bf v} | =t\}} |\nabla |\nabla {\bf v}||
\bigg|\frac{\partial {\bf v}}{\partial \nu}\bigg|^2 d\hh (x)
\end{align}
for a.e. $t>0$. Here, we have exploited the fact that
$$ \frac{\partial v^\a}{\partial \nu} = - \frac{\nabla v^\a \cdot \nabla |\nabla {\bf
v}|}{|\nabla |\nabla {\bf v}||} \quad \quad  \hbox{on } \,\, 
\{|\nabla {\bf v} | =t\}\,, \quad \hbox{for a.e. $t>0$}.$$
\par\noindent
Next, we make use of the fact that, for each $\a = 1, \dots , N$,
\begin{multline}\label{23}
\Delta v^\a \frac{\partial v^\a}{\partial \nu } - \sum _{i,j=1}^n
v_{x_i x_j}^\a v_{x_i}^\a \nu _j \\ = {\rm div }_T
\bigg(\frac{\partial v^\a}{\partial \nu } \nabla _T v^\a\bigg) -
{\rm tr}\mathcal B \bigg(\frac{\partial v^\a}{\partial \nu }\bigg)^2
- \mathcal B (\nabla _T \,v^\a , \nabla _T \,v^\a) - 2 \nabla _T\,
v^\a \cdot \nabla _T\, \frac{\partial v^\a}{\partial \nu }
 \qquad \quad
\hbox{on $\partial \Omega$,}
\end{multline}
where 
 ${\rm div }_T$ and
$\nabla _T $ denote the divergence operator and the gradient
operator on $\partial \Omega$, respectively \cite[Equation
(3,1,1,2)]{Grisvard}.
Coupling \eqref{23} with the  condition $\v =0$ on $\partial \Omega$
tells us that
\begin{align}\label{24}
\Delta v^\a \frac{\partial v^\a}{\partial \nu }- \sum _{i,j=1}^N
v_{x_i x_j}^\a v_{x_i}^\a \nu _j = -{\rm tr}\mathcal B
\bigg(\frac{\partial v^\a}{\partial \nu }\bigg)^2
%
 \qquad \quad
\hbox{on $\partial \Omega $\,.}
\end{align}
Therefore,
\begin{multline}\label{29}
\sum _{\alpha = 1}^N\int_{
 \partial \Omega \cap\partial \{|\nabla{\bf v} |>t\}
 }
a(|\nabla {\bf v}|) \Big(\Delta v^\a \frac{\partial v^\a}{\partial
\nu }- \sum _{i,j=1}^n   v_{x_i x_j}^\a v_{x_i}^\a \nu _j\Big) d\hh
(x)
\\
=
- \sum _{\alpha = 1}^N\int_{
 \partial \Omega \cap\partial \{|\nabla{\bf v} |>t\}
 }
a(|\nabla {\bf v}|)
 {\rm tr}\mathcal B (x)\bigg(\frac{\partial v^\a}{\partial \nu }\bigg)^2
 d\hh
(x)
\\ \geq -  \int_{\partial \Omega \cap \partial \{|\nabla{\bf v} |>t\}} a(|\nabla \v|)|\nabla \v |^2 |{\rm tr}\mathcal B(x)| d\hh (x)
 \qquad \hbox{for a.e. $t>0$}.
\end{multline}
By Young's inequality and  the inequality  $|\Delta  \v| \leq
|\nabla ^2 \v|$, we have that
\begin{multline}\label{27}
\sum _{\alpha = 1}^N \int _{\{|\nabla {\bf v} | >t\}} \Delta
v^\alpha  \,{\rm div} (a(|\nabla {\bf v}|)\nabla v^\alpha )dx
  \leq \frac{1 + \min \{i_a , 0\}}2 \int _{\{|\nabla {\bf v} |
>t\}}a(|\nabla {\bf v}|) |\nabla ^2  \v|^2dx \\ + \frac 2{1 + \min \{i_a , 0\}} \int _{\{|\nabla {\bf v} |
>t\}} \frac 1{a(|\nabla {\bf v}|)} |{\rm {\bf div}} (a(|\nabla {\bf v}|)\nabla \v)|^2
dx
\end{multline}
for a.e. $t>0$. Combining \eqref{21}, \eqref{22}, \eqref{26},
\eqref{29} and \eqref{27} yields
\begin{multline}\label{30}
 \int _{\{|\nabla {\bf v} | =t\}}  |\nabla
|\nabla {\bf v} || \bigg(a(t)t + a'(t)\bigg|\frac{\partial {\bf
v}}{\partial \nu}\bigg|^2\bigg)\, d\hh (x) + \frac{1 + \min \{i_a ,
0\}}2 \int _{\{|\nabla {\bf v} |
>t\}}a(|\nabla {\bf v}|) |\nabla ^2
{\bf v}|^2dx \\
 \leq t  \int _{\{|\nabla {\bf v} | =t\}}|{\rm {\bf div}} (a(|\nabla {\bf v}|)\nabla \v )| d\hh
 (x) + \frac 2{1 + \min \{i_a , 0\}} \int _{\{|\nabla {\bf v} |
>t\}} \frac 1{a(|\nabla {\bf v}|)} |{\rm {\bf div}} (a(|\nabla {\bf v}|)\nabla \v|^2
dx
\\
+  \int_{\partial \Omega \cap \partial \{|\nabla {\bf v} |
>t\}} a(|\nabla {\bf v}|)|\nabla {\bf v} |^2 |{\rm tr}\mathcal B(x)|\,d\hh (x)
\end{multline}
 for a.e. $t>0$.
Observe  that
\begin{align}\label{31}
a(t) t +  a'(t)\bigg|\frac{\partial {\bf v}}{\partial \nu}\bigg|^2
\geq a(t) t + \min \{0, a'(t)\} t^2 \geq (1 + \min \{i_a , 0\})
a(t)t
 \quad   \hbox{on $\{|\nabla {\bf
 v}
| =t\}$}.
\end{align}
From
\eqref{30} and \eqref{31} we deduce that
\begin{align}\label{32new}
(1 + \min \{i_a , 0\}) b(t)& \int _{\{|\nabla \v | =t\}}  |\nabla
|\nabla \v || \, d\hh (x)  + \frac{1 + \min \{i_a , 0\}}2 \int
_{\{|\nabla \v |
>t\}}a(|\nabla \v|) |\nabla ^2
\v|^2 dx  \\ \nonumber &
 \leq t  \int _{\{|\nabla {\bf v} |
=t\}}|{\rm {\bf div}} (a(|\nabla {\bf v}|)\nabla \v)| d\hh
 (x)
 \\ \nonumber & \quad + \frac 2{1 + \min \{i_a , 0\}} \int _{\{|\nabla {\bf v} |
>t\}} \frac 1{a(|\nabla {\bf v}|)} |{\rm {\bf div}} (a(|\nabla {\bf v}|)\nabla \v)|^2
dx
\\ \nonumber & \quad +
 \int_{\partial \Omega \cap \partial\{|\nabla{\bf v} |>t\}}
a(|\nabla \v|)|\nabla \v |^2 \,|{\rm tr}\mathcal B(x)|  d\hh (x)
\qquad \hbox{for a.e. $t>0$}.
\end{align}
Hence, since $b(t)$ is an increasing function, and hence also
$a(t)t^2$ is an increasing function,
\begin{align}\label{32new'}
(1 + \min \{i_a , 0\}) b(t)& \int _{\{|\nabla \v | =t\}}  |\nabla
|\nabla \v || \, d\hh (x)
 \leq t  \int _{\{|\nabla {\bf v} |
=t\}}|{\rm div} (a(|\nabla {\bf v}|)\nabla \v)| d\hh
 (x)
 \\ \nonumber & \quad + \frac 2{1 + \min \{i_a , 0\}}
  \frac {\|\nabla {\bf v} \|_{L^\infty (\Omega , \mathbb R^{Nn})}}{ b(t)}  \int _{\{|\nabla {\bf v} |
>t\}}
 |{\rm div} (a(|\nabla {\bf v}|)\nabla \v)|^2
dx
\\ \nonumber & \quad + a\big(\|\nabla
{\bf v} \|_{L^\infty (\Omega , \mathbb R^{Nn})}\big)\|\nabla {\bf v}
\|_{L^\infty (\Omega , \mathbb R^{Nn})}^2
 \int_{\partial \Omega \cap \partial\{|\nabla{\bf v} |>t\}}
 |{\rm tr}\mathcal B(x)|  d\hh (x)
\end{align}
for a.e. $t>0$. Inequality \eqref{37old} follows.
\par
Let us next focus on \eqref{37}. Inasmuch as  $a(t)t^2$ is an
increasing function,
\begin{align}\label{33}
\int_{\partial \Omega \cap \partial\{|\nabla{\bf v} |>t\}}&
a(|\nabla \v|)|\nabla \v |^2 |{\rm tr}\mathcal B(x)| d\hh (x)
\\ \nonumber & \leq 2 \int_{\partial \Omega \cap \partial\{|\nabla{\bf v} |>t\}}
\big(a(|\nabla \v|)^{1/2}|\nabla \v |- a(t)^{1/2}t\big)^2 |{\rm
tr}\mathcal B(x)| d\hh (x) \\ \nonumber & \quad  + 2 a(t)t^2
\int_{\partial \Omega \cap
\partial\{|\nabla{\bf v} |>t\}} |{\rm tr}\mathcal B(x)| d\hh (x)\,.
\end{align}
Denote, for simplicity, the distribution function $\mu _{|\nabla
\v|}$ of $|\nabla \v|$ by $\mu : [0 , \infty ) \to [0, \mo ]$ \quad \hbox{for a.e. $t>0$}.
Set
 $\delta = \frac{n-1}{nr'} -
\frac{n-2}n$, and observe that $\delta >0$ since $r>n-1$. Thanks to
our assumptions on the function $a$ and to the chain rule for
vector-valued Sobolev functions \cite[Theorem 2.1]{MarcusMizel}, the
function
 $\max\{a(|\nabla \v|)^{1/2}|\nabla \v |-
a(t)^{1/2}t , 0\}$ belongs to $W^{1,2}(\Omega )$. H\"older's
inequality and an application of Lemma \ref{trace} with $v$ replaced
with $\max\{a(|\nabla \v|)^{1/2}|\nabla \v |- a(t)^{1/2}t , 0\}$
tell us that
\begin{align}\label{34}
 & \int_{\partial \Omega \cap \partial\{|\nabla{\bf v} |>t\}}
\big(a(|\nabla \v|)^{1/2}|\nabla \v |- a(t)^{1/2}t\big)^2 |{\rm
tr}\mathcal B(x)| d\hh (x) \\ \nonumber & \leq \bigg(\int_{\partial
\Omega \cap
\partial\{|\nabla{\bf v} |>t\}} \big(a(|\nabla
\v|)^{1/2}|\nabla \v |- a(t)^{1/2}t\big)^{2r'}  d\hh
(x)\bigg)^{\frac 1{r'}} \bigg( \int_{\partial \Omega \cap
\partial\{|\nabla{\bf v} |>t\}}|{\rm tr}\mathcal B(x)|^{r} d\hh
(x)\bigg)^{\frac 1{r}}
\\ \nonumber &
\leq C \mu (t)^\delta \|{\rm tr}\mathcal B\|_{L^{r}(\partial \Omega
)} \int _{\{|\nabla \v|
>t\}}\left|\nabla \big[a(|\nabla \v|)^{1/2}|\nabla \v
|\big]\right|^2dx
\\ \nonumber & = C\mu (t)^
\delta
\|{\rm tr}\mathcal B\|_{L^{r}(\partial \Omega )} \int _{\{|\nabla \v
|>t\}}\Big(\tfrac 12 a'(|\nabla \v|)a(|\nabla \v|)^{-1/2}|\nabla \v
| + a(|\nabla \v|)^{1/2}\Big)^2|\nabla |\nabla \v||^2dx
\\ \nonumber & \leq  C'\mu (t)^\delta
\|{\rm tr}\mathcal B\|_{L^{r}(\partial \Omega )} \int _{\{|\nabla \v
|
>t\}}a(|\nabla \v|)|\nabla ^2 \v|^2dx \quad \hbox{if $t >|\nabla
\v|^*(\m2 )$} \quad \hbox{for a.e. $t>0$},
\end{align}
for some positive constants $C$, depending on the constant in
\eqref{isopboundary} and on $r$, and $C'$  depending on the same
quantities and on $s_a$.
%
%
%
 Observe that, in  the last
inequality in \eqref{34}, we have employed the inequality
%
 $|\nabla |\nabla \v|| \leq |\nabla ^2 \v|$.
Set
$$\beta =\Big( \frac{1 + \min \{i_a , 0\}}{4C'\|{\rm tr}\mathcal B\|_{L^{r}(\partial
\Omega )}}\Big)^{\frac 1\delta},$$ where $C'$ is the constant
appearing in \eqref{34},
 \begin{equation}\label{alpha}
 \alpha = \min \{\beta /\mo , 1/2\},
 \end{equation}
  and $t_{\v} =|\nabla \v |^*(\alpha \mo )$. Thus,
$\alpha$ depends on the quantities specified in the statement, and
%
\begin{equation}\label{35}
\frac{1 + \min \{i_a , 0\}}2 - 2C'\|{\rm tr}\mathcal
B\|_{L^{r}(\partial \Omega )} \mu (t)^{\delta} \geq 0 \quad \quad
\hbox{if $t > t_\v$.}
\end{equation}
Inequality \eqref{37} follows from \eqref{32new}, \eqref{33},
\eqref{34} and \eqref{35}, since $\tfrac{1 + \min \{i_a , 0\}}2 <1$.
\par
The assertion concerning the case when $\Omega$ is convex follows
via the same argument, on observing that the
%
leftmost  side of \eqref{29} can be estimated from below just by
$0$, inasmuch as ${\rm tr}\mathcal B \leq 0$ on  $\partial \Omega$
in this case.
%
\qed

\begin{lemma}\label{keystepneumann}
Let $\Omega$ be a domain  with
$\partial \Omega \in C^2$ in $\rn$, $n \geq 2$, and let $a$ be as in Theorem
\ref{dirichletc2}.  Assume that $\v\in C^\infty (\Omega
, \rN)\cap C^2(\overline \Omega , \rN)$, and $\frac{\partial
\v}{\partial \nu} = 0$ on $\partial \Omega$. Let $\mathcal B$ denote
the second fundamental form on $\partial \Omega$, and let $|\mathcal
B|$ be its operator norm, namely
$$|\mathcal
B (x)| = \sup_{0 \neq \zeta \in \R^{n-1}} \frac{|\mathcal B(x) (\zeta ,
\zeta )|}{|\zeta |^2} \quad \hbox{for $x \in \partial \Omega $.}$$
 Then
\begin{multline}\label{37oldneumann}
\frac {(1 + \min \{i_a , 0\})^2}2
  b(t) \int _{\{|\nabla \v | =t\}}  |\nabla
|\nabla \v || \, d\hh (x)
 \leq t  \int _{\{|\nabla \v | =t\}}|{\rm {\bf div}}
(a(|\nabla {\bf v}|)\nabla \v)| d\hh (x) \\ +
 \frac {\|\nabla {\bf v} \|_{L^\infty (\Omega , \mathbb R^{Nn})}}{ b(t)}  \int _{\{|\nabla \v |>t\}}
|{\rm {\bf div}} (a(|\nabla {\bf v}|)\nabla \v)|^2  dx \\ +
a\big(\|\nabla {\bf v} \|_{L^\infty (\Omega , \mathbb
R^{Nn})}\big)\|\nabla {\bf v} \|_{L^\infty (\Omega , \mathbb
R^{Nn})}^2 \int_{\partial \Omega \cap
\partial\{|\nabla{\bf v} |>t\}}|\mathcal B (x)| d\hh (x)
\end{multline}
for a.e. $t > 0$. Moreover, if  $r >n-1$, then
\begin{multline}\label{37neumann}
\frac {(1 + \min \{i_a , 0\})^2}2
  b(t) \int _{\{|\nabla \v | =t\}}  |\nabla
|\nabla \v || \, d\hh (x)
 \leq t  \int _{\{|\nabla \v | =t\}}|{\rm {\bf div}}
(a(|\nabla {\bf v}|)\nabla \v)| d\hh (x) \\ +
  \int _{\{|\nabla \v |>t\}}
\frac {1}{ a(|\nabla {\bf v}|)}|{\rm {\bf div}} (a(|\nabla {\bf
v}|)\nabla \v)|^2  dx +  a(t)t^2 \int_{\partial \Omega \cap
\partial\{|\nabla{\bf v} |>t\}}| \mathcal B (x)| d\hh (x)
\end{multline}
for a.e. $t \geq t_\v$, where $t_\v = |\nabla \v|^*(\alpha
|\Omega|)$, and
 $\alpha \in (0, \tfrac 12]$ is a constant depending
on $i_a$, $s_a$, $n$, $r$, $\|\, |\mathcal B |\,\|_{L^r(\partial
\Omega )}$, $|\Omega |$, and on the constant  in  inequality
\eqref{isopboundary}.
%
%
\par\noindent
If $\Omega$ is convex, the integral involving $\mathcal B$ can be
dropped on the right-hand sides of
 inequalities \eqref{37oldneumann} and \eqref{37neumann}, and the constant $\alpha$ neither depends on $r$, nor on  $\|\,|\mathcal
 B|\,
\|_{L^r(\partial \Omega )}$.
  \end{lemma}
\par\noindent
{\bf Sketch of the proof}. The proof is completely analogous to that
of Lemma \ref{keystep}. One has just to observe that, by \eqref{23}
and the  condition $\frac {\partial {\bf v}}{\partial \nu }
=0$ on $\partial \Omega$, equation \eqref{24} has to be replaced
 with
\begin{align}\label{24neumann}
\Delta v^\a \frac{\partial v^\a}{\partial \nu }- \sum _{i,j=1}^n
v_{x_i x_j}^\a v_{x_i}^\a \nu _j = - \mathcal B (\nabla _T\, v^\a ,
\nabla _T\, v^\a)
 \qquad \quad
\hbox{on  $\partial \Omega $\,,}
\end{align}
which, in particular, implies that
\begin{align}\label{24neumann'}
\bigg|\Delta v^\a \frac{\partial v^\a}{\partial \nu }- \sum
_{i,j=1}^N v_{x_i x_j}^\a v_{x_i}^\a \nu _j\bigg| \leq |\mathcal B|
|\nabla v^\a|^2
 \qquad \quad
\hbox{on $ \partial \Omega $\,.}
\end{align}
\par\noindent
The conclusion concerning convex domains $\Omega$ holds owing to the
fact that $\mathcal B \leq 0$
on $\partial \Omega$ in this case. \qed

\section{Proof of the main results}\label{sec4}

 Let $B$ be the Young function defined by \eqref{B}, and let $B_n$ be its Sobolev conjugate given by
\eqref{sobolevconj}.  Assume that $f \in L^{\widetilde {B_n}}(\Omega , \rN
)$.
A weak solution to the Dirichlet problem \eqref{eqdirichlet} is a function $u \in
W^{1,B}_0(\Omega , \rN )$ such that
\begin{equation}\label{weakdirichlet}
\int _\Omega a(|\nabla \u|) \nabla \u \cdot \nabla \phi \,dx = \int
_\Omega {\bf f} \cdot \phi \, dx
\end{equation}
for every $\phi \in W^{1,B}_0(\Omega , \rN)$.
\par\noindent
 Assume now, in addition that $\Omega $ has a Lipschitz boundary.
A weak solution to the Neumann  problem \eqref{eqneumann} is a function $\u \in
W^{1,B}(\Omega , \rN)$ such that
\begin{equation}\label{weakneumann}
\int _\Omega a(|\nabla \u|) \nabla \u \cdot \nabla \phi \,dx = \int
_\Omega {\bf f} \cdot \phi \,dx
\end{equation}
for every $\phi \in W^{1,B}(\Omega , \rN)$.
\par\noindent

Note that the left-hand sides of \eqref{weakdirichlet} and
\eqref{weakneumann} are well defined by inequalities
\eqref{holder} and \eqref{conj}. The right-hand sides  are also
well defined, owing to the Sobolev inequality \eqref{embedding0} and
 inequality   \eqref{holder} with $B$ replaced
with $B_n$. In particular, the right-hand sides of
\eqref{weakdirichlet} and \eqref{weakneumann} are well defined
if ${\bf f} \in L^{n, 1}(\Omega , \mathbb R^{N})$,
since $L^{n, 1}(\Omega ,\rN ) \to L^{\widetilde {B_n}}(\Omega ,
\rN)$, as shown in \cite[Remark 2.12]{CMlipschitz}.
\par
The following existence and uniqueness result holds for weak
solutions to problems  \eqref{eqdirichlet} and \eqref{eqneumann}.

\begin{theorem}\label{esistenza}
Let $\Omega $ be a domain in $\rn$, $n \geq 2$, and let $N\geq 1$.
Assume that $a : (0, \infty ) \to (0, \infty )$ is of class $C^1$,
and fulfils
\eqref{infsup}. Let ${\bf f} \in  L^{n, 1}(\Omega , \mathbb R^{N})$.
Then there exists a unique solution $\u \in W^{1,B}_0(\Omega , \rN )$
to problem \eqref{eqdirichlet}.
\par\noindent
Assume, in addition, that $\Omega$ has a Lipschitz boundary, and ${\bf
f}$ fulfills \eqref{intf0}.
Then there exists a  solution $\u \in W^{1, B} (\Omega , \rN )$ to
problem \eqref{eqneumann}, which is  unique up to additive  constant
vectors in $\rN$. In particular, there exists a unique solution in $
W^{1, B}_\bot (\Omega , \rN)$.
\end{theorem}

A proof of Theorem  \ref{esistenza} in the case when $N=1$ can be
found in \cite{CMlipschitz}; the proof for $N>1$ is completely
analogous.
\par
The next Proposition provides us with a basic energy estimate for
the weak  solutions to problems  \eqref{eqdirichlet} and
\eqref{eqneumann}.

\begin{proposition}\label{energia}
Let $\Omega $ be a domain in $\rn$, $n \geq 2$, and
let $N \geq 1$. Assume that $a$ is as in Theorem \ref{esistenza}.
%
Let ${\bf f} \in  L^{n, 1}(\Omega ,
\mathbb R^{Nn})$.
 \par\noindent (i) Let $\u \in W^{1,B}_0(\Omega ,
\rN)$ be the weak solution to problem \eqref{eqdirichlet}. Then
\begin{equation}\label{energia1}
\int _\Omega B(|\nabla \u|) dx \leq C \,  \|{\bf
f}\|_{L^{n,1}(\Omega , \rN)}\, b^{-1}\big(\|{\bf f}\|_{L^{n,1}(\Omega
, \rN )}\big),
\end{equation}
where $C=C'|\Omega|$, and $C'$ is a constant depending on   $n$,  $N$ and $i_a$.
\par\noindent (ii) Assume, in addition, that $\Omega$ has a Lipschitz boundary, and
${\bf f}$ fulfills \eqref{intf0}. Let $\u \in W^{1,B}(\Omega , \rN)$
be a weak solution to problem \eqref{eqneumann}. Then inequality
holds for some constant  $C$ depending on $n$, $N$, $i_a$ and on the
constant in \eqref{relisop}
\end{proposition}
{\bf Proof}. (i) Making use of $\u$ as test function $\phi$ in the
definition of weak solution \eqref{weakdirichlet} tells us that
$$\int _\o a(|\nabla \u|) |\nabla \u|^2 \, dx = \int _\o
{\bf f} \cdot \u \, dx.$$
 By the first inequality in \eqref{B5},  H\"older's inequality in
 Lorentz spaces \eqref{holderlorentz}, and \eqref{lorentz11}, there
 exist constants $C$ and $C'$, depending on $n$, such that
\begin{equation}\label{energy1}
\int _\o B(|\nabla \u|) \, dx \leq C \|{\bf f}\|_{L^{n,1}(\Omega ,
\rN)} \|\u\|_{L^{n', \infty}(\Omega , \rN)} \leq  C' \|{\bf
f}\|_{L^{n,1}(\Omega , \rN)} \|\u\|_{L^{n'}(\Omega , \rN)}.
\end{equation}
By the Poincar\'e inequality in $W^{1,1}_0(\Omega , \rN)$, there
exists a constant $C=C(n, N)$ such that
\begin{equation}\label{energy2}
 \|\u\|_{L^{n'}(\Omega , \rN)} \leq C \int _\o |\nabla \u|  \, dx.
\end{equation}
On the other hand, Jensen's inequality entails that
\begin{equation}\label{energy3}
B\bigg(\frac 1{\mo }\int _\o |\nabla \u|  \, dx\bigg) \leq \frac
1{\mo }\int _\o B(|\nabla \u|) \, dx.
\end{equation}
Combining inequalities \eqref{energy1}--\eqref{energy3}, and making
use of the second inequality in \eqref{youngprop} yields
\begin{equation}\label{energy4}
\frac 1{\mo }\int _\o B(|\nabla \u|) \, dx \leq \widetilde B \big( 2
C \|{\bf f}\|_{L^{n,1}(\Omega , \rN)}\big).
\end{equation}
Since $b^{-1}$ is an increasing function, equation \eqref{Btilde}
ensures that $\widetilde B (t) \leq t b^{-1}(t)$ for $t \geq 0$.
Thus, \eqref{energia1} follows from \eqref{energy4}, via
\eqref{B7bis}.
\par\noindent
(ii) The proof follows along the same lines as above. One has just
to make use of the fact that inequality \eqref{energy2} holds, for
every $\u \in W^{1,B}_\bot (\Omega , \rN)$, with a constant $C$
depending on $n$, $N$  and on the constant  in \eqref{relisop}
\cite[Theorem 5.2.3]{Mazlibro}, and that any solution $u$ to
\eqref{eqneumann} differs from the solution in $W^{1,B}_\bot (\Omega
, \rN)$ by a constant vector in $\rN$.
 \qed

\medskip
\par\noindent
 We are now in a position to
prove Theorem \ref{dirichletc2}.
\smallskip
\par\noindent {\bf
Proof of Theorem\ref{dirichletc2}}. We split the proof in steps.
\par\noindent
{\bf Step 1}. We assume in addition, for the time being,  that
\begin{equation}\label{domainsmooth}
\partial \Omega \in C^\infty ,
\end{equation}
and there exist positive constants $c$ and $C$ such that
\begin{equation}\label{nondeg}
c \leq a(t) \leq C \quad \hbox{for $t \geq 0$\,.}
\end{equation}
\par\noindent Since  ${\bf f} \in L^{n,1} (\Omega , \rN)$,
in particular ${\bf f} \in L^2 (\Omega , \rN)$,  owing to
\eqref{lorentz12}. A result by Elcrat and Meyers implies that the
weak solution ${\bf u}$ to problem \eqref{eqdirichlet} belongs to
$W^{2,2}(\Omega )$ \cite[Theorem 8.2]{BF}. Notice that the
hypotheses of that result are fulfilled under our additional
assumptions \eqref{domainsmooth}--\eqref{nondeg}, owing to equation
\eqref{elliptic1}. Thus, ${\bf u} \in W^{1,2}_0(\Omega , \rN) \cap
W^{2,2}(\Omega , \rN)$. By  standard approximation,  there exists a
sequence $\{{\bf u}_k\} \subset C^\infty (\Omega , \rN)\cap
C^2(\overline \Omega , \rN)$ such that $\u_k = 0$ on $\partial
\Omega$,
\begin{equation}\label{convuk}
{\bf u}_k \to {\bf u} \quad \hbox{in $W^{1,2}_0(\Omega ,\rN)$,}
\quad {\bf u}_k \to {\bf u} \quad \hbox{in $W^{2,2}(\Omega , \rN)$,}
\quad \nabla {\bf  u}_k \to \nabla {\bf u}\quad  \hbox{a.e. in
$\Omega $},
\end{equation}
as $k \to \infty$. Furthermore,  $|\nabla {\bf u}_k| \in
W^{1,2}(\Omega )$ and $|\nabla |\nabla {\bf u}_k|| \leq |\nabla ^2
\u_k|$ a.e. in $\Omega$,
%
%
 by the chain rule for
vector-valued Sobolev functions \cite[Theorem 2.1]{MarcusMizel}. Thus, owing
to the compactness of the trace embedding ${\rm Tr} : W^{1,2}(\Omega
) \to L^1(\partial \Omega )$, we may also assume that
\begin{equation}\label{convtr}
{\rm Tr}\,|\nabla {\bf u}_k| \to {\rm
Tr}\, |\nabla {\bf u}|
\qquad \hbox{$\hh$ a.e. on $\partial \Omega $,}
\end{equation}
as $k \to \infty$.
 We
claim that
\begin{equation}\label{operator}
- {\rm {\bf div}} (a(|\nabla {\bf u_k}|)\nabla {\bf u_k} ) \to {\bf
f}
\quad \hbox{in $L^2(\Omega , \rN)$},
\end{equation}
as $k \to \infty$. Let us verify this claim. First, since $\nabla
{\bf u} \in W^{1,2}(\Omega , \mathbb R^{Nn})$, an application of the
chain rule for vector-valued Sobolev functions again
tells us that, for each  $\alpha = 1, \dots N$,
\begin{equation}\label{chain}
{\rm {div}} (a(|\nabla {\bf u}|)\nabla { u}^\alpha ) =
 \frac{a'(|\nabla  {\bf u}|)}{|\nabla  {\bf u}|} \sum _{\beta =1}^N \sum_{i,j=1}^n  u^\alpha _{x_i} u^\beta_{x_j}  u^\beta_{x_j x_i}
  \chi_{\{\nabla  {\bf u} \neq 0\}} + a(|\nabla  {\bf u}|)\Delta  u^\alpha
   \quad \hbox{a.e. in $\Omega $,}
 \end{equation}
 and that the same equation holds with $\u$ replaced with $\u_k$.
 Here, and in what follows, we adhere the convention that $0 \cdot \infty
 = 0$, so that $\frac{\chi_{\{\nabla  {\bf u} \neq 0\}}}{|\nabla
 \u|} =0$ in $\{\nabla  {\bf u} = 0\}$.
 Now, for each $k \in \N$ and $\alpha = 1, \dots N$,
 \begin{align}\label{chain1}
\bigg(&  \int_\Omega \big| {\rm { div}} (a(|\nabla { \u_k}|)\nabla {
u_k^\alpha } ) + f^\alpha (x)\big|^2 \, dx \bigg)^{\frac 12} =
\bigg( \int_\Omega  \big| {\rm { div}} (a(|\nabla { \u_k}|)\nabla {
u_k^\alpha } ) - {\rm { div}} (a(|\nabla {\bf u}|)\nabla { u^\alpha}
)\big|^2 \, dx \bigg)^{\frac 12}   \\ \nonumber &  =
 \bigg( \int_\Omega \bigg|
 \frac{a'(|\nabla  {\bf u_k}|)}{|\nabla  {\bf u_k}|} \sum _{\beta =1}^N \sum_{i,j=1}^n  {u^\alpha _k} _{x_i} {u^\beta _k}_{x_j}  {u^\beta_k}_{x_j x_i}
  \chi_{\{\nabla  {\bf u}_k \neq 0\}} + a(|\nabla  {\bf u}_k|)\Delta
  u^\alpha_k
\\ \nonumber &
\quad -
 \frac{a'(|\nabla  {\bf u}|)}{|\nabla  {\bf u}|} \sum _{\beta =1}^N \sum_{i,j=1}^n  u^\alpha _{x_i} u^\beta_{x_j}  u^\beta_{x_j x_i}
  \chi_{\{\nabla  {\bf u} \neq 0\}} - a(|\nabla  {\bf u}|)\Delta  u^\alpha \bigg|^2\, dx \bigg)^{\frac 12}
\\ \nonumber  & \leq
\bigg(
 \int_\Omega \big| a(|\nabla {\bf u_k} |)(\Delta u^\alpha  _k - \Delta u^\alpha) \big|^2\, dx  \bigg)^{\frac 12}
 +
 \bigg( \int_\Omega \big| (a(|\nabla {\bf u_k}|)- a(|\nabla {\bf u}|))  \Delta u^\alpha  \big|^2\, dx  \bigg)^{\frac 12}
\\ \nonumber  & \quad +
\bigg( \int_\Omega   \bigg|\frac{a'(|\nabla  {\bf u_k}|)}{|\nabla
{\bf u_k}|} \sum _{\beta =1}^N \sum_{i,j=1}^n  {u^\alpha _k} _{x_i}
{u^\beta _k}_{x_j}
  \chi_{\{\nabla  {\bf u}_k \neq 0\}} ({ { u}_k^\beta }_{x_j x_i} -  u^\beta_{x_j x_i} )\bigg|^2\, dx  \bigg)^{\frac 12}
 \\ \nonumber
& \quad +
 \bigg( \int_\Omega   \bigg|\sum _{\beta =1}^N \bigg(\frac{a'(|\nabla  {\bf u_k}|)}{|\nabla
{\bf u_k}|}  \sum_{i,j=1}^n  {u^\alpha _k} _{x_i} {u^\beta _k}_{x_j}
  \chi_{\{\nabla  {\bf u}_k \neq 0\}} -
  \frac{a'(|\nabla  {\bf u}|)}{|\nabla  {\bf u}|}  \sum_{i,j=1}^n  u^\alpha _{x_i} u^\beta_{x_j}
  \chi_{\{\nabla  {\bf u} \neq 0\}}\bigg)u^\beta _{x_j x_i} \bigg|^2\, dx  \bigg)^{\frac 12}.
 \end{align}
 Since the functions $a(t)$ and $a'(t)t$ are  bounded,
  the first and the third addend on the rightmost side of \eqref{chain1}
   converge to $0$ as $k\to \infty$, inasmuch as ${ { u}_k^\beta }_{x_j x_i} \to u^\beta_{x_j x_i}$
   in $L^2(\Omega)$ as $k \to \infty$, for $\beta =1, \dots N$
   and $i, j =1, \dots n$.
%
%
The boundedness of the functions $a(t)$ and $a'(t)t$, and the
   convergence of $\nabla \u_k$ to $\nabla \u$ a.e. in $\Omega$
   implies that the second and the fourth addend also converge to
   $0$  by the dominated convergence theorem for integrals. Hence,
   \eqref{operator} follows.

\smallskip
\par\noindent
{\bf Step 2}. Let $\{\u _k\}$ be the sequence considered in Step 1.
For each $k \in \N$, the function  $\u _k$ satisfies the same
assumptions as the function $\v$ in Lemma \ref{keystep}. Hence,
inequality \eqref{37} holds with $\v$ replaced with $\u_k$. This
tells us that
\begin{align}\label{32k}
C b(t)\int _{\{|\nabla \u _k | =t\}}  |\nabla |\nabla \u_k || \,
d\hh (x) &  \leq t  \int _{\{|\nabla {\bf u}_k | =t\}}|{\rm {\bf
div}} (a(|\nabla {\bf u}_k|)\nabla \u_k)| d\hh
 (x)
 \\ \nonumber & \quad +  \int _{\{|\nabla {\bf u}_k |
>t\}} \frac 1{a(|\nabla {\bf u}_k|)} |{\rm {\bf div}} (a(|\nabla {\bf u}_k|)\nabla
\u_k)|^2 \\ \nonumber & \quad +   a(t)t^2 \int_{\partial \Omega \cap
\partial\{|\nabla{\bf u}_k |>t\}}|{\rm tr} \mathcal B (x)| d\hh (x)  \quad \hbox{for a.e. $t> t_{\u
_k}$},
\end{align}
where $C = \tfrac {(1 + \min \{i_a , 0\})^2}2$, and $t_{\u _k}$ is
defined analogously to $t_\v$, with $\v$ replaced with $\u _k$.
 We claim that inequality \eqref{32k}
continues to hold with $\u_k$ replaced with $\u$, namely that
\begin{align}\label{32u}
 C b(t)\int _{\{|\nabla \u  | =t\}}  |\nabla
|\nabla \u  || \, d\hh (x) &  \leq t  \int _{\{|\nabla {\bf u} |
=t\}}|{\bf f}(x)| d\hh
 (x)
  + \int _{\{|\nabla {\bf u} |
>t\}} \frac 1{a(|\nabla {\bf u}|)} |{\bf f}(x)|^2
dx \\ \nonumber & \quad +   a(t)t^2 \int_{\partial \Omega \cap
\partial ^M \{|\nabla{\bf u} |>t\}}|{\rm tr} \mathcal B (x)| d\hh (x)  \qquad \hbox{for a.e. $t> t_{\u}$}.
\end{align}
To verify this claim, observe that $ t_{\u _k} \to t_{\u}$ as $k \to
\infty$, fix any $t>t_{\u}$ and $h>0$, and, for sufficiently large
$k$, integrate inequality \eqref{32k} over the interval $(t, t+h)$,
and make use of the coarea formula \eqref{coareanabla} to obtain
\begin{align}\label{32h}
 C \int _{\{t < |\nabla \u _k | <{t+h}\}} b(|\nabla \u _k |) |\nabla
|\nabla \u_k ||^2 \, dx &  \leq \int _{\{t < |\nabla \u _k |
<{t+h}\}} |\nabla \u _k ||\nabla |\nabla \u_k || \,|{\rm {\bf div}}
(a(|\nabla {\bf u}_k|)\nabla \u_k)|\, dx
 \\ \nonumber & \quad +  \int _t^{t+h} \int _{\{|\nabla {\bf u}_k |
>\tau\}} \frac 1{a(|\nabla {\bf u}_k|)} |{\rm {\bf div}} (a(|\nabla {\bf u}_k|)\nabla \u_k)|^2
dx\, d\tau
 \\ \nonumber & \quad +  \int _t^{t+h} a(\tau)\tau ^2 \int_{\partial \Omega
\cap
\partial\{|\nabla{\bf u}_k |>\tau \}}|{\rm tr} \mathcal B (x)|d\hh (x)\,d\tau\, .
\end{align}
We have that
\begin{align}\label{conv1}
 & \int _{\{t < |\nabla \u _k | <{t+h}\}} b(|\nabla \u _k |) |\nabla
|\nabla \u_k ||^2 \, dx  =
 \int _\Omega \chi_{_{\{t < |\nabla \u _k | <{t+h}\}}}(x)b(|\nabla \u _k |)\sum _{i=1}^n\big(
 \sum _{\beta =1}^N
\sum_{j=1}^n  \frac{{u^\beta_k}_{x_j}}{|\nabla \u _k |}
{u^\beta_k}_{x_j x_i} \big)^2\, dx
\\ \nonumber & =
\int _\Omega \chi_{_{\{t < |\nabla \u _k | <{t+h}\}}}(x) b(|\nabla
\u _k |)\sum _{i=1}^n\Big(
 \sum _{\beta =1}^N
\sum_{j=1}^n  \frac{{u^\beta_k}_{x_j}}{|\nabla \u _k |}
({u^\beta_k}_{x_j x_i} - u^\beta _{x_j x_i } ) + \sum _{\beta =1}^N
\sum_{j=1}^n \frac{{u^\beta_k}_{x_j}}{|\nabla \u _k |}u^\beta _{x_j
x_i }\Big)^2\, dx
\\ \nonumber
 & =
\int _\Omega \chi_{_{\{t < |\nabla \u _k | <{t+h}\}}}(x) b(|\nabla
\u _k |)\sum _{i=1}^n\Big(
 \sum _{\beta =1}^N
\sum_{j=1}^n  \frac{{u^\beta_k}_{x_j}}{|\nabla \u _k |}
({u^\beta_k}_{x_j x_i} - u^\beta _{x_j x_i } )\Big)^2\,dx \\
\nonumber & + 2\int _\Omega \chi_{_{\{t < |\nabla \u _k |
<{t+h}\}}}(x) b(|\nabla \u _k |)\sum _{i=1}^n \Big(\sum _{\beta
=1}^N \sum_{j=1}^n \frac{{u^\beta_k}_{x_j}}{|\nabla \u _k |}
({u^\beta_k}_{x_j x_i} - u^\beta _{x_j x_i } )\Big) \Big(\sum
_{\beta =1}^N \sum_{j=1}^n \frac{{u^\beta_k}_{x_j}}{|\nabla \u _k |}
  u^\beta _{x_j x_i }\Big) \,dx
\\ \nonumber &
 + \int _\Omega \chi_{_{\{t < |\nabla \u _k | <{t+h}\}}}(x) b(|\nabla \u _k |)\sum
_{i=1}^n\Big(\sum _{\beta =1}^N \sum_{j=1}^n
\frac{{u^\beta_k}_{x_j}}{|\nabla \u _k |}u^\beta _{x_j x_i
}\Big)^2\, dx.
\end{align}
Note that the first equality holds by the chain rule for
vector-valued functions. Since $b$ is an increasing function,
$$\chi_{_{\{t < |\nabla \u _k | <{t+h}\}}}(x)b(|\nabla \u _k |)\leq b(t+h) \quad \hbox{for $x \in \Omega$,} $$ for every
$k \in \N$. Moreover, ${|{u^\beta_k}_{x_j}|}/{|\nabla \u _k |} \leq
1 $ for every $k \in \N$, $\beta = 1, \dots , N$, $j=1, \dots , n$.
Thus, the first integral on the rightmost side of \eqref{conv1}
converges to $0$ as $k \to \infty$, since ${ { u}_k^\beta }_{x_j
x_i} \to u^\beta_{x_j x_i}$ in $L^2(\Omega)$, for $\beta =1, \dots
N$
   and $i, j =1, \dots n$.
 The same
observation, combined with  H\"older's inequality, ensures that also
the second integral  converges to $0$ as $k \to \infty$. Since
$\nabla \u_k \to \nabla \u$ a.e. in $\Omega$, the last integral in
\eqref{conv1} tends to
\begin{equation}\label{convlim}
\int _{{\{t < |\nabla \u | <{t+h}\}}} b(|\nabla \u  |)\sum
_{i=1}^n\Big(\sum _{\beta =1}^N \sum_{j}^n
\frac{u^\beta_{x_j}}{|\nabla \u  |}u^\beta _{x_j x_i }\Big)^2\, dx,
\end{equation}
 by to the dominated convergence theorem for
integrals, and the expression  \eqref{convlim} agrees with $\int
_{\{t < |\nabla \u  | <{t+h}\}} b(|\nabla \u  |) \,|\nabla |\nabla
\u ||^2 \, dx$.  Thus, we have shown that
\begin{equation}\label{conv2}
\int _{\{t < |\nabla \u _k | <{t+h}\}} b(|\nabla \u _k |) |\nabla
|\nabla \u_k ||^2 \, dx \to \int _{\{t < |\nabla \u  | <{t+h}\}}
b(|\nabla \u  |) \,|\nabla |\nabla \u ||^2 \, dx
\end{equation}
as $k \to \infty$. A similar argument implies, via \eqref{operator},
that
\begin{equation}\label{conv3}
\int _{\{t < |\nabla \u _k | <{t+h}\}} |\nabla \u _k |\, |\nabla
|\nabla \u_k ||\, |{\rm {\bf div}} (a(|\nabla {\bf u}_k|)\nabla
\u_k)|\, dx \to \int _{\{t < |\nabla \u  | <{t+h}\}} |\nabla \u
|\,|\nabla |\nabla \u  || \, |{\bf f}(x)|\, dx
\end{equation}
as $k \to \infty$. Moreover, equation \eqref{operator} and the
boundedness of $\tfrac 1a$ entail that the sequence
$$\int _{\{|\nabla {\bf u}_k |
>\tau\}} \frac 1{a(|\nabla {\bf u}_k|)} |{\rm {\bf div}} (a(|\nabla {\bf u}_k|)\nabla \u_k)|^2
dx$$ is uniformly bounded for $\tau >0$, and that,  for every $\tau
>0$,
\begin{equation}\label{conv4}
\int _{\{|\nabla {\bf u}_k |
>\tau\}} \frac 1{a(|\nabla {\bf u}_k|)} |{\rm {\bf div}} (a(|\nabla {\bf u}_k|)\nabla \u_k)|^2
dx \to \int _{\{|\nabla {\bf u} |
>\tau\}} \frac 1{a(|\nabla {\bf u}|)} |{\bf f}(x)|^2
dx
\end{equation}
as $k \to \infty$.
Consequently,
\begin{equation}\label{conv5}
\int _t^{t+h} \int _{\{|\nabla {\bf u}_k |
>\tau\}} \frac 1{a(|\nabla {\bf u}_k|)} |{\rm {\bf div}} (a(|\nabla {\bf u}_k|)\nabla \u_k)|^2
dx\, d\tau \to \int _t^{t+h} \int _{\{|\nabla {\bf u} |
>\tau\}} \frac 1{a(|\nabla {\bf u}|)} |{\bf f}(x)|^2
dx\, d\tau
\end{equation}
as $k \to \infty$. Let us finally focus on the  last integral on the
right-hand side of \eqref{32k}. For a.e. $\tau >0$,
\begin{equation}\label{frontiere}
\partial \Omega \cap
\partial\{|\nabla{\bf u}_k |>t\} = \{{\rm Tr}\, |\nabla \u_k |> \tau
\}\quad \hbox{up to subsets of $\partial \Omega $ of $\hh$ measure
zero}
\end{equation}
for $k \in \N$, and
\begin{equation}\label{frontiere1}
\partial \Omega \cap
\partial ^M \{|\nabla{\bf u} |>t\} = \{{\rm Tr}\, |\nabla \u |> \tau
\}\quad \hbox{up to subsets of $\partial \Omega $ of $\hh$ measure
zero.}
\end{equation}
Equations \eqref{frontiere}  follow, for instance, from a close
inspection of the proof of \cite[Lemma 6.5.1/2]{Mazlibro}. By
\eqref{convtr}, for a.e. $\tau
>0$,
\begin{equation}\label{frontiere2}
\chi_{\{{\rm Tr}\, |\nabla \u_k |> \tau \}}(x) |{\rm tr} \mathcal B
(x)| \to \chi_{\{{\rm Tr}\, |\nabla \u |> \tau \}}|{\rm tr} \mathcal
B (x)| \quad \hbox{$\hh$ a.e. on $\partial \Omega$.}
\end{equation}
Hence, by the dominated convergence theorem for integrals,
\begin{equation}\label{frontiere3}
\int _{\partial \Omega} \chi_{\{{\rm Tr}\, |\nabla \u_k |> \tau
\}}(x) |{\rm tr} \mathcal B (x)| \, d\hh (x) \to \int _{\partial
\Omega} \chi_{\{{\rm Tr}\, |\nabla \u |> \tau \}}(x) |{\rm tr}
\mathcal B (x)| \, d\hh (x),
\end{equation}
and the first integral in \eqref{frontiere3} is uniformly bounded
for $\tau
>0$. Thus,
\begin{multline}\label{frontiere4} \int _t^{t+h}
a(\tau)\tau ^2 \int_{\partial \Omega } \chi_{\{{\rm Tr}\, |\nabla
\u_k |> \tau \}}(x)\, |{\rm tr} \mathcal B (x)| d\hh (x)\,d\tau \\ \to
\int _t^{t+h} a(\tau)\tau ^2 \int_{\partial \Omega} \chi_{\{{\rm
Tr}\, |\nabla \u |> \tau \}}(x) \, |{\rm tr} \mathcal B (x)| d\hh
(x)\,d\tau
\end{multline}
as $k \to \infty$, whence, by \eqref{frontiere} and
\eqref{frontiere1},
\begin{equation}\label{frontiere5} \int _t^{t+h}
a(\tau)\tau ^2 \int_{\partial \Omega \cap
\partial\{|\nabla{\bf u}_k |>\tau \}}  |{\rm tr} \mathcal B (x)| d\hh (x)\,d\tau \to
\int _t^{t+h} a(\tau)\tau ^2 \int_{\partial \Omega \cap
\partial ^M \{|\nabla{\bf u} |>\tau \}} |{\rm tr} \mathcal B (x)| d\hh (x)\,d\tau
\end{equation}
as $k \to \infty$. Combining \eqref{32h}, \eqref{conv2},
\eqref{conv3}, \eqref{conv5} and \eqref{frontiere5} tells that
\begin{align}\label{32uh}
 C \int _{\{t < |\nabla \u  | <{t+h}\}} b(|\nabla \u  |) |\nabla
|\nabla \u ||^2 \, dx &  \leq \int _{\{t < |\nabla \u  | <{t+h}\}}
|\nabla \u ||\nabla |\nabla \u  || |{\bf f}(x)|\, dx
 \\ \nonumber & \quad + \int _t^{t+h} \int _{\{|\nabla {\bf u} |
>\tau\}} \frac 1{a(|\nabla {\bf u}|)} |{\bf f}(x)|^2
dx\, d\tau
 \\ \nonumber & \quad +
\int _t^{t+h} a(\tau)\tau ^2 \int_{\partial \Omega \cap
\partial ^M \{|\nabla{\bf u} |>\tau \}} |{\rm tr} \mathcal B (x)| d\hh (x)\,d\tau
.
\end{align}
Dividing through by $h$ in \eqref{32uh}, making use of the coarea
formula again, and passing to the limit as $h \to 0^+$ yields
\eqref{32u}.

\smallskip
\par\noindent
{\bf Step 3}. Here, we
show that, given $r>n-1$,
\begin{equation}\label{estimate}
\|\nabla \u\|_{L^\infty (\Omega )} \leq C b^{-1}\big(\|{\bf
f}\|_{L^{n,1}(\Omega , \rN)}\big)
\end{equation}
for some constant $C$ depending on $i_a$, $s_a$, $n$, $N$, $r$,
$\|{\rm tr} \mathcal B \|_{L^r(\partial \Omega )}$, $|\Omega |$, and
on the constant in  \eqref{isopboundary}. More precisely, hereafter
dependence on $i_a$ and $|\Omega|$ will mean just through a lower
bound, and dependence on $s_a$, $\|{\rm tr} \mathcal B
\|_{L^r(\partial \Omega )}$, and on the constant in
\eqref{isopboundary} is just through an upper bound.
%
\par\noindent
By the Hardy-Littlewood inequality \eqref{HL},
\begin{equation}\label{estk}
\int_{\partial \Omega \cap
\partial ^M \{|\nabla{\bf u} |>t \}}|{\rm tr} \mathcal B (x)| d\hh (x) \leq
\int _0^{\hh (\partial \Omega \cap \partial ^M \{|\nabla{\bf u} |>t
\})} ({\rm tr} \mathcal B)^*(r)dr \qquad \hbox{for a.e. $t>0$},
\end{equation}
where $({\rm tr} \mathcal B)^*$ denotes the decreasing rearrangement
of ${\rm tr} \mathcal B$ with respect to the measure $\hh$ on
$\partial \Omega$. Since $|\nabla \u|$ is a suitably represented
Sobolev function, for a.e. $t>0$,
$$ \Omega \cap \partial ^M \{|\nabla \u |
>t\} = \{|\nabla \u |
=t\}, \quad \hbox{up to sets of $\hh$ measure zero}$$ (see e.g.
\cite{BZ}.) Thus,
\begin{align}\label{36}
\hh (\partial \Omega \cap \partial ^M \{|\nabla \u |
>t\})\leq C \hh ( \{|\nabla \u |
=t\}) \qquad \hbox{for a.e. $t \geq |\nabla \u|^*(|\Omega | /2)$,}
\end{align}
where $C$ is the constant in \eqref{isopboundary}. Denote  the
distribution function $\mu _{|\nabla \u|}$ of $|\nabla \u|$, defined
as in \eqref{mu}, simply by  $\mu$.
By \eqref{relisop},
\begin{align}\label{36bis}
\mu (t)^{1/n'}\leq C \hh ( \{|\nabla \u | =t\}) \qquad \qquad
 \hbox{for a.e. $t \geq |\nabla \u|^*(|\Omega | /2)$,}
\end{align}
where $C$ is a constant depending on $n$ and on the constant in
\eqref{isopboundary}. From \eqref{36} and \eqref{36bis}, we obtain
that
\begin{align}\label{new1}
\int _0^{\hh (\partial \Omega \cap \partial ^M \{|\nabla \u |>t\})}
({\rm tr} \mathcal B)^*(r)dr & \leq \int _0^{C \hh (\{|\nabla \u |=t\})} ({\rm tr} \mathcal B)^*(r)dr \\
\nonumber &= C
\hh (\{|\nabla \u |=t\})\, ({\rm tr} \mathcal B)^{**}(C\hh (\{|\nabla \u |=t\})) \\
\nonumber & \leq C \hh (\{|\nabla \u |=t\}) \,({\rm tr} \mathcal
B)^{**}\big(C'\mu (t)^{1/n'}\big)
\end{align}
 for a.e. $t \geq |\nabla \u|^*(|\Omega |/2)$.
 Observe
 that the last inequality holds since $({\rm tr} \mathcal B)^{**}$ is a non-increasing
 function.
Coupling \eqref{32u} with \eqref{new1} tells us that
\begin{align}\label{32final}
 C b(t)\int _{\{|\nabla \u  | =t\}}  |\nabla
|\nabla \u  || \, d\hh (x)
&  \leq t  \int _{\{|\nabla {\bf u} | =t\}}|{\bf f}(x)| d\hh
 (x)
  +   \int _{\{|\nabla {\bf u} |
>t\}} \frac 1{a(|\nabla {\bf u}|)} |{\bf f}(x)|^2
dx \\ \nonumber & \quad +    a(t)t^2  \hh (\{|\nabla \u |=t\})\,
({\rm tr} \mathcal B)^{**}\big(C'\mu (t)^{1/n'}\big)  \,\, \hbox{for
a.e. $t> t_{\u}$,}
\end{align}
for some constants $C=C(\Omega, \min\{i_a , 0\})$ and $C'=C'(\Omega
)$.
\par\noindent Now, we distinguish into the cases when
  $a$ is non-decreasing or
non-increasing.
 \par\noindent First, assume that $a$ is
non-decreasing. Then we infer from \eqref{32final} that
\begin{align}\label{32incr}
 C b(t)\int _{\{|\nabla \u | =t\}}  |\nabla
|\nabla \u || \, d\hh (x) &  \leq t  \int _{\{|\nabla {\bf u} |
=t\}}|{\bf f}(x)| d\hh
 (x)
 \\ \nonumber & \quad + \frac 1{a(t)} \int _{\{|\nabla {\bf u} |
>t\}} |{\bf f}(x)|^2
dx
\\ \nonumber & \quad +  a(t)t^2  \hh (\{|\nabla \u |=t\})
({\rm tr} \mathcal B)^{**}\big(C'\mu (t)^{1/n'}\big) \,\, \hbox{for
a.e. $t> t_\u$}.
\end{align}
By H\"older's inequality,   \eqref{coareadiff} and
 \eqref{38},
\begin{align}\label{39}
 \int _{\{|\nabla \u | =t\}}|{\bf f}(x)|\, d\hh (x)  & \leq \bigg(\int _{\{|\nabla \u |
=t\}}\frac{|{\bf f}(x)|^2}{|\nabla |\nabla \u ||} d\hh
(x)\bigg)^{1/2} \bigg(\int _{\{|\nabla \u | =t\}}|\nabla |\nabla \u
|| d\hh (x)\bigg)^{1/2} \\ \nonumber & \leq \bigg(- \frac{d}{dt}
\int _{\{|\nabla \u |
>t\}}|{\bf f}(x)|^2 dx\bigg)^{1/2}
\bigg(- \frac{d}{dt} \int _{\{|\nabla \u |
>t\}}|\nabla
|\nabla \u ||^2 dx\bigg)^{1/2}
\end{align}
for a.e. $t>0$. An analogous chain as in \eqref{39}, with $|{\bf
f}(x)|$ replaced with $1$ yields
\begin{align}\label{40}
\hh (\{|\nabla \u | =t\})
 \leq (- \mu ' (t))^{1/2}
\bigg(- \frac{d}{dt} \int _{\{|\nabla \u |
>t\}}|\nabla
|\nabla \u ||^2 dx\bigg)^{1/2} \qquad \hbox{for a.e. $t>0$}.
\end{align}
By the Hardy-Littlewwod inequality \eqref{HL},
\begin{equation}\label{41}
\int _{\{|\nabla \u |
>t\}}|{\bf f}(x)|^2  dx \leq \int _0^{\mu (t)} |{\bf f}|^*(r)^2 dr \qquad
\hbox{for $t>0$}.
\end{equation}
Inequalities \eqref{32incr} --
%
\eqref{41}, and inequality  \eqref{tal2}
applied with $v=|\nabla \u|$ entail that
\begin{align}\label{42}
 C b(t) & \bigg(- \frac{d}{dt} \int _{\{|\nabla \u |
>t\}}|\nabla
|\nabla \u ||^2 dx\bigg)
 \\ \nonumber & \leq t  \bigg(- \frac{d}{dt} \int
_{\{|\nabla \u |
>t\}}|{\bf f}(x)|^2 dx\bigg)^{1/2}
\bigg(- \frac{d}{dt} \int _{\{|\nabla \u |
>t\}}|\nabla
|\nabla \u ||^2 dx\bigg)^{1/2}
\\ \nonumber &
 \quad +
  \frac {1}{ a(t)}
 (-\mu '(t))^{1/2}\mu (t)^{-1/n'} \int _0^{\mu (t)} |{\bf f}|^*(r)^2 dr
 \bigg(- \frac{d}{dt} \int
_{\{|\nabla \u |>t \}}|\nabla |\nabla \u ||^2 dx\bigg)^{1/2}
\\ \nonumber &
 \quad +   a(t)t^2   (- \mu ' (t))^{1/2} ({\rm tr} \mathcal B)^{**}\big(C'\mu (t)^{1/n'}\big)
\bigg(- \frac{d}{dt} \int _{\{|\nabla \u |
>t\}}|\nabla
|\nabla \u ||^2 dx\bigg)^{1/2}
\end{align}
 for a.e. $t > t_\u$, for some constants $C=C(\Omega, \min\{i_a , 0\})$ and $C'=C'(\Omega
)$. By \eqref{tal2}
 with $v=|\nabla \u|$, we have that  $- \frac{d}{dt} \int _{\{|\nabla \u |
>t\}}|\nabla
|\nabla \u ||^2 dx >0$ for a.e. $t > t_\u$. Hence, we  may divide
 through  by $- \frac{d}{dt} \int _{\{|\nabla \u |
>t\}}|\nabla
|\nabla \u ||^2 dx$ in \eqref{42}, and exploit  \eqref{tal2}
 with $v=|\nabla \u|$ again to obtain
\begin{align}\label{43}
 C b(t) &
  \leq   t  (-\mu '(t))^{1/2}\mu (t)^{-1/n'}\bigg(- \frac{d}{dt} \int
_{\{|\nabla \u |
>t\}}|{\bf f}(x)|^2 dx\bigg)^{1/2}
\\ \nonumber &
 +
 \frac {1}{  a(t)}
(-\mu '(t))\mu (t)^{-2/n'} \int _0^{\mu (t)} |{\bf f}|^*(r)^2 dr
\\ \nonumber &
 +
   a(t)t^2
(-\mu '(t))\mu (t)^{-1/n'} ({\rm tr} \mathcal B)^{**}\big(C'\mu
(t)^{1/n'}\big)\quad \hbox{for a.e. $t > t_\u$,}
\end{align}
 for some constants $C=C(\Omega, \min\{i_a , 0\})$ and $C'=C'(\Omega
)$.
Since $|\nabla \u|$ is a Sobolev function, the function $|\nabla
\u|^*$ is (locally absolutely) continuous \cite[Lemma 6.6]{CEG}, and
$|\nabla \u|^*(\mu (t) )=t$ for $t>0$. Define the function $\phi
_{\bf f} : (0, \mo ) \to [0, \infty )$ as
\begin{equation}\label{46}
\phi _{\bf f}(s) = \bigg( \frac{d}{ds} \int _{\{|\nabla \u |
>|\nabla \u |^*(s)\}}|{\bf f}(x)|^2 dx\bigg)^{1/2} \qquad \quad \hbox{for a.e. $s
\in (0, \mo )$.}
\end{equation}
As a consequence,
\begin{equation}\label{44}
\bigg(- \frac{d}{dt} \int _{\{|\nabla \u |>t \}}|{\bf f}(x)|^2
dx\bigg)^{1/2} = (-\mu '(t))^{1/2} \phi _{\bf f}(\mu (t)) \qquad
\quad \hbox{for a.e. $t>0$,}
\end{equation}
and, by \cite[Proposition 3.4]{CMlipschitz},
\begin{equation}\label{45}
\int _0^s \phi _{\bf f} ^*(r)^2 dr \leq \int _0^s |{\bf f}|
^*(r)^2 dr \qquad \quad \hbox{for $s \in (0, \mo )$.}
\end{equation}
We thus deduce from inequality \eqref{43} that
%
\begin{multline}\label{47}
 C a(t) b(t)
  \leq    b(t) (-\mu '(t))\mu (t)^{-1/n'} \phi _{\bf f} (\mu
(t))
\\  +
(-\mu '(t))\mu (t)^{-2/n'} \int _0^{\mu (t)} |{\bf f}|^*(r)^2 dr +
  b(t)^2
(-\mu '(t))\mu (t)^{-1/n'}({\rm tr} \mathcal B)^{**}\big(C'\mu
(t)^{1/n'}\big)
\end{multline}
for a.e. $t > t_\u$. Let $t_\u \leq  t_0 < T < \|\nabla
\u\|_{L^\infty (\Omega , \R^{Nn})}$, and let $H$ be the function
defined by \eqref{Hdef}. On estimating $b(t)$ by $b(T)$ for $t \in
(t_0, T)$ on the right-hand side of \eqref{47}, and integrating the
resulting inequality over $(t_0, T)$ yields
\begin{align}\label{47'}
 C H(T)
& \leq C H(t_0) +   b(T) \int _{t_0}^T (-\mu '(t))\mu (t)^{-1/n'}
\phi _{\bf f} (\mu (t)) dt
\\ \nonumber & \quad +
   \int _{t_0}^T
(-\mu '(t))\mu (t)^{-2/n'} \int _0^{\mu (t)} |{\bf f}|^*(r)^2 dr\, dt
+
 b(T)^2
\int _{t_0}^T (-\mu '(t))\mu (t)^{-1/n'}({\rm tr} \mathcal
B)^{**}\big(C'\mu (t)^{1/n'}\big)dt
\\ \nonumber
& \leq C H(t_0) +   b(T) \int _{\mu (T)}^{\mu (t_0)} s^{-1/n'} \phi
_{\bf f}
  (s)\,ds
\\ \nonumber & \quad +
   \int _{\mu (T)}^{\mu (t_0)}
s^{-2/n'} \int _0^{s} |{\bf f}|^*(r)^2 dr \, ds +
 b(T)^2
\int _{\mu (T)^{1/n'}}^{\mu (t_0)^{1/n'}} ({\rm tr} \mathcal B)^{**}(C's) s^{\frac
1{n-1}} \frac{ds}{s}.
\end{align}
Hence, owing to  \eqref{H},
\begin{align}\label{47new}
b(T)^2 & \leq C b(t_0)^2 + C  b(T) \int _{0}^{\mu (t_0)} s^{-1/n'}
\phi _{\bf f}
  (s)\,ds
\\ \nonumber & \quad +
 C  \int _{0}^{\mu (t_0)}
s^{-2/n'} \int _0^{s} |{\bf f}|^*(r)^2 dr \, ds +
 C b(T)^2
\int _{0}^{\mu (t_0)^{1/n'}} ({\rm tr} \mathcal B)^{**}(C's)
s^{\frac 1{n-1}} \frac{ds}{s}
\end{align}
for some constants $C=C(\Omega , i_a , s_a)$ and $C'=C'(\Omega )$.
Note that the last integral is actually finite, since ${\rm tr}
\mathcal B \in L^r(\partial \Omega )$, and  $L^r(\partial \Omega
)\to L^{n-1, 1}(\partial \Omega )$ for $r
>n-1$, by \eqref{lorentz12}.
%
Define the function $G: [0, \infty ) \to [0, \infty )$ as
\begin{equation}\label{G}
G(s) = C \int _{0}^{s^{1/n'}} ({\rm tr} \mathcal B)^{**}(C'r)
r^{\frac 1{n-1}} \frac{dr}{r} \qquad \hbox{for $s \geq 0$,}
\end{equation}
 where $C$ and $C'$ are as in \eqref{47new}.
Set $s_0 = \min \{\alpha |\Omega | ,
G^{-1}(\frac 1{2C} )\}$,
 where $\alpha$ is given by \eqref{alpha}, with $\v$ replaced with
 $\u$,
and choose
$$t_0 = |\nabla \u|^* (s_0).$$
 One has that $t_0 \geq t_\u$, inasmuch as  $s_0 \leq \alpha \mo$.
 Moreover, since $\mu (t_0) \leq G^{-1}(\frac 1{2C} )$,
$$C  \int _{0}^{\mu (t_0)^{1/n'}}
({\rm tr} \mathcal B)^{**}(C'r) r^{\frac 1{n-1}} \frac{dr}{r} \leq
\tfrac 12.$$ From \eqref{47new}  we thus infer that
\begin{align}\label{47''}
 b(T)^2
  \leq C b(t_0)^2 + C  b(T) \int _{0}^{|\Omega |} s^{-1/n'} \phi _{\bf f}
  (s)\,ds
 +
 C  \int _{0}^{|\Omega |}
s^{-2/n'} \int _0^{s} |{\bf f}|^*(r)^2 dr \, ds
\end{align}
for some constant $C=C(\Omega , i_a, s_a)$. By \eqref{45} and
\cite[Lemma 3.5]{CMlipschitz}, there exists a constant $C=C(n)$
such that
\begin{equation}\label{50}
\int _{0}^{|\Omega |}
    s^{-1/n'} \phi_{\bf f} (s)ds \leq C  \|{\bf f}\|_{L^{n,1}(\Omega ,
    \rN)}.
    \end{equation}
   Moreover,  by \cite[Lemma 3.6]{CMlipschitz}, there
exists a constant $C=C(n)$
 such that
\begin{equation}\label{51}
\int _{0}^{|\Omega |}s^{-2/n'} \int _0^{s} |{\bf f}|^*(r)^2 dr\,ds
\leq C \|{\bf f}\|_{L^{n,1}(\Omega , \rN)}^2.
\end{equation}
Owing to \eqref{47''}--\eqref{51}, there exists a constant
$C=C(\Omega ,  i_a, s_a)$ such that
\begin{align}\label{52}
b(T)^2  \leq C b(t_0)^2 + C
 b(T)
\|{\bf f}\|_{L^{n,1}(\Omega , \rN)}  +
 C \|{\bf f}\|_{L^{n,1}(\Omega , \rN)}^2.
%
\end{align}
Thus,
\begin{align}\label{53}
b(T)  \leq C b(t_0) + C \|{\bf f}\|_{L^{n,1}(\Omega , \rN )}
 \end{align}
 for some constant $C=C(\Omega ,  i_a, s_a)$.
 Next, let $\beta , \psi: [0, \infty ) \to [0, \infty )$   be the functions defined by $\beta (t) =
 b(t)t$ for $t \geq 0$ and $\psi (s ) = s b^{-1}(s )$ for $s \geq 0$. Proposition \ref{energia} and inequality \eqref{B5} ensure that
 \begin{align}\label{54}
C \psi\big(\|{\bf f}\|_{L^{n,1}(\Omega , \rN)}\big) \geq \int
_\Omega \beta (|\nabla \u|) dx \geq \int _{\{|\nabla \u| \geq t_0\}}
\beta (|\nabla \u|) dx \geq \beta (t_0) \lim _{t \to t_0^-}\mu (t)
\geq \beta (t_0) s_0,
%
%
\end{align}
for some constant $C=C(\Omega,  i_a , s_a)$, whence, by \eqref{B7},
\begin{align}\label{55}
 \beta (t_0)  \leq  \psi \big(C \|{\bf f}\|_{L^{n,1}(\Omega , \rN)}\big),
\end{align}
for some constant $C=C(\Omega,  r, i_a , s_a)$. Since $b(\beta
^{-1}(\psi (s )))= s $ for $s \geq 0$, inequality \eqref{55} implies
that
\begin{align}\label{56}
 b(t_0)  \leq C
 \|{\bf f}\|_{L^{n,1}(\Omega , \rN )}.
\end{align}
Hence, by \eqref{53},
\begin{equation}\label{56bis}
b(T)  \leq  C
 \|{\bf f}\|_{L^{n,1}(\Omega , \rN)}
 \end{equation}
 for some constant $C=C(\Omega, r,  i_a , s_a)$. Taking the limit as
 $T \to \|\nabla \u\|_{L^\infty (\Omega , \mathbb R^{Nn})}$ in \eqref{56bis}, and making use
 of
 \eqref{B7bis}, yields inequality \eqref{estimate}. An inspection of
 the proof shows that the constant in \eqref{estimate} actually
 depends on the specified quantities.
\par
Assume next that $a$ is non-increasing. From \eqref{32final} we
deduce that
\begin{align}\label{32decr}
 C b(t)\int _{\{|\nabla \u | =t\}}  |\nabla
|\nabla \u || \, d\hh (x)
&  \leq t  \int _{\{|\nabla {\bf u} | =t\}}|{\bf f}(x)| d\hh
 (x)
  + \frac 1{ a(\|\nabla \u \|_{L^\infty (\Omega)})}  \int _{\{|\nabla {\bf u} |
>t\}} |{\bf f}(x)|^2
dx
\\ \nonumber & \quad +   a(t)t^2 \, \hh (\{|\nabla \u |=t\})\,
({\rm tr}\mathcal B)^{**}\big(C'\mu (t)^{1/n'}\big) \,\, \hbox{for
a.e. $t>t_\u$}.
\end{align}
Observe that, although $\|\nabla \u \|_{L^\infty (\Omega , \mathbb
R^{Nn})}$ is not yet known to be finite at this stage, the quantity
$\frac 1{ a(\|\nabla \u \|_{L^\infty (\Omega , \mathbb R^{Nn})})} $
is finite, since assumption \eqref{nondeg} is still in force. On
starting from \eqref{32decr}, instead of \eqref{32incr}, and arguing
as in the proof of \eqref{47}, one can now show that
\begin{align}\label{sharp}
C a(\|\nabla \u \|_{L^\infty (\Omega , \mathbb R^{Nn})}) b(t) & \leq
t
 a(\|\nabla \u \|_{L^\infty (\Omega , \mathbb R^{Nn})})
 (-\mu '(t))\mu (t)^{-1/n'} \phi _{\bf f} (\mu
(t))
\\ \nonumber & \quad +
(-\mu '(t))\mu (t)^{-2/n'} \int _0^{\mu (t)} |{\bf f}|^*(r)^2 dr
\\ \nonumber & \quad +
  a(t) t^2 a(\|\nabla \u \|_{L^\infty (\Omega , \mathbb R^{Nn})})
(-\mu '(t))\mu (t)^{-1/n'}({\rm tr}\mathcal B)^{**}\big(C'\mu
(t)^{1/n'}\big)
\end{align}
for a.e. $t > t_\u$, for some constants $C=C(\Omega, \min\{i_a ,
0\})$ and $C'=C'(\Omega )$. Let us fix $t_0$ and $T$ as above.  For
every $t \in (t_0 , T)$, the expression $ t a(\|\nabla \u
\|_{L^\infty (\Omega , \mathbb R^{Nn})})$ on the right-hand side of
\eqref{sharp} can be estimated from above by $T a(\|\nabla \u
\|_{L^\infty (\Omega , \mathbb R^{Nn})})$. Also,  owing to
\eqref{B5}, the quantity $a(t) t^2$ can be bounded by $C B(T)$ for
some constant $C=C(s_a)$. Integrating the resulting inequality over
$(t_0 , T)$ tells us that
\begin{multline}\label{47decr}
B(T) a(\|\nabla \u \|_{L^\infty (\Omega , \R^{Nn})})  \leq C
 B(t_0) a(\|\nabla \u \|_{L^\infty (\Omega , \mathbb R^{Nn})}) + C T
a(\|\nabla \u \|_{L^\infty (\Omega , \mathbb R^{Nn})}) \int
_{0}^{\mu (t_0)} s^{-1/n'} \phi _{\bf f}
  (s)\,ds
\\  +
 C  \int _{0}^{\mu (t_0)}
s^{-2/n'} \int _0^{s} |{\bf f}|^*(r)^2 dr \, ds   +
 Ca(\|\nabla \u \|_{L^\infty (\Omega , \R^{Nn})}) B(T)
\int _{0}^{\mu (t_0)^{1/n'}} ({\rm tr}\mathcal B)^{**}(C's) s^{\frac
1{n-1}} \frac{ds}{s}
\end{multline}
for some constants $C=C(\Omega , i_a , s_a)$ and $C'=C'(\Omega )$.
Exploiting inequality \eqref{47decr} instead of \eqref{47new}, and
arguing as in the proof of \eqref{52} lead  to
\begin{multline}\label{52decr}
B(T) a(\|\nabla \u \|_{L^\infty (\Omega , \R^{Nn})}) \\ \leq C
B(t_0)a(\|\nabla \u \|_{L^\infty (\Omega , \R^{Nn} )}) + C T
 a(\|\nabla \u \|_{L^\infty (\Omega , \R^{Nn})})
 \|{\bf f}\|_{L^{n,1}(\Omega , \rN)}  +
 C \|{\bf f}\|_{L^{n,1}(\Omega , \rN)}^2
\end{multline}
for some constant $C=C(\o , i_a , s_a)$. Dividing through by $T$,
and recalling \eqref{B5} entail that
\begin{multline}\label{53decr}
b(T) a(\|\nabla \u \|_{L^\infty (\Omega , \R^{Nn})}) \\ \leq C
\frac{B(t_0)}{T} a(\|\nabla \u \|_{L^\infty (\Omega , \R^{Nn})}) + C
 a(\|\nabla \u \|_{L^\infty (\Omega , \R^{Nn})})
 \|{\bf f}\|_{L^{n,1}(\Omega , \rN)}  +
 \frac{C}T \|{\bf f}\|_{L^{n,1}(\Omega , \rN)}^2
\end{multline}
for some constant $C=C(\o ,  i_a , s_a)$. The limit as $T$ goes to
$\|\nabla \u \|_{L^\infty (\Omega , \mathbb R^{Nn})}$ of the
right-hand side of \eqref{53decr} is obviously finite, and hence the
limit of the left-hand side is finite as well. Thus, $ \|\nabla \u
\|_{L^\infty (\Omega , \mathbb R^{Nn})} < \infty$, since $\lim _{T
\to \infty } b(T) = \infty$. Taking $T= \|\nabla \u \|_{L^\infty
(\Omega , \mathbb R^{Nn})}$ in \eqref{53decr}, and multiplying
through by $\|\nabla \u \|_{L^\infty (\Omega , \mathbb R^{Nn})}$
yields
\begin{multline}\label{54decr}
b(\|\nabla \u \|_{L^\infty (\Omega , \mathbb R^{Nn})})^2 \\ \leq C
B(t_0) a(\|\nabla \u \|_{L^\infty (\Omega , \mathbb R^{Nn})}) + C
 b(\|\nabla \u \|_{L^\infty (\Omega , \mathbb R^{Nn})})
  \|{\bf f}\|_{L^{n,1}(\Omega , \rN)}  +
 C \|{\bf f}\|_{L^{n,1}(\Omega , \rN)}^2.
\end{multline}
Observe that, by \eqref{B5} and \eqref{56},
\begin{multline}\label{55decr}
B(t_0) a(\|\nabla \u \|_{L^\infty (\Omega , \mathbb R^{Nn})})  \leq
C t_0 b(t_0)a(\|\nabla \u \|_{L^\infty (\Omega , \mathbb R^{Nn})})
\\ \leq  C b(t_0)b(\|\nabla \u \|_{L^\infty (\Omega , \mathbb R^{Nn})}) \leq  C'
b(\|\nabla \u \|_{L^\infty (\Omega , \mathbb R^{Nn})}) \|{\bf
f}\|_{L^{n,1}(\Omega , \rN)},
\end{multline}
for some constants $C=C(s_a)$ and $C' = C'(\o , r, i_a , s_a)$.
 Coupling
\eqref{54decr} with \eqref{55decr} tells us that
$$b(\|\nabla \u \|_{L^\infty (\Omega , \mathbb R^{Nn})})^2 \leq  C
 b(\|\nabla \u \|_{L^\infty (\Omega , \mathbb R^{Nn})})
\|{\bf f}\|_{L^{n,1}(\Omega , \rN)}  +
 C \|{\bf f}\|_{L^{n,1}(\Omega , \rN)}^2
%
$$
 for some constant $C=C(\o , r, i_a , s_a)$. Hence,
\begin{equation}\label{estimater}
b(\|\nabla \u \|_{L^\infty (\Omega , \mathbb R^{Nn})}) \leq C \|{\bf
f}\|_{L^{n,1}(\Omega , \rN)}
%
\end{equation}
 for some constant $C= C(\o , r, i_a , s_a)$, and  \eqref{estimate} follows also in this case,
  with a constant $C$ depending on the specified quantities.

\par\noindent {\bf Step 4}. The present step exploits  some variants of the arguments of Steps 1-3 in order to
 show that inequality \eqref{estimate} holds with a constant $C$
 which only depends on $i_a$, $s_a$, $n$, $|\Omega |$, the constant in \eqref{isopboundary}, and on  ${\rm tr} \mathcal B$ just
 through (an upper bound for) the norm
 $\|{\rm tr}
B\|_{L^{n-1,1}(\partial \Omega )}$, instead of a
stronger norm $\|{\rm tr} B\|_{L^r(\partial \Omega )}$ with $r>n-1$.
%
%
%
%
The piece of information that was missing until this stage, and
makes this further step possible, is that  the solution  $\u$ is now
already known to satisfy
\begin{equation}\label{gradientinf}
\|\nabla \u \|_{L^\infty (\Omega , \mathbb R^{Nn})} < \infty,
\end{equation}
and hence we can exploit inequality \eqref{37old} in the place of
\eqref{37}. By  \eqref{gradientinf}, ${\bf u} \in
W^{1,\infty}_0(\Omega , \rN) \cap W^{2,2}(\Omega , \rN)$.
%
%
%
Hence,  there exists a sequence $\{{\bf u}_k\} \subset C^\infty
(\Omega , \rN)\cap C^2(\overline \Omega , \rN)$ fulfilling \eqref{convuk}--\eqref{operator},  such that
 $\u_k = 0$
on $\partial \Omega$, and, in addition,
\begin{equation}\label{convukinf}
 \|\nabla \u _k \|_{L^\infty (\Omega , \mathbb
R^{Nn})} \to \|\nabla \u \|_{L^\infty (\Omega , \mathbb R^{Nn})} \quad \hbox{as $k \to \infty$.}
\end{equation}
\par\noindent
Inequality \eqref{37old} holds with ${\bf v}$ replaced with ${\bf
u}_k$. An analogous argument as in Step 2 shows that the same
inequality continues to hold for $\u$, that is
\begin{align}\label{37oldu}
C  b(t) \int _{\{|\nabla \u | =t\}}  |\nabla |\nabla \u ||& \, d\hh
(x)
  \leq t  \int _{\{|\nabla \u | =t\}} |{\bf f} (x)|d\hh (x)  +
 \frac {\|\nabla {\bf u} \|_{L^\infty (\Omega , \mathbb R^{Nn})}}{ b(t)}  \int _{\{|\nabla \u |>t\}}
|{\bf f} (x)|^2  dx \\ \nonumber & \quad  + a\big(\|\nabla {\bf u}
\|_{L^\infty (\Omega , \mathbb R^{Nn})}\big)\|\nabla {\bf u}
\|_{L^\infty (\Omega , \mathbb R^{Nn})}^2 \int_{\partial \Omega \cap
\partial\{|\nabla{\bf u} |>t\}}|{\rm tr} \mathcal B (x)| d\hh (x)
\end{align}
for a.e. $t > 0$, where $\tfrac {(1 + \min \{i_a , 0\})^2}2$. We now
start from \eqref{37oldu}, make use of arguments similar to -- and
even simpler than -- those which lead to either \eqref{47} or
\eqref{sharp} from \eqref{32u} (in particular, now we do not need to
distinguish into the cases when $a$ is non-decreasing or
non-increasing), and show that
\begin{align}\label{47inf}
 C b(t)^2 &
  \leq   b\big(\|\nabla \u\|_{L^\infty (\Omega , \mathbb R^{Nn})}\big)\|\nabla
\u\|_{L^\infty (\Omega , \mathbb R^{Nn})} (-\mu '(t))\mu (t)^{-1/n'}
\phi _{\bf f}(\mu (t))
\\ \nonumber &
 +
  \|\nabla \u\|_{L^\infty (\Omega , \mathbb R^{Nn})}
(-\mu '(t))\mu (t)^{-2/n'} \int _0^{\mu (t)} |{\bf f}|^*(r)^2 dr
   \\ \nonumber & +
  b\big(\|\nabla \u\|_{L^\infty (\Omega , \mathbb
R^{Nn})}\big)^2\|\nabla \u\|_{L^\infty (\Omega ,\mathbb R^{Nn})} (-
\mu ' (t))\mu (t)^{-1/n'}({\rm tr} \mathcal B)^{**}\big(C'\mu (t)^{1/n'}\big)\\
\nonumber & \quad \quad\quad
\quad\quad\quad\quad\quad\quad\quad\qquad \qquad\qquad
\qquad\hbox{for a.e. $t\in [|\nabla \u|^*(|\Omega | /2),\|\nabla
\u\|_{L^\infty (\Omega ,\mathbb R^{Nn})}]$,}
\end{align}
for some positive constants $C=C(\Omega , \min\{i_a
, 0\})$ and $C'=C'(\Omega )$. Moreover, the dependence on $\Omega$
is only through the constant in \eqref{isopboundary}.
\par\noindent
 Let $F$ be the function
defined by \eqref{F}. Given $t_1 \in [|\nabla u|^*(|\Omega |
/2),\|\nabla u\|_{L^\infty (\Omega , \R^{Nn} )}]$, an integration in
\eqref{47inf} yields, via \eqref{B1},
\begin{align}\label{48}
  F(|\nabla
\u|^*(s)) &
 \leq C  F(t_1) +
  C b(\|\nabla \u\|_{L^\infty (\Omega , \R^{Nn}
)})\|\nabla\u\|_{L^\infty (\Omega,  \R^{Nn} )}
 \int _{0}^{\mu (t_1)}
    r^{-1/n'} \phi_{\bf f} (r)dr
\\ \nonumber &
 \quad +
  C \|\nabla\u\|_{L^\infty (\Omega ,  \R^{Nn})}
 \int _{0}^{\mu (t_1)}r^{-2/n'} \int _0^{r} |{\bf f}|^*(\rho)^2
d\rho\,dr
   \\ \nonumber & \quad +  C F\big(\|\nabla\u\|_{L^\infty (\Omega ,  \R^{Nn} )}\big)\int _{0}^{\mu (t_1)^{1/n'}}
({\rm tr} \mathcal B)^{**}(C'r) r^{\frac 1{n-1}} \frac{dr}{r} \quad
\hbox{for $s \in [0, \mu (t_1))$,}
\end{align}
for some  constants  $C=C(\Omega , i_a, s_a)$ and $C'=C'(\Omega )$.
Let $G$ be the function defined as in \eqref{G}, save that  now $C$
and $C'$ are the constants appearing in \eqref{48}.
Set $s_1 = \min \{\frac{|\Omega |}2 , G^{-1}\big(\tfrac
1{2C}\big)\}$, and choose
$$t_1 = |\nabla\u|^* (s_1).$$
Since $\mu (t_1) \leq G^{-1}\big(\tfrac 1{2C}\big)$,
$$ C \int _{0}^{\mu (t_1)^{1/n'}}
({\rm tr} \mathcal B)^{**}(C'r) r^{\frac 1{n-1}} \frac{dr}{r} \leq
\tfrac 12.$$ From this choice of $t_1$, and the choice  $s=0$ in
\eqref{48} we infer that
\begin{align}\label{new2}
F\big(\|\nabla \u\|_{L^\infty (\Omega , \R^{Nn} )}\big) & \leq C F(t_1) +
 C
 b(\|\nabla
\u\|_{L^\infty (\Omega ,  \R^{Nn} )})\|\nabla \u\|_{L^\infty (\Omega
,  \R^{Nn} )}
 \int _{0}^{\mu (t_1)}
    r^{-1/n'} \phi_{\bf f} (r)dr
\\ \nonumber &
 +
 C \|\nabla \u\|_{L^\infty (\Omega , \R^{Nn} )}
 \int _{0}^{\mu (t_1)}r^{-2/n'} \int _0^{r} |{\bf f}|^*(\rho)^2
d\rho\,dr
\end{align}
for some constant $C=C(\Omega , i_a, s_a)$.
%
%
%
From \eqref{new2}, via \eqref{50}, \eqref{51} and \eqref{B1},
%
%
we deduce that there exists a constant $C=C(\Omega , i_a, s_a)$ such
that
\begin{align}\label{52old}
b(\|\nabla \u\|_{L^\infty (\Omega , \R^{Nn} )})^2  \leq C b(t_1)^2 + C
 b(\|\nabla \u\|_{L^\infty (\Omega ,  \R^{Nn}
)})
 \|{\bf f}\|_{L^{n,1}(\Omega , \rN )}  +
 C \|{\bf f}\|_{L^{n,1}(\Omega , \rN )}^2.
\end{align}
An inspection of the proof shows that, in fact, the dependence of
$C$ on $\Omega$ is only through $|\Omega |$ and on the constant in
\eqref{isopboundary}. Starting from \eqref{52old} instead of
\eqref{52}, and arguing as in Step 3, yield \eqref{estimate} with a
constant $C$ depending on $i_a$, $s_a$, $n$, $N$, $|\Omega|$,
$\|{\rm tr} \mathcal B\|_{L^{n-1,1}(\partial \Omega )}$  and on the
constant in \eqref{isopboundary}.

\smallskip
\par\noindent
{\bf Step 5}.
 Here we remove the additional assumption
\eqref{domainsmooth}.
\par\noindent
 Since the space  $C^\infty (U) \cap W^{2}L^{n-1,1}(U)$ is dense in  $W^{2}L^{n-1,1}(U)$ for every open set $U \subset \R ^{n-1}$, there exists a sequence $\{\o _m\}_{m \in \N}$ of domains
$\Omega _m \supset \Omega$ such that $\partial \Omega _m \in
C^\infty$, $|\Omega _m \setminus \Omega| \to 0$, $\Omega _m \to
\Omega$ with respect to the Hausdorff distance, and $\|{\rm tr}
\mathcal B_m\|_{L^{n-1, 1}(\partial \Omega _m)} \leq C$ for some
constant $C=C(\Omega)$, where ${\rm tr} \mathcal B_m$ denotes the
trace of the second fundamental form on $\partial \Omega _m$.
%
%
The sequence $\{\Omega _m\}_{m \in \N}$ can be chosen in such a way
that the constant in \eqref{isopboundary}, and hence the constants
in \eqref{11} and \eqref{tal2}, with $\Omega $ replaced with $\Omega
_m$, are bounded, up to a multiplicative constant independent of
$m$, by the corresponding  constants for $\Omega$. This fact depends, in particular, on
the embedding $W^{2}L^{n-1,1}(U) \to W^{1, \infty}(U)$ for $U
\subset \R ^{n-1}$, which entails  the convergence of the Lipschitz
constants of the functions whose graphs locally agree with $\partial
\Omega _m$ to the Lipschitz constant of the function whose graph coincides  with $\partial
\Omega$ .
%
%
%
%
%
%
%
%
Let ${\bf f}$ be continued by $0$ in $\Omega _m \setminus \Omega$,
and  let $\u_m$ be the solution to \eqref{eqdirichlet} with $\Omega$
replaced with $\Omega_m$. Owing to the estimates for $\u_m$ in
$W^{2,2}_{\rm loc}(\Omega _m)$ \cite[Theorem 8.1]{BF}, for every
open set $\Omega '$ such that $\overline {\Omega '} \subset \Omega$,
there exists a constant $C$ such that
\begin{equation}\label{w22m}
\|\u_m\|_{W^{2,2}(\Omega ')}  \leq C
\end{equation}
for $m \in \N$. Furthermore,   by
estimate \eqref{estimate} (in the form established in Step 4) with
$\Omega$ replaced with $\Omega _m$ and $\u$ replaced with $\u_m$ ,
there exists a a constant $C$ such that
\begin{equation}\label{w1m}
\|\nabla \u_m\|_{L^\infty(\Omega,  \R^{Nn} )}  \leq C
\end{equation}
for $m \in \N$. Note that the constant $C$ in \eqref{w22m} and
\eqref{w1m} is independent of $m$.
Let $s \in [1, \frac {2n}{n-2})$. If $\partial \Omega '$ is smooth,
 then  the embedding $W^{2,2}(\o ', \rN) \to W^{1,s}(\o ',
\rN)$ is compact. Thus, by \eqref{w22m} and \eqref{w1m}, there
exists a function 
%
%
$\u \in
W^{1,\infty}(\Omega , \rN)$, and a subsequence of $\{\u_m\}$, still
denoted by $\{\u_m\}$, such that
$$\u_m \to \u \quad\quad  \hbox{in $W^{1,s}_{\rm loc}(\Omega , \rN )$ }$$
 and
\begin{equation}\label{convgrad}
\nabla \u_m \to \nabla \u \quad \hbox{a.e. in $\Omega $.}
\end{equation}
Since $\u_m =0$ on $\partial \Omega _m$, and $\Omega _m \to \Omega$
in the Hausdorff distance, one can deduce  from \eqref{w1m} that
$\u=0$ on $\partial \Omega$. Thus, in particular, $\u \in
W^{1,B}_0(\Omega , \rN )$. The function $\u$ is the weak solution to
the Dirichlet problem \eqref{eqdirichlet}. Indeed, inasmuch as  $\o
\subset \o _m$ for each $m \in \N$,
\begin{equation}\label{weakdirm}
\int _{\Omega _m} a (|\nabla \u_m|) \nabla \u_m \cdot \nabla \phi
\,dx = \int _{\Omega _m} {\bf f} \cdot \phi \,dx
\end{equation}
for every $\phi \in C_0^\infty (\Omega , \rN)$. Passing to the limit
as $m \to \infty$ in \eqref{weakdirm} yields
\begin{equation}\label{weakdir}
\int _{\Omega } a (|\nabla \u|) \nabla \u \cdot \nabla \phi \,dx =
\int _{\Omega } {\bf f}\cdot  \phi \, dx
\end{equation}
for every $\phi \in C_0^\infty (\Omega , \rN)$, owing to
\eqref{convgrad} and \eqref{w1m}, via the dominated convergence
theorem for integrals. Since $B \in \Delta _2$, the space
$C_0^\infty (\Omega , \rN)$ is dense in $W^{1,B}_0(\Omega , \rN )$.
Thus, \eqref{weakdir} holds  for every  $\phi \in W^{1,B}_0(\Omega ,
\rN )$ as well. Note that, by \eqref{convgrad}, the solution $u$
satisfies
\begin{equation}\label{estimatebis}
\|\nabla \u\|_{L^\infty (\Omega , \R^{Nn})} \leq C b^{-1}\big(\|{\bf
f}\|_{L^{n,1}(\Omega , \rN)}\big),
\end{equation}
since such an estimate is fulfilled, with $\u$ replaced with $\u_m$,
by \eqref{estimate}. Here, the constant $C$ depends on $i_a$, $s_a$,
$n$, $N$, $|\Omega|$, $\|{\rm tr} \mathcal B\|_{L^{n-1,1}(\partial
\Omega )}$ and on the constant in \eqref{isopboundary}.
\par\noindent {\bf Step 6}
We conclude by  removing  assumption
\eqref{nondeg}.
\par\noindent
Let $\{a_\varepsilon \}_{\varepsilon \in (0,1)} $ be the family of
functions  defined in Lemma \ref{approx}, and let $b_\varepsilon$
and $B_\varepsilon$ be as in its statement. Let
 $\u$ be the   weak solution in $W^{1, B}_0 (\Omega , \rN )$
to   problem \eqref{eqdirichlet}, and let
  $\u_{\varepsilon }$ denote the
solution in $W^{1, B_\varepsilon }_0 (\Omega , \rN )$ to the problem
\begin{equation}\label{eqdirichleteps}
\begin{cases}
- {\rm div} (a_\varepsilon(|\nabla \u_\varepsilon  |)\nabla \u_\varepsilon  ) = {\bf f} (x)
 & {\rm in}\,\,\, \o\,, \\
 \u_\varepsilon =0 &
{\rm on}\,\,\,
\partial \o \,.
\end{cases}
\end{equation}
 We claim that \begin{equation}\label{aeepsmis}
\nabla \u_{\varepsilon} \to \nabla \u \qquad \hbox{in measure}
 \end{equation}
as $\varepsilon \to 0^+$, and hence there exists a sequence
$\varepsilon _k \to 0$ such
 that
  \begin{equation}\label{aeeps}
\nabla \u_{\varepsilon _k} \to \nabla \u \qquad \hbox{a.e. in
$\Omega$}
 \end{equation}
 as $k \to \infty$.
By the previous steps, there exists a constant $C=C(i_a , s_a,
\Omega)$ (in particular, independent of $\varepsilon$, owing to
\eqref{indici}), such that
\begin{equation}\label{estimateeps}
b_\varepsilon \big(C \|\nabla \u_\varepsilon\|_{L^\infty (\Omega ,
\R^{Nn})}\big) \leq \|{\bf f}\|_{L^{n,1}(\Omega , \rN)}.
\end{equation}
Thus, by the definition of $b_\varepsilon$ and by \eqref{convb}, it
is easily seen that
\begin{equation}\label{bounde}
\|\nabla \u_\varepsilon\|_{L^\infty (\Omega , \R^{Nn})} \leq C
\end{equation}
for some constant $C$ independent of $\varepsilon$, whence
\begin{equation}\label{boundeps}
\int _\Omega B(|\nabla \u_\varepsilon |)\, dx \leq C
\end{equation}
for some constant $C$ independent of $\varepsilon$.
\par\noindent
We preliminarily observe that, although the function $\u$ need not
belong to $W^{1, B_\varepsilon}_0(\Omega , \rN)$ in the case when
$a$ is non-increasing, it can still be used as a test function in
the weak formulation of problem \eqref{eqdirichleteps}.
Indeed, by \eqref{bounde},  $a_\varepsilon(|\nabla \u_\varepsilon
|)\nabla \u_\varepsilon \in L^\infty (\Omega , \rN)$. Therefore,
since $\u\in W^{1,B}_0(\Omega , \rN)$, and the latter space is
embedded into $W^{1,1}_0(\Omega , \rN)$, the function $\u$ can be
approximated by a sequence $\{\u_k\} \in C^\infty _0(\Omega )$ of
functions such that $\u_k \to \u$ in $L^{n'}(\Omega , \rN)$, and
hence in $L^{n', \infty}(\Omega , \rN)$, and  $\nabla \u_k \to
\nabla \u$ in $L^1(\Omega , \rN)$. This allows one to employ $\u_k$
as a test function in the weak formulation of problem
\eqref{eqdirichleteps}, and then pass to the limit as $k \to
\infty$.
\par\noindent The test function $\phi = \u - \u_{\varepsilon
}$ can thus be used both in the weak formulation of problem
\eqref{eqdirichlet}, and in that of problem \eqref{eqdirichleteps}.
Subtracting the resulting equations yields
\begin{multline}\label{eps1}
\int_\Omega [a_{\varepsilon }(|\nabla \u_{\varepsilon }|)\nabla
\u_{\varepsilon } - a(|\nabla \u_\varepsilon|)\nabla \u_\varepsilon
]\cdot (\nabla \u - \nabla \u_{\varepsilon })\, dx \\ = \int_\Omega
[a(|\nabla \u|)\nabla \u- a(|\nabla \u_\varepsilon|)\nabla
\u_\varepsilon ]\cdot (\nabla \u - \nabla \u_{\varepsilon })\, dx.
\end{multline}
Fix any $\sigma \in (0, 1)$. By the definition of Young conjugate,
\begin{align}\label{eps10}
\bigg|\int_\Omega & [a_{\varepsilon }(|\nabla \u_\varepsilon|)\nabla
\u_\varepsilon- a(|\nabla \u_\varepsilon|)\nabla \u_\varepsilon
]\cdot (\nabla \u - \nabla \u_{\varepsilon })\, dx\bigg|
 \leq \int_\Omega |a_{\varepsilon }(|\nabla \u_\varepsilon|)\nabla \u_\varepsilon -
a(|\nabla \u_\varepsilon|)\nabla \u_\varepsilon |\, |\nabla \u -
\nabla \u_{\varepsilon }|\, dx
\\ \nonumber & \leq
\int_\Omega \widetilde B \big( \tfrac 1\sigma |a_{\varepsilon
}(|\nabla \u _\varepsilon|)\nabla \u_\varepsilon - a(|\nabla
\u_\varepsilon|)\nabla \u_\varepsilon |\big)\,dx + \int_\Omega B
\big(\sigma |\nabla \u - \nabla \u_{\varepsilon }|\big)\, dx.
\end{align}
Since $B$ is a Young function of class $\Delta _2$ and $\sigma \in
(0,1)$, there exists a constant $C=C(B)$ such that $B(\sigma (t+s))
\leq C \sigma \big(B(t) +  B(s)\big)$ for $t, s \geq 0$. Hence,
owing to \eqref{boundeps}, there exist positive constants $C$ and
$C'$, independent of $\varepsilon$, such that
\begin{equation}\label{eps11}
\int_\Omega  B \big(\sigma |\nabla \u - \nabla \u_{\varepsilon
}|\big)\, dx \leq  C\sigma  \bigg(\int_\Omega  B \big( |\nabla
\u_{\varepsilon }|\big)\, dx + \int_\Omega  B \big(|\nabla
\u|\big)\, dx\bigg) \leq C'\sigma\,.
\end{equation}
Next, fix any $\delta >0$. Let $C$ be the constant appearing in
\eqref{bounde}, and let $t>C$. We have that
%
\begin{align}\label{eps12}
\int_\Omega \widetilde B & \big(\sigma |a_{\varepsilon }(|\nabla \u
_\varepsilon|)\nabla \u_\varepsilon - a(|\nabla \u_\varepsilon|)\nabla \u_\varepsilon |\big)\,dx \\
\nonumber &= \int_{\{|\nabla u_\varepsilon|\leq t\}} \widetilde B
\big( \sigma |a_{\varepsilon }(|\nabla \u_\varepsilon|)\nabla
\u_\varepsilon- a(|\nabla \u_\varepsilon|)\nabla \u_\varepsilon
|\big)\,dx
\\ \nonumber & \quad +
\int_{\{|\nabla u_\varepsilon|>t\}} \widetilde B \big( \sigma
|a_{\varepsilon }(|\nabla \u_\varepsilon|)\nabla \u_\varepsilon-
a(|\nabla \u_\varepsilon|)\nabla \u_\varepsilon |\big)\,dx.
\end{align}
By the choice of $t$, the last integral vanishes for every
$\varepsilon \in (0, 1)$.
%
%
%
%
%
On the other hand, by
\eqref{conva},
\begin{align}\label{eps14}
\int_{\{|\nabla \u_\varepsilon|\leq t\}} \widetilde B \big( \sigma
|a_{\varepsilon }(|\nabla \u_\varepsilon|)\nabla \u_\varepsilon-
a(|\nabla \u_\varepsilon|)\nabla \u_\varepsilon |\big)\,dx <
\delta\,,
\end{align}
if $\varepsilon $ is sufficently small. Thanks to the arbirtrariness
of $\sigma $ and $\delta$, we infer from
\eqref{eps10}--\eqref{eps14} that
$$\lim _{\varepsilon \to 0^+}
\int_\Omega   [a_{\varepsilon }(|\nabla \u_\varepsilon|)\nabla
\u_\varepsilon- a(|\nabla \u_\varepsilon|)\nabla \u_\varepsilon
]\cdot (\nabla \u - \nabla \u_{\varepsilon })\, dx = 0,$$ and hence,
by \eqref{eps1},
\begin{align}\label{eps15}
\lim _{\varepsilon \to 0^+} \int_\Omega [a(|\nabla \u|)\nabla \u  -
a (|\nabla \u_\varepsilon|)\nabla \u_\varepsilon ]\cdot (\nabla \u
-\nabla \u_\varepsilon )\, dx
 =0.
\end{align}
%
%
We now follow an argument from \cite{BBGGPV}. Fix any $\delta >0$.
Given $t, \tau
>0$, we have that
\begin{multline}\label{417}
|\{|\nabla \u
-\nabla \u_\varepsilon | >t \}| \\
\le |\{|\nabla \u_{\varepsilon }|>\tau \}|+|\{| \nabla \u|>\tau\}|+
|\{ |\nabla \u -\nabla \u_\varepsilon |>t, |\nabla \u_{\varepsilon
}|\le \tau , |\nabla \u|\le \tau \}|.
\end{multline}
 If $\tau$ is
sufficiently large, then
\begin{equation}\label{419u}
 |\{|\nabla  \u|>\tau \}|< \delta ,
\end{equation}
and, by \eqref{bounde},
\begin{equation}\label{419}
 |\{|\nabla  \u_{\varepsilon }|>\tau \}|=0
%
 \quad \hbox{for $\varepsilon \in (0,1)$.}
\end{equation}
Next, define
$$
\vartheta (t, \tau)=\inf\{[a(| \xi|)\xi -a(|\eta |)\eta ]\cdot
(\xi-\eta):
 |\xi - \eta |\geq t, |\xi |\le \tau , |\eta |\le \tau \},
$$
and observe that $\vartheta(t, \tau) >0$,  by Lemma \ref{approxbis}.
 Thus, since
\begin{multline*}
\int_\Omega [a (|\nabla \u |)\nabla \u  - a (|\nabla
\u_\varepsilon|)\nabla \u_\varepsilon ]\cdot (\nabla \u -\nabla
\u_\varepsilon )\, dx
\\ \geq \vartheta (t , \tau )|\{ |\nabla \u
-\nabla \u_\varepsilon |>t, |\nabla \u_{\varepsilon }|\le \tau ,
|\nabla \u |\le \tau \}|,
\end{multline*}
 by \eqref{eps15}
 $$|\{ |\nabla \u
-\nabla \u_\varepsilon |>t, |\nabla \u_{\varepsilon }|\le \tau ,
|\nabla \u|\le \tau \}|< \delta $$ if $\varepsilon$ is sufficiently
small. Consequently, by \eqref{417}, \eqref{419}, and \eqref{419u},
$$|\{|\nabla \u
-\nabla \u_\varepsilon |>t \}|< 2 \delta$$
 if $\varepsilon$ is sufficiently
small. This proves \eqref{aeepsmis}. Inequality
\eqref{gradienteqdirichlet} follows from \eqref{aeeps} and
\eqref{estimateeps}.

\qed

\bigskip
\par\noindent
{\bf Proof of Theorem \ref{dirichletconvex}}. The proof is the same
as that of Theorem \ref{dirichletc2}, save that the simplified
versions of inequalities  \eqref{37old} and \eqref{37}  without the term depending on
${\rm tr} \mathcal B$, described in
the last part of Lemma \ref{keystep}, have to be exploited in Steps 2 and 4,
respectively. Also, a sequence of (smooth) convex approximating
domains $\Omega _m$ has to be employed in Step 5. \qed

\bigskip
\par\noindent
{\bf Proof of Theorem \ref{neumannc2}}.   The outline of the proof
is analogous to that of Theorem \ref{dirichletc2} for the Dirichlet
problem. Hereafter, we point out some minor  required variants.
\par\noindent
{\bf Step 1}.
Under assumptions \eqref{domainsmooth}--\eqref{nondeg}, the solution
$\u$ to the Neumann problem \eqref{eqneumann} belongs to $W^{2,2}(\o
, \rN)$. A proof of this fact parallels  that of \cite[Theorem
8.2]{BF} in the case of homogeneous Dirichlet boundary conditions.
One can make use of this piece of information to conclude that
$a(|\nabla \u|)\nabla \u \in W^{1,2}(\Omega , \rN)$,
 and hence derive from an
integration by parts in \eqref{weakneumann} that the boundary
condition $\frac{\partial \u}{\partial \nu}=0$ in fulfilled in the
sense of traces.
\par\noindent
As a consequence, a sequence $\{{\bf u}_k\} \subset C^\infty (\Omega
, \rN)\cap C^2(\overline \Omega , \rN)$ can be constructed such that
\begin{equation}\label{convukneu}
 {\bf u}_k \to {\bf u} \quad \hbox{in $W^{2,2}(\Omega , \rN)$,}
\quad \nabla {\bf  u}_k \to \nabla {\bf u}\quad  \hbox{a.e. in
$\Omega $}, \,\, \hbox{and}\,\, \frac{\partial \u_k}{\partial \nu}=0 \,\, \hbox{on}\,\, \partial \o.
\end{equation}
Such
construction can be accomplished as follows. First, one can
(locally) reduce the problem in some neighborhood of each point $x_0
\in
\partial \Omega$ to the case when $\partial \Omega$ is flat via a
change of variables. In order to preserve the boundary condition
$\frac{\partial \u}{\partial \nu}=0$, the new system of
(curvilinear) orthogonal coordinates $(y_1, \dots , y_n)$ can be
chosen in such a way that the level surfaces $\{y_n = c_n\}$, with
$c_n \in \R$, agree with the level surfaces of the distance function
to $\partial \Omega$, and the curves $\{y_1= c _1, \cdots , y_{n-1}=
c _{n-1}\}$, with $c_1 \in \R$, $\dots$ , $c_{n-1}$, are orthogonal
to these level surfaces. Second,  the function $\u$ can be extended
to a function $\widetilde \u$  beyond the flattened boundary of
$\Omega$ by reflection, so that $\widetilde \u$
 is symmetric with respect to the  boundary. The function $\widetilde \u$ is now defined in a complete
neighborhood $U$ of $x_0$. The fact that $\frac{\partial
\u}{\partial \nu}=0$ on the flattened boundary ensures that such an
extension is twice weakly differentiable, and hence belongs to
$W^{2,2}(U)$. Standard mollification of $\widetilde \u$ by a
symmetric kernel provides an approximation of $\widetilde \u$ in
$W^{2,2}(U)$ by a sequence of smooth functions $\widetilde{\bf u}_k$
which satisfy
 $$\widetilde{\bf u}_k \to {\bf u} \quad \hbox{in $W^{2,2}(U , \rN)$,}
\quad \nabla \widetilde{\bf  u}_k \to \nabla {\bf u}\quad \hbox{a.e.
in $U $,}$$ and are symmetric about the boundary of $\Omega$. The
latter property ensures that  $\frac{\partial
\widetilde\u_k}{\partial \nu}=0$ on the boundary of $\o$. The
function $\u _k$ is then just defined as the restriction of
$\widetilde \u _k$ to $\Omega$.
\par\noindent
Via the same argument as in the case of the Dirichlet problem, on
can prove that the sequence $\u _k$ just obtained also fulfills
\eqref{convtr} and \eqref{operator}.
\par\noindent
{\bf Step 2}. Here one shows that inequality \eqref{37neumann} is
fulfilled when $\v$ equals the solution $\u$  to \eqref{eqneumann}.
This follows on applying \eqref{37neumann} of Lemma
\ref{keystepneumann} with $\v=\u _k$ (defined in Step 1), and
passing to the limit as $k \to \infty$  via the same argument as in
the Dirichlet case.
\par\noindent
{\bf Step 3}. This step is exactly the same as in the Dirichlet
case, save that $|\mathcal B|$ replaces ${\rm tr} \mathcal B$
everywhere.
\par\noindent
{\bf Step 4}. Here one applies inequality \eqref{37oldneumann} with
$\v=\u _k$, and obtains the same inequality for $\u$ on passing to
the limit as $k \to \infty$ as in the case of the Dirichlet problem.
\par\noindent
{\bf Step 5}. We construct the sequence of domains $\Omega _m$ as in
the case of Dirichlet problems, and obtain a corresponding sequence
$\{\u_m\}$ of solutions to the Neumann problems in $\Omega _m$
satisfying \eqref{w22m}, \eqref{w1m}, and  \eqref{weakdirm} for
every
  function $\phi \in {\rm Lip}(\rn , \rN)$.
Thanks to \eqref{w1m}, passing to the limit as $m \to \infty$ yields
\eqref{weakdir} for every $\phi \in {\rm Lip}(\rn , \rN)$. Since
$\Omega$ has a Lipschitz boundary, the space of the restrictions to
$\Omega$ of the functions from ${\rm Lip}(\rn , \rN)$ is dense in
$W^{1,B}(\o , \rN)$. Hence, \eqref{weakdir} also holds for every
$\phi \in W^{1,B}(\o , \rN)$.
\par\noindent
{\bf Step 6}. This step is the same as in the case of the Dirichlet
problem. \qed

\bigskip
\par\noindent
{\bf Proof of Theorem \ref{neumannconvex}}. The proof consists in a
slight modification of that of Theorem \ref{neumannc2}.
Specifically,  the  versions of inequalities \eqref{37oldneumann}
and \eqref{37neumann} where the term depending on $|\mathcal B|$ is
dropped, described in the last part of Lemma \ref{keystepneumann},
play  a role in Steps 2 and 4, respectively. Moreover, the
approximating domains $\Omega _m$  in Step 5 have to be chosen
convex.  \qed

\bigskip
\bigskip
\par\noindent
{\bf Acknowledgements}. The authors wish to thank Dominic Breit for
pointing out the interpretation of Theorems \ref{dirichletconvex}
and \ref{neumannconvex} stated in Corollaries \ref{hilberhaar} and
\ref{hilberhaarneumann}, respectively.
\par\noindent
This research was partially supported by the PRIN research project
``Geometric aspects of partial differential equations and related
topics" (2008) of MIUR (Italian Ministry of University and
Research).

%
 {}{}{}{}{}{}
 {}{}{}{}{}{}{}{}{}{}{}{}{}{}{}{}{}{}{}{}{}{}{}{}{}{}{}{}{}{}{}{}{}{}{}{}{}{}
 {}{}{}{}{}{}
 {}{}{}{}{}{}{}{}{}{}{}{}{}{}{}{}{}{}{}{}{}{}{}{}{}{}{}{}{}{}{}{}{}{}{}{}{}{}


\medskip
\begin{thebibliography}{99}
\bibitem[AF]{AF}   E.Acerbi \& N.Fusco,
Regularity for minimizers of nonquadratic functionals: the case
$1<p<2$, \emph{J. Math. Anal. Appl.} {\bf 140} (1989), 115--135.

\bibitem[BC]{BeiraoCrispo1} H.Beir\~{a}o da Veiga \& F.Crispo, On the global $W^{2,q}$ regularity for
nonlinear $N$-systems of the $p$-Laplacian type in $n$ space
variables, \emph{Nonlinear Anal.} {\bf 75} (2012), 4346--4354.
%



\bibitem[BBGGPV]{BBGGPV} P.B\'enilan, L.Boccardo, T.Gallou\"et, R.Gariepy, M.Pierre \& J.L.Vazquez,
An $L^1$-theory of existence and uniqueness of solutions of
nonlinear elliptic equations, \emph{Ann. Sc. Norm. Sup. Pisa} {\bf
22} (1995), 241--273.
\bibitem[BF]{BF} A.Bensoussan \& J.Frehse,  ``Regularity results for nonlinear elliptic systems and applications", Springer-Verlag, Berlin, 2002.
\bibitem[BS]{BS} C.Bennett \&  R.Sharpley, ``Interpolation of operators",
Academic Press, Boston, 1988.
\bibitem[BSV]{BreitSV} D.Breit, B.Stroffolini \& A.Verde,  A general regularity theorem for functionals with $\varphi$-growth,
\emph{J. Math. Anal. Appl.}, {\bf 383} (2011), 226--233.
\bibitem[BZ]{BZ} J.E.Brothers and  W.P.Ziemer, Minimal rearrangements of Sobolev
functions, \emph{J. Reine Angew. Math} {\bf 384} (1988), 153--179.

\bibitem[CDiB]{CDiB} Y.Z.Chen \& E.Di Benedetto, Boundary estimates
for solutions of nonlinear degenerate parabolic systems,  \emph{J.
Reine Angew. Math.} {\bf 395} (1989), 102--131

\bibitem[Ci2]{Cmaxim} A.Cianchi, Maximizing the $L^\infty$ norm of the gradient of
solutions to the Poisson equation, \emph{J. Geom. Anal.} {\bf 2}
(1992), 499--515.
\bibitem[Ci3]{Csharp} A.Cianchi, A sharp embedding theorem for Orlicz-Sobolev spaces,
 \emph{Indiana Univ. Math. J.} {\bf 45} (1996), 39--65.
 \bibitem[Ci4]{Cbound} A.Cianchi, Boundednees of solutions to
variational problems under general growth conditions, \emph{Comm.
Part. Diff. Equat.} \textbf{22} (1997), 1629--1646.
\bibitem[CEG]{CEG} A.Cianchi, D.E.Edmunds \& P.Gurka, On weighted
Poincar\'e inequalities, \emph{Math. Nachr.} {\bf 180} (1996),
 15--41.
 \bibitem[CM1]{CMlipschitz} A.Cianchi \& V.Maz'ya,
 Global Lipschitz regularity for a class of quasilinear elliptic
equations, \emph{Comm. Part. Diff. Equat.} {\bf 36} (2011),
100--133.
\bibitem[CM2]{CMJEMS} A.Cianchi  \& V.Maz'ya, Gradient regularity via rearrangements for $p$-Laplacian type elliptic problems,
\emph{J. Europ. Math. Soc.}, to appear.
\bibitem[CP]{CP} A.Cianchi \& L.Pick, Sobolev embeddings into BMO, VMO and $L^\infty$, \emph{Ark. Math.} {\bf 36} (1998),
317--340.
\bibitem[Di]{Di} E.Di Benedetto, $C^{1+\alpha }$ local regularity of weak solutions of degenerate
elliptic equations, \emph{Nonlinear Anal.}  {\bf 7} (1983),
827--850.
\bibitem[DT]{DT} D.T.Donaldson \& N.S.Trudinger, Orlicz-Sobolev spaces and embedding theorems,
\emph{J. Funct. Anal.} {\bf 8} (1971), 52--75.

\bibitem[DeG]{DeG} E.De Giorgi, Un esempio di estremali discontinue per un  problema variazionale di tipo ellittico,
 \emph{Boll. Un. Mat. Ital.} {\bf 1} (1968), 135--137.

\bibitem[DSV]{DieningSV} L.Diening, B.Stroffolini \& A.Verde,
Everywhere regularity of functionals with $\phi$-growth,
\emph{Manus. Math.} {\bf 129} (2009), 449--481.


\bibitem[DGK]{DGK} F.Duzaar, J.F.Grotowski \& M.Kronz, Partial and full boundary regularity
for minimizers of functionals with nonquadratic growth, \emph{J.
Convex Anal.} {\bf 11} (2004), 437--476.

\bibitem[DM1]{DM} F.Duzaar \& G.Mingione, Gradient estimates via non-linear potentials, \emph{Amer. J.
Math.}  {\bf 133} (2011), 1093--1149.
\bibitem[DM2]{DMnew} F.Duzaar \& G.Mingione,
Local Lipschitz regularity for degenerate elliptic systems,
 \emph{Ann. Inst. Henri Poincar\'e}  {\bf  27} (2010), 1361--1396.

\bibitem[DM3]{DMcont} F.Duzaar \& G.Mingione, Gradient continuity
estimates,  \emph{ Calc. Var. Part. Diff. Equat.}  {\bf 39} (2010),
379--418.


\bibitem[Ev]{Ev} L.C.Evans, A new proof of local $C^{1,\alpha}$
regularity  for solutions of certain degenerate elliptic P.D.E.,
\emph{J. Diff. Eq.} {\bf 45} (1982), 356--373.


\bibitem[Fo]{Fo} M.Foss,
Global regularity for almost minimizers of nonconvex variational
problems,
 \emph{Ann. Mat. Pura Appl.}  {\bf  187} (2008), 263--321.

\bibitem[FPV]{FPV} M.Foss, A.Passarelli di Napoli \& A.Verde,
Global Lipschitz regularity for almost minimizers of asymptotically
convex variational problems,
 \emph{Ann. Mat. Pura Appl.}  {\bf  189} (2010), 127--162.





\bibitem[FM]{FuchsMingione} M.Fuchs \& G.Mingione, Full $C^{1,\alpha}$-regularity for free and
constrained local minimizers of elliptic variational integrals with
nearly linear growth,  \emph{Manus. Math.} {\bf 102} (2000),
227--250.

\bibitem[Gia]{Gia} M.Giaquinta,  ``Multiple integrals in the calculus of variations and nonlinear elliptic systems", Annals of Mathematical Studies,
 Princeton University Press, Princeton, NJ, 1983.

\bibitem[GiaMo]{GiaMo} M.Giaquinta \& G.Modica,  Almost-everywhere regularity for solutions of nonlinear elliptic systems,
 \emph{Manuscr. Math.} {\bf 28} (1979), 109--158.


\bibitem[Gi]{Giusti} E.Giusti,  ``Direct methods in the calculus of variations", World
Scientific, River Edge, NJ, 2003.



\bibitem[GM]{GM} E.Giusti \& M.Miranda,  Un esempio di soluzioni discontinue per un  problema  di minimo relativo ad
un integrale regolare del calcolo delle variazioni,
 \emph{Boll. Un. Mat. Ital.} {\bf 1} (1968), 219--226.


\bibitem[Gr]{Grisvard} P.Grisvard,  ``Elliptic problems in nonsmooth domains", Pitman, Boston, MA, 1985.


\bibitem[HKW]{HKW} S.Hildebrandt, H.Kaul \& K.-O.Widman, An existence theorem for harmonic mappings of Riemannian manifolds,
{\em Acta Math.} {\bf 138}  (1977), 1--16.





\bibitem[Iv]{Iv} P.-A.Ivert, Regularit\"atsuntersuchungen von L\"osungen elliptischer Systeme von quasilinearen Differentialgleichungen zweiter Ordnung,
{\em Manuscr. Math.} {\bf 30}  (1979), 53--88.

\bibitem[JM]{JM} J.Jost \& M.Meier, Boundary regularity for minima of certain quadratic functionals,
{\em Math. Ann.} {\bf 262}  (1983), 549--561.


\bibitem[KrM1]{KrM} J.Kristensen \& G.Mingione, The singular set of minima of integral functionals,  \emph{Arch. Ration. Mech. Anal.} {\bf 180} (2006),
331--398.

\bibitem[KrM2]{KrM1} J.Kristensen \& G.Mingione, Boundary regularity in variational problems,  \emph{Arch. Ration. Mech. Anal.} {\bf 198} (2010),
369--455.

\bibitem[KuM]{KuM} T.Kuusi \& G.Mingione,
Linear potentials in nonlinear potential theory,
 \emph{Arch. Ration. Mech. Anal.}, to appear.



\bibitem[LU1]{LU1} O.A.Ladyzenskaya \& N.N.Ural'ceva,
Quasilinear elliptic equations and variational problems with many
indepedent variables, \emph{Usp. Mat. Nauk.} {\bf 16} (1961), 19--92
(Russian); English translation: \emph{Russian Math. Surveys} {\bf
16} (1961), 17--91.
\bibitem[LU2]{LU2} O.A.Ladyzenskaya \& N.N.Ural'ceva,  ``Linear and quasilinear elliptic equations", Academic Press, New York, 1968.
\bibitem[Le]{Le} J.L.Lewis, Regularity of the derivatives of solutions
to certain degenerate elliptic equations, \emph{Indiana Univ.
Math. J.} {\bf 32} (1983),  849--858.
\bibitem[Li1]{Li1} G.M.Lieberman, The Dirichlet problem for quasilinear elliptic
equations with continuously differentiable data, \emph{ Comm.
Part. Diff. Eq.} {\bf 11} (1986), 167--229.

\bibitem[Mar]{Mar} P.Marcellini,  Everywhere regularity for a class of elliptic systems without growth conditions, \emph{Annali
 Scuola Norm. Sup. Pisa} {\bf 23} (1996), 1--25.
 \bibitem[MM]{MarcusMizel} M.Marcus \& V.J.Mizel, \emph{Absolute continuity of tracks and mappings of Sobolev
spaces}, Arch. Ration. Mech. Anal. {\bf 45}  (1972),  294--320.


 \bibitem[Ma1]{Macounter} V.G.Maz'ya, Examples of nonregular solutions of quasilinear elliptic equations with analytic
 coefficients, {\em Funkc. Anal. Prilozen.} {\bf 2} (1968),
 53--57 (Russian); English translation:
{\em Funct. Anal. Appl.} {\bf 2} (1968), 230--234.


\bibitem[Ma2]{Magradient} V.G.Maz'ya,  The boundedness of the first derivatives of the
solution of the Dirichlet problem in a region with smooth nonregular
boundary,  \emph{Vestnik Leningrad. Univ.} {\bf 24} (1969), 72--79
(Russian); English translation: \emph{Vestnik Leningrad. Univ.
Math.} {\bf 2} (1975), 59--67.
\bibitem[Ma3]{Ma2} V.G.Maz'ya, On weak solutions of the Dirichlet and Neumann problems,
{\em Trusdy Moskov. Mat. Obsc.} {\bf 20} (1969), 137--172 (Russian);
English translation: \emph{Trans. Moscow Math. Soc.} {\bf 20}
(1969), 135--172.
%
\bibitem[Ma4]{Maconvex} V.G.Maz'ya,  On the boundedness of the first
derivatives for solutions to the Neumann-Laplace problem in a convex
domain, \emph{J. Math. Sci. (N.Y.)} {\bf 159} (2009), 104--112.
\bibitem[Ma5]{Mazlibro} V.G.Maz'ya,  ``Sobolev spaces with applications to elliptic partial differential equations", Springer, Heidelberg, 2011.

\bibitem[MS]{MingSiepe} G.Mingione \& F.Siepe,  Full $C^{1,\alpha}$-regularity for minimizers of integral functionals with LlogL-growth,
\emph{ Z. Anal. Anw.} {\bf 18} (1999), 1083--1100
\bibitem[Mi1]{Mi1} G.Mingione,  The singular set
of solutions to non-differentiable elliptic systems,  \emph{Arch.
Ration. Mech. Anal.} {\bf 166} (2003),  287--301.
\bibitem[Mi2]{Mi2} G.Mingione,  Bounds for the singular set of
solutions to non linear elliptic systems, \emph{Calc. Var. Part.
Diff. Equat.} {\bf 18} (2003), 373--400.

\bibitem[Ne]{Ne} J.Necas, Example of an irregular solution to a nonlinear elliptic system with analytic coefficients and conditions for regularity,
in Theor. Nonlin. Oper., Constr. Aspects. Proc. 4th Int. Summer
School. Akademie-Verlag, Berlin, 1975, 197--206.
\bibitem[SY]{SY} V.Sver\'{a}k \& X.Yan, Non-Lipschitz minimizers of smooth uniformly convex variational integrals,
\emph{Proc. Natl. Acad. Sci. USA} {\bf 99} (2002),   15269--15276.
\bibitem[To]{To} P.Tolksdorf, Regularity for a more general class of quasilinear elliptic equations,
\emph{J. Diff. Equat.} {\bf 51} (1983), 126--150.
\bibitem[Uhl]{Uhl} K.Uhlenbeck, Regularity for a class of non-linear elliptic
systems, \emph{Acta Math.} {\bf 138} (1977), 219--240.
\bibitem[Ur]{Ur} N.N.Ural'ceva, Degenerate quasilinear elliptic systems,
\emph{Zap. Naucn. Sem. Leningrad. Otdel. Mat. Inst. Steklov. (LOMI)}
{\bf 7} (1968), 184--222 (Russian).
\bibitem[Zi]{Z} W.P.Ziemer, \lq \lq Weakly differentiable
functions",  Springer, Berlin, 1989.

\end{thebibliography}
\end{document}